\documentclass[11pt, reqno]{amsart}
\usepackage{amsmath,amssymb,amsfonts,amsthm}

\usepackage{float}
\usepackage[colorlinks=true,linkcolor=blue, allcolors=blue]{hyperref}
\usepackage{tabularx}
\usepackage{adjustbox}
\usepackage[a4paper,margin=1in]{geometry}
\usepackage{tikz}
\usepackage[T1]{fontenc}
\usepackage{lmodern}
\usepackage{amsmath,amssymb}

\usepackage{bbm}
\usepackage{csquotes}
\usepackage{bbm}

\newcommand{\gram}{\mathrm{Gram}}

\newcommand{\cA}{\mathcal{A}}
\newcommand{\cB}{\mathcal{B}}
\newcommand{\cC}{\mathcal{C}}
\newcommand{\cH}{\mathcal{H}}
\newcommand{\cM}{\mathcal{M}}
\newcommand{\cR}{\mathcal{R}}
\newcommand{\cT}{\mathcal{T}}

\newcommand{\bZ}{\mathbb{Z}}

\newcommand{\bu}{\mathbf{u}}
\newcommand{\bp}{\boldsymbol\pi}
\newcommand{\bw}{\mathbf{w}}
\newcommand{\bk}{\mathbbm{k}}

\newcommand{\bs}{\mathbf{s}}
\newcommand{\bt}{\mathbf{t}}

\newcommand{\Hom}{\mathrm{Hom}}

\newcommand{\BFS}{\mathrm{BFS}}

\newcommand{\remove}[1]{{}}

\newcommand{\wrev}[1]{{#1}^{\mathrm{rev}}}
\newcommand{\winv}[1]{{#1}^{\mathrm{inv}}}

\numberwithin{equation}{section}
\newtheorem{theorem}{Theorem}[section]
\newtheorem{lemma}[theorem]{Lemma}
\newtheorem{definition}[theorem]{Definition}

\newtheorem{proposition}[theorem]{Proposition}

\newtheorem{conjecture}[theorem]{Conjecture}

\theoremstyle{remark}
\newtheorem{remark}[theorem]{Remark}
\newtheorem{example}[theorem]{Example}

\makeatletter
\@namedef{subjclassname@2020}{\textup{2020} Mathematics Subject Classification}
\makeatother

\title[The BMM conjecture holds for all complex reflection groups]{The BMM symmetrising trace conjecture holds for all complex reflection groups}

\author{Eirini Chavli}
\address{Institut Denis Poisson,
Universit\'e de Tours,
Parc de Grandmont,
37200 Tours, France.}
\email{eirini.chavli@univ-tours.fr}

\author{Maria Chlouveraki}
\address{Department of Mathematics,
	National and Kapodistrian University of Athens,
	Panepistimioupolis,
	157 84 Athens, Greece.}
\email{mchlouve@math.uoa.gr}

\author{G\"otz Pfeiffer}
\address{School of Mathematical and Statistical Sciences, University of Galway, University Road, Galway, H91 TK33, Ireland.}\email{goetz.pfeiffer@universityofgalway.ie}

\author{Kostas Psaromiligkos}
\address{Department of Mathematics,
	National and Kapodistrian University of Athens,
	Panepistimioupolis,
	157 84 Athens, Greece.}\email{kostaspsa@math.uoa.gr}

\date{}

\begin{document}

	\maketitle

	\begin{abstract}
We prove the symmetrising trace conjecture of Brou\'e, Malle and Michel
for the generic Hecke algebras of all exceptional complex reflection
groups, and we treat the remaining case of the infinite series, so that
the conjecture now holds for all complex reflection groups. A key
ingredient is a new family of bases, constructed via breadth-first
search, which generalise the standard basis of an Iwahori--Hecke
algebra; they provide a uniform proof of the freeness conjecture of
Brou\'e, Malle and Rouquier for these groups, and a model for
computing in the Hecke algebras. The canonical trace is then obtained
from other bases: generalisations of the $z$-bases of earlier work for
the smaller groups, and bases built along towers of parabolic subgroups
for the largest ones. The verification for each group reduces to a
finite computation, carried out using programs developed with the
assistance of a large language model.
\end{abstract}
	
	\maketitle

	\section{Introduction}

Two conjectures lie at the foundation of the theory of Hecke algebras
attached to complex reflection groups. The first, due to Brou\'e,
Malle and Rouquier \cite{BMR}, asserts that the generic Hecke algebra
$\mathcal H(W)$ of a complex reflection group $W$ is a free module over
its ring of definition, of rank $|W|$. The second, due to Brou\'e,
Malle and Michel \cite{BMM}, asserts that $\mathcal H(W)$ carries a
canonical symmetrising trace: a linear form that makes $\mathcal H(W)$
a symmetric algebra, specialises to the canonical trace on the group
algebra of $W$, and satisfies a compatibility condition with the
inversion of the parameters; these conditions determine it uniquely.
In this paper we complete the proof of the second conjecture for all  complex reflection groups.

Both conjectures generalise properties of the real case. A complex
reflection group is a finite group of linear transformations of a
complex vector space generated by pseudo-reflections, that is, by
elements of finite order whose fixed space is a hyperplane. Those that
can be realised over the real numbers are precisely the finite Coxeter
groups, among which are the Weyl groups. For a Weyl group, the
Iwahori--Hecke algebra arises in the representation theory of the
corresponding finite reductive groups, as the endomorphism algebra of
the permutation module on the cosets of a Borel subgroup, and it can
equivalently be defined by deforming the Coxeter presentation of the
group; the latter definition makes sense for every finite Coxeter
group. In this setting, the length function and Matsumoto's theorem
provide a standard basis $(T_w)_{w\in W}$, which shows that the algebra
is free of rank $|W|$, and the linear form $T_w\mapsto\delta_{1w}$ is a
canonical symmetrising trace. Neither tool is available for complex
reflection groups in general, which is why the two statements become
conjectures there.

By the classification of Shephard and Todd \cite{ShTo}, every
irreducible complex reflection group either belongs to the infinite
series $G(de,e,n)$ or is one of $34$ exceptional groups
$G_4,\dots,G_{37}$, of which $G_{23}$, $G_{28}$, $G_{30}$ and
$G_{35},G_{36},G_{37}$ are real. The braid group $B(W)$ and the Hecke
algebra $\mathcal H(W)$ of an arbitrary complex reflection group were
introduced in \cite{BMR}, as part of a programme in which certain
complex reflection groups would serve as Weyl groups of hypothetical
objects, the Spetses, extending the theory of finite reductive groups
\cite{BMM}. Since then, the two conjectures have been assumed in much
of the work on these algebras and on related structures, such as
Cherednik algebras. In particular, the Schur elements of
$\mathcal H(W)$, which govern its modular representation theory, are
defined through the canonical symmetrising trace, and have been explicitly computed for every complex reflection group assuming its existence.

The BMR freeness conjecture is now a theorem, established case by case by
many authors over two decades; we give a detailed account in Section~\ref{sec-2.2}.
The BMM symmetrising trace conjecture has proved considerably harder. For the infinite
series $G(de,e,n)$, a symmetrising trace is known \cite{BreMa, MM98}, and its canonicality follows from
\cite{BMM, GIM} except when $n=2$ and $e$ is even; we treat this
remaining case in Section~\ref{sec-7}.  For the non-real
exceptional groups, the conjecture was first proved by Malle and Michel
\cite{MM10} for $G_4$, $G_{12}$, $G_{22}$ and $G_{24}$, and then, in a
series of papers \cite{BCCK, BCC, ChPf}, for all the groups
$G_4,\dots,G_{15}$, which are the exceptional groups of rank $2$ in the
tetrahedral and octahedral families. This left open the groups of the
icosahedral family $G_{16},\dots,G_{21}$ and the groups of higher rank
$G_{25}$, $G_{26}$, $G_{27}$, $G_{29}$, $G_{31}$, $G_{32}$, $G_{33}$
and $G_{34}$.

There are two main reasons why these groups resisted.
The first is
their size: the generic Hecke algebra of $G_{19}$ has rank $3600$ and
depends on ten parameters, and that of $G_{34}$ has rank $39191040$.
The second is more fundamental. All previous proofs follow the same
strategy: one chooses a basis $\mathcal B$ of $\mathcal H(W)$
containing $1$ and consisting of lifts of the elements of $W$ to the
braid group, defines $\tau$ as the linear form that sends an element to
its coefficient on $1$, and shows that the Gram matrix
$(\tau(bb'))_{b,b'\in\mathcal B}$ is symmetric with invertible
determinant. Whether this works depends entirely on the basis. 
In \cite{ChPf}, suitable bases were found among the so-called $z$-bases,
whose elements are products of powers of a central element
$\mathbf z$, powers of a generator, and coset representatives, and
whose Gram matrices have a block structure that makes their
determinants computable. For larger groups such bases become rare: for
$G_{16}$, we show that none of the $405504$ candidates of this form
that we consider is a basis. We therefore need more candidate bases and a way to identify the good ones without computing in the generic algebra, where calculations with multivariate polynomials are computationally expensive or even infeasible.

In this article, in order to make the computations feasible, we carry them out not over the ring of
definition $R(W)$ of $\mathcal H(W)$, but after specialising to $1$ the
constant terms of the Hecke relations, which reduces the number of
parameters. The results are then lifted to the generic algebra by means
of the grading of $\mathcal H(W)$ introduced in \cite{BMM}: we deduce
that bases and symmetrising traces can be detected after
specialisation, provided that all the objects involved are homogeneous.
Before any exact computation, we further discard most
candidates by working in specialisations of $\mathcal H(W)$ over finite fields  and using two \emph{modular tests}, which give necessary conditions for a set to be a basis and for its coefficient on 1 to be a symmetrising trace.
In this way, we prove  Theorem \ref{thm:Gj} and Theorem \ref{thm:top}, which, together with the classical case of the real groups, yield the following.

\begin{theorem}\label{thm:main}
	The BMM symmetrising trace conjecture holds for all exceptional complex
	reflection groups.
\end{theorem}

Combined with Section~\ref{sec-7}, where we treat the last open case of the
infinite series, this yields the following theorem as corollary.

\begin{theorem}
The BMM symmetrising trace conjecture holds for all complex reflection
groups.
\end{theorem}

Our starting point for proving Theorem \ref{thm:main} is a new family of bases. For a complex reflection
group there is no length function, and hence no standard basis; but the
standard basis of an Iwahori--Hecke algebra can be described without
any reference to length, as the set of words read off a breadth-first
search through the group, and this description makes sense in general.
Given a parabolic subgroup $P$ of $W$ and an ordering $o$ of a
generating set $A$, we obtain in this way a set $\mathrm{BFS}_o$ of $|W|$
elements of $\mathcal H(W)$, which for real reflection groups is the
standard basis whatever the ordering. To decide whether
$\mathrm{BFS}_o$ is a basis, we complete its coset table, which records
the products of its elements with the generators, by a fully automated procedure (Section~\ref{coset}). We obtain:

\begin{theorem}\label{thm:bfs}
	For every exceptional complex reflection group $W$, there exist a generating set $A$ of $W$ and an ordering $o$ of $A$ such that $\mathrm{BFS}_o$ is an $R(W)$-basis of $\mathcal H(W)$.
\end{theorem}

This incidentally gives a new proof of the BMR freeness conjecture for the
exceptional groups, by a single construction and a single verification
procedure. Earlier proofs treated the groups one at a time, or
uniformly only within a family, such as the groups of rank $2$ in
\cite{Tsu} or the $2$-reflection groups in \cite{MaPf}. Moreover, a
completed coset table is a faithful matrix model of $\mathcal H(W)$: it
allows us to express any element in the basis $\mathrm{BFS}_o$, and all
subsequent computations rely on it.

The $\mathrm{BFS}_o$ bases themselves rarely define the canonical
trace: for $G_{18}$ and $G_{19}$, for example, none of those we found
does. To find bases that do, we use the modular tests to select candidates among larger families, and the proof of Theorem~\ref{thm:main} takes two different forms, according to the order of the group.

For the groups $G_{4},\dots,G_{29}$, we generalise $z$-bases to
\emph{$z^k$-bases}, in which only the powers of $\mathbf z^k$ appear,
for a divisor $k$ of the order of the centre. Candidates correspond to
spanning trees of a coset graph, and can be enumerated. For $G_{16}$,
$G_{17}$, $G_{18}$ and $G_{19}$ we find $z^k$-bases with $k=2,4,6,12$
respectively, and for the other groups we find $z$-bases.

For the groups $G_{31}$, $G_{32}$, $G_{33}$ and $G_{34}$, whose orders
make this approach impractical, we use instead a tower of parabolic
subgroups, starting from a cyclic one, with bases obtained by
breadth-first search at each level. We give a criterion that reduces
the existence of a symmetrising trace at each level to the
corresponding statement one level below, together with the
invertibility of a matrix with entries in the parabolic Hecke
subalgebra, indexed by the cosets of the parabolic subgroup and
therefore much smaller than the Gram matrix.  At the top level of
$G_{32}$, no breadth-first basis defines a trace function, and we
prescribe the words representing two of the cosets before completing
the tree by breadth-first search. This adapts to coset trees an idea
first used for $G_7$ in \cite[\S4.2.4]{BCCK}, where a basis is modified
by replacing a few of its elements. As for the smaller groups, all
verifications are carried out after specialisation.

All computations are carried out by programs available on the
project's webpage \cite{web}. These programs are an integral part of
the proofs. They also provide, for each exceptional group, an efficient
way to compute the canonical symmetrising trace of any element of
$\mathcal H(W)$, which we hope will be useful in applications. 

The paper is organised as follows. Section~\ref{prel} recalls complex reflection
groups, their braid groups and Hecke algebras, and the conjectures and
their parabolic versions. Section~\ref{sec:reduce} shows that it suffices to compute
after specialising some of the parameters, using a grading of
$\mathcal H(W)$; it introduces $\mathrm{pr}$-admissible collections,
which yield the conjecture, and the modular tests. Section~\ref{standard} is devoted
to the $\mathrm{BFS}_o$ bases and to coset tables, and contains the
proof of Theorem~\ref{thm:bfs}. In Sections~\ref{bases} and~\ref{rem} we construct
$\mathrm{pr}$-admissible collections, for $G_4,\dots,G_{29}$ using
$z^k$-bases, and for $G_{31},\dots,G_{34}$ using towers of parabolic
subgroups, which completes the proof of Theorem~\ref{thm:main}.  
In Section~\ref{sec-7}, we prove the conjecture for the groups $G(de,e,2)$ with
$e$ even, which completes the proof for the infinite series.
In the Appendix, we give a table with the exact values of the Gram matrix determinants.

\subsection*{Disclosure}
This project began in October 2025, with the aim of completing the proof of the BMM symmetrising trace conjecture for the rank two groups. The previous papers~\cite{MM10, BCCK, BCC, ChPf1, ChPf} had dealt with the groups $G_4,\ldots,G_{15}$ and $G_{22}$, leaving the groups $G_{16},\ldots,G_{21}$ of the icosahedral family. Our initial goal was to extend the methods of~\cite{ChPf} to these remaining groups, which had remained out of reach because of their size. A different approach was therefore needed. We began by searching for suitable bases, but the higher orders of the generators gave rise to a much larger number of candidates.

Over the following year, we carried out extensive mathematical and computational experimentation. An important role of a large language model (LLM) was to make the exploration of candidates much faster. Without it, the time required to write the necessary programs would have forced us to make guesses about which candidates to investigate first. Instead, by generating programs rapidly from our hand-written programs for smaller groups and our descriptions of the algorithms, using the LLM allowed us to test a wide range of candidates in parallel. The mathematical strategy, including methods for making the computations efficient, such as modular computations and changes of rings, was developed by us. The programs also include an implementation of polynomial arithmetic following an idea of M\"uller~\cite{Mul}. The figures and tables were generated with similar assistance.

Moreover, two useful mathematical ideas emerged from our examination of the LLM's output:\\
1. Our use of the grading of~\cite{BMM} in order to lift results to the generic algebra (Section~\ref{sec:reduce}). \\
2. Our pairing sets for $G_{32}$: The main tool proving the BMM conjecture is the pairing method, formulated in Lemma~\ref{pair} and Proposition~\ref{symcond2}. Its application requires the choice of suitable sets to be paired. For $G_{32}$, after our initial attempts to find such sets had failed, we recognised in the LLM's output an application of an idea developed by the first two authors in~\cite{BCCK}, which led us to the required sets (Section~\ref{subsec:Y}).

\subsection*{Acknowledgements} This research is supported by the Onassis Foundation. We thank Panagiotis Mertikopoulos, who suggested early in this project that we use large language models in our work, and Zafeirakis Zafeirakopoulos for discussions on computing determinants of large matrices, which inspired the modular tests of Section~\ref{modtests}.

	\section{Preliminaries}\label{prel}
	
	\subsection{Hecke algebras of complex reflection groups}\label{sec-2.1}
	Let $V$ be a finite dimensional complex vector space.
	A {\em complex reflection group} $W$ is a finite subgroup of ${\rm GL}(V)$ that is generated by {\em pseudo-reflections}, which are non-trivial elements fixing a hyperplane in $V$, called the {\em reflecting hyperplane} of the given pseudo-reflection.
	We can and shall always assume that $V$ has no non-zero pointwise $W$-fixed points.
	
	If $W$ acts irreducibly on $V$, then $W$ is an \emph{irreducible complex reflection group} and it either belongs to the infinite series $G(de,e,n)$ or is one of the $34$ exceptional groups $G_4,\ldots,G_{37}$ in the Shephard--Todd classification \cite{ShTo}. In this case, the dimension of $V$ is called the \emph{rank} of $W$. Every complex reflection group is isomorphic to a product of irreducible complex reflection groups.
	
	Let  $\cA$ denote the set of reflecting hyperplanes of $W$. Set $\mathcal{M}:=V\setminus \cup_{H\in \mathcal{A}}H$. We denote by $p:\cM\twoheadrightarrow \cM/W$ the canonical surjection. 	The following definition needs a choice of basepoint $x_0\in \cM$ but gives isomorphic groups independently of the choice \cite{BMR}:
	
	\begin{definition}{\rm
			The  \emph{pure braid group} of $W$ is $P(W)=\pi_1(\cM,x_0)$ and the \emph{braid group} of $W$ is $B(W)=\pi_1(\cM/W,p(x_0))$.}
	\end{definition}

	Moreover, Brou\'e, Malle and Rouquier associated to every pseudo-reflection in $W$ homotopy classes in $B(W)$ that are called \emph{braided reflections}. By  \cite[Theorem 0.1]{Bes2}, the braid group $B(W)$ admits a presentation with the generators being braided reflections and the relations being homogeneous relations between positive words in the generating elements, that is, braid relations. 
	A Coxeter-like presentation of $W$ is derived from this Artin-like presentation of $B(W)$, thanks to the following short exact sequence:
	$$1\rightarrow P(W) \rightarrow B(W) \rightarrow W \rightarrow 1.$$ 
	We denote by $S$ the set of generators in that presentation and we call them the \emph{standard generators} of $W$.

	From now on, we will denote by $Z(G)$ the centre of any group $G$.
	If $W$ is an irreducible complex reflection group, then $Z(W)$ is isomorphic to a finite subgroup of $\mathbb{C}^{\times}$, and it is thus a cyclic group. 
	We denote by $\boldsymbol{\pi}$ and $\boldsymbol{\rm{z}}$ the homotopy classes of the loops $t\mapsto x_0 \,{\rm exp}(2\pi i t)$ and $t\mapsto  x_0\,{\rm exp}(2\pi i t/|Z(W)|)$ respectively. In particular, $\boldsymbol\pi=\mathbf z^{|Z(W)|}$. Brou\'e, Malle and Rouquier \cite[Lemma 2.22 (2)]{BMR} proved that  $\boldsymbol{\pi} \in Z(P(W))$ and $\boldsymbol{{\rm z}} \in Z(B(W))$, and they conjectured that each of the two elements  generates the centre it belongs to as a cyclic group. The fact that $Z(P(W)) = \langle \boldsymbol{\pi} \rangle$ was proved by Digne, Marin and Michel \cite[Theorem 1.2]{DMM}, while it was Bessis \cite[Theorem 12.8]{Bes3} who proved that
	$Z(B(W)) = \langle \boldsymbol{{\rm z}} \rangle$.

	For each $H\in \cA$, the subgroup $W_H\subseteq W$ stabilizing $H$ is cyclic. Let $e_H$ denote the order of $W_H$. The unique generator  $s_H$  of $W_H$ that has determinant $e^{2\pi i/e_H}$ is called the {\em distinguished} pseudo-reflection of $H$.  It is immediate that the set $\cR=\{s_H \mid H\in \cA\}$ generates $W$.
	
	Now, for each distinguished pseudo-reflection $s\in \mathcal{R}$, we denote its order by $e_s$ and  we choose a
	set of $e_s$ indeterminates $u_{s,1},\dots, u_{s,e_s}$, such that
	$u_{s,j}=u_{t,j}$ if $s, t \in \mathcal{R}$ are conjugate in $W$. 
	The following definition \cite{BMR} generalises the notion of Iwahori--Hecke algebras from real to complex reflection groups.
	
	\begin{definition}{\rm
			The \emph{generic Hecke algebra} $\cH(W)$ associated to $W$ is the
			quotient of the group algebra $\mathbb{Z}[u_{s,j},u_{s,j}^{-1}][B(W)]$  by the ideal $I$
			generated by the elements of the form
			\begin{equation}
				(\mathbf{s}-u_{s,1})(\mathbf{s}-u_{s,2})\cdots (\mathbf{s}-u_{s,e_s}),
				\label{Hecker}
			\end{equation}
			where $s$ runs over a set of representatives of the $W$-conjugacy classes of $\mathcal{R}$ 
			and $\mathbf{s}$ is a braided reflection associated to $s$.}
	\end{definition}

	An equivalent definition of $\cH(W)$ follows by expanding the relations \eqref{Hecker}.
	More precisely, we have
	\begin{equation}
		{\mathbf{s} }^{e_{s}}=a_{{s},e_{s}-1}{\mathbf{s}}^{e_{s}-1}+a_{{s},e_{s}-2}{\mathbf{s}}^{e_{s}-2}+\dots+a_{{s},0},
		\label{Hecker2}
	\end{equation}
	where $a_{{s},e_{s}-k}:=(-1)^{k-1}f_k({u}_{{s},1},\dots,{u}_{{s},e_{s}})$ with $f_k$ denoting the $k$-th elementary symmetric polynomial, for $k=1,\ldots,e_{s}$. We call the above relations the \emph{positive Hecke relations}. 
	Since $a_{s,0}$ is invertible in $\mathbb{Z}[u_{s,j},u_{s,j}^{-1}]$, we obtain that ${\mathbf{s}}$ is invertible in $\cH(W)$ with
	\begin{equation}\label{invhecke}
		\mathbf{s}^{-1}=a_{s,0}^{-1}\,\mathbf{s}^{e_{s}-1}-a_{s,0}^{-1}\,a_{{s},e_{s}-1}\mathbf{s}^{e_{s}-2}-a_{s,0}^{-1}\,a_{{s},e_{s}-2}{\bf s}^{e_{s}-3}-\dots- a_{s,0}^{-1}\,a_{{s},1}.
	\end{equation}
	We call  relations \eqref{invhecke} the \emph{inverse Hecke relations}.
	Note that, from now on, the elements of $B(W)$ and their images inside $\cH(W)$
	will be denoted by the same (bold) symbols.  In particular, for $s\in S$ we also
	write $\mathbf s$ for the corresponding generator of $\cH(W)$, and we set
	$\mathbf S=\{\mathbf s\mid s\in S\}$. By the inverse Hecke relations,
	$\mathbf s^{-1}$ is a polynomial in $\mathbf s$, so $\mathbf S$ generates
	$\cH(W)$ as an  algebra.

	In this paper we will consider $\cH(W)$  over the ring $R(W)=\mathbb Z\big[a_{s,0}^{\pm1}, a_{s,1},\dots,a_{s,e_s-1}\big]$. Since the $a_{s,i}$ are, up to sign, the elementary
	symmetric polynomials in the $u_{s,j}$, they are algebraically
	independent, and $R(W)$ is a subring of
	$\mathbb Z[u_{s,j},u_{s,j}^{-1}]$. The relations \eqref{Hecker2}, together with
	the braid relations, define an $R(W)$-algebra, which becomes $\cH(W)$
	after extending scalars to $\mathbb Z[u_{s,j},u_{s,j}^{-1}]$; all our
	results transfer along this base change (see also
	\cite[Proposition 2.3]{MarNew}).
	Moreover, often we need to work over different rings, where the parameters are specialised or over different characteristic. We then define a \emph{specialisation}, that is, a ring homomorphism $\theta: R(W) \rightarrow R'$ where $R'$ is the new ring we want to work on, and we denote by $\cH_\theta$ the specialised Hecke algebra $\cH(W)\otimes_{\theta} R'$.
	For example, if we take $\theta: R(W) \rightarrow \mathbb{C}$ to be the map $a_{s,i} \mapsto \delta_{0i}$,
	then $\cH_\theta \cong \mathbb{C}[W]$; this is why we say that the Hecke algebra is a deformation of the group algebra of $W$.

	\subsection{The two fundamental conjectures on generic Hecke algebras}\label{sec-2.2}
	
	Iwahori--Hecke algebras have a \emph{standard basis} denoted by $(T_w)_{w \in W}$, indexed by the elements of the finite Coxeter group $W$ (see \cite[IV, \S 2]{Bou05}). 
	Namely, if $w=s_{i_1}\ldots s_{i_r}$ is a reduced expression for the element $w$ as product of the simple reflections generating $W$ (that is, $r$ is minimal), then $T_w:={\mathbf{s}_{i_1}}\ldots {\mathbf{s}_{i_r}}$.
	Therefore, in the real case, $\cH(W)$ is a free $R(W)$-module of rank $|W|$. 
		Brou\'e, Malle and Rouquier conjectured such a result for non-real complex reflection groups \cite[\S 4]{BMR}:
	
	\begin{conjecture}[BMR freeness conjecture]
		The algebra $\cH(W)$ is a free $R(W)$-module of rank $|W|$. 
	\end{conjecture}
	
	The BMR freeness conjecture became a theorem 20 years after it was stated. 
	The last cases for which it was proved were $G_{17},\,G_{18},\,G_{19}$ by Tsuchioka \cite{Tsu} (with a computer method applicable to all rank $2$ groups, that is, $G_4,\ldots,G_{22}$), where the bases are described implicitly as the irreducible words of a
	machine-generated rewriting system.
		In all other non-real cases, explicit bases were obtained as follows: for
	\begin{itemize}
		\item the groups of the infinite series $G(de,e,n)$ by \cite{ArKo, BM, Ar};\smallbreak
		\item the group $G_4$ by \cite{BM, Fun, Mar41, Ch18} (4 independent proofs); \smallbreak
		\item the groups $G_8,\,G_{16}$ by \cite{Ch18}; \smallbreak
		\item the groups $G_5,\ldots, G_{15}$ by \cite{Ch17}; \smallbreak
		\item the groups $G_{20}, \,G_{21}$ by \cite{MarNew};\smallbreak
		\item the groups $G_{25}, \,G_{32}$ by \cite{Mar41};\smallbreak
		\item the group $G_{26}$ by  \cite{Mar43};\smallbreak
		\item the groups $G_{12}, G_{22},G_{24},G_{27}, G_{29}, G_{31}, G_{33}, G_{34}$ by \cite{MaPf}. 	\end{itemize}

	The following criterion is frequently used in proofs of the BMR freeness conjecture. 
		For a proof, see \cite[proof of Theorem 4.24]{BMR} or \cite[Proposition 2.4(1)]{Mar43}. 
		\begin{proposition}\label{span}
			Let $X$ be a set of $|W|$ elements of $\cH(W)$. If $X$ spans $\cH(W)$ as an $R(W)$-module, then $X$ is an $R(W)$-basis of $\cH(W)$. In particular, $\cH(W)$ is a free $R(W)$-module of rank $|W|$.
		\end{proposition}

	Iwahori--Hecke algebras associated to real reflection groups are also examples of \emph{symmetric algebras}, namely they admit a linear map $\tau:  \cH(W) \rightarrow R(W)$ such that
	$\tau(hh') = \tau(h'h)$ for all $h, h' \in \cH(W)$ (\emph{i.e.,} $\tau$ is a \emph{trace function}) and the bilinear form
	$(h, h') \mapsto \tau(hh')$ is non-degenerate. Equivalently, the Gram matrix $\gram_{\tau}(\cB,\cB):=\bigl(\tau(bb')\bigr)_{b,b'\in\cB}$ with respect to some, and hence every, $R(W)$-basis $\mathcal{B}$ of $\mathcal H(W)$ is symmetric and has invertible determinant in $R(W)$. Such a linear map is called \emph{symmetrising trace}.
	A symmetrising trace for Iwahori--Hecke algebras is given by $\tau(T_w)=\delta_{1w}$  \cite[IV, \S 2]{Bou05} and it is called the \emph{canonical symmetrising trace} on $\cH(W)$, because it corresponds to the canonical symmetrising trace on the group algebra. Brou\'e, Malle and Michel  conjectured the existence of a canonical symmetrising trace also for non-real complex reflection groups \cite[\S2.1, Assumption 2(1)]{BMM}:
	
	\begin{conjecture}[BMM symmetrising trace conjecture]\label{bmmconj}
		The algebra $\cH(W)$ admits a linear form $\tau: \cH(W)\rightarrow R(W)$ such that 
		\begin{enumerate}
			\item $\tau$ is a symmetrising trace, 
			\item $\tau$ specialises to the canonical symmetrising trace of the group algebra $\bZ[W]$ upon specialising $a_{s,0}\mapsto 1$ and $a_{s,i}\mapsto 0$ for all $i\neq 0$,
			\item for all ${\bf b} \in B(W)$, $\tau$ satisfies
			$$\tau({{\bf b}^{-1}})^{\vee}=\frac{\tau({{\bf b}\bp})}{\tau({\bp})},$$
			where $\vee$ denotes the automorphism of  $R(W)$ given by $a_{s,0} \mapsto a_{s,0}^{-1}$ and $a_{s,i} \mapsto -a_{s,0}^{-1}a_{s,e_s-i}$ for $1\le i\le e_s-1$.
		\end{enumerate}
	\end{conjecture}

	Note that, by  \cite[2.1]{BMM}, since the BMR freeness conjecture holds,  there exists at most
	one symmetrising trace satisfying Conditions (2) and (3) of Conjecture \ref{bmmconj}, meaning that $\tau$ is unique. We call $\tau$ the \emph{canonical symmetrising trace on  $\mathcal{H}(W)$}.

	The BMM symmetrising trace conjecture has been proved for the following non-real exceptional complex reflection groups:
	\begin{itemize}
		\item the group $G_4$ by \cite{MM10, Mar46, BCCK} (3 independent proofs); \smallbreak
		\item the groups $G_5, G_6, G_7, G_8$ by \cite{BCCK}; \smallbreak
		\item the group $G_{13}$ by \cite{BCC}; \smallbreak
		\item the groups $G_9, G_{10}, G_{11}, G_{14}, G_{15}$ by \cite{ChPf}; \smallbreak
		\item the groups $G_{12}, G_{22}, G_{24}$ by \cite{MM10}.
	\end{itemize}	
	A symmetrising trace satisfying conditions (1) and (2) has been constructed for  $G(de,e,n)$ by \cite{BreMa, MM98}. Condition (3) for that trace derives from \cite[Lemma 2.7]{BMM} and  \cite[Theorem 5.2]{GIM} for all cases except for when $n=2$ and $e$ is even.

	A basis $\cB$ for $\cH(W)$ naturally defines linear forms $\lambda^{\cB}_b:\cH(W)\rightarrow R(W)$ for all $b \in \cB$ such that,  for all $h\in \cH(W)$,
	$$h=\sum_{b\in \cB}\lambda^{\cB}_b(h) b.$$	
	The map $b\mapsto \lambda_b^{\cB}$ linearly extends to an isomorphism of free $R(W)$-modules from $\cH(W)$ to the space of linear forms $\cH(W)^{*}=\Hom_{R(W)}(\cH(W),R(W))$. The set $\{\lambda_b^\cB\mid b\in \cB\}$ is an $R(W)$-basis of $\cH(W)^{*}$, and any linear form $E:\cH(W)\rightarrow R(W)$ can be easily written out in this basis as
	$$E=\sum_{b\in \cB} E(b)\lambda^{\cB}_b$$ 
	Thus, the canonical symmetrising trace 
	$\tau$, if it exists, must also be an $R(W)$-linear combination of $\lambda^{\cB}_b$.
	Because of Condition (2) of Conjecture \ref{bmmconj}, we are particularly interested in constructing bases that contain $1$ and for which $\tau=\lambda^{\cB}_1$.

	A \emph{section} of the natural map $B(W)\twoheadrightarrow W$ is a map
$W\to B(W)$, $w\mapsto\mathbf w$, such that $\mathbf w$ maps to $w$ for
every $w\in W$ and $\mathbf 1=1$. If the elements $\mathbf w$, for
$w\in W$, form an $R(W)$-basis $\mathcal B$ of $\mathcal H(W)$, we call
$\mathcal B$ a \emph{section basis}. For a section basis $\mathcal B$,
the linear form $\lambda^{\mathcal B}_1$ specialises to the canonical
trace of the group algebra under $a_{s,0}\mapsto1$, $a_{s,i}\mapsto0$,
since $\mathcal B$ specialises to the basis $W$; thus
$\lambda^{\mathcal B}_1$ satisfies Condition~(2) of
Conjecture~\ref{bmmconj}. The following is
\cite[Conjecture~2.6]{MM10}.

\begin{conjecture}[The lifting conjecture]
There exists a section $w\mapsto\mathbf w$ of
$B(W)\twoheadrightarrow W$ such that $\tau(\mathbf w)=\delta_{1w}$ for
all $w\in W$, where $\tau$ is the canonical symmetrising trace.
\end{conjecture}

If $\mathcal B$ is a section basis and $\lambda^{\mathcal B}_1$ is a
symmetrising trace satisfying Condition~(2), then, by
\cite[Proposition~2.7]{MM10}, Condition~(3) of
Conjecture~\ref{bmmconj} for $\lambda^{\mathcal B}_1$ is equivalent to
\begin{equation}\label{extra con2}
\lambda^{\mathcal B}_1(\mathbf w^{-1}\boldsymbol\pi)=0
\quad\text{for all } w\in W\setminus\{1\}.
\end{equation}
In that case, $\lambda^{\mathcal B}_1$ is the canonical symmetrising
trace, and the section of $\mathcal B$ satisfies the lifting conjecture.

	\subsection{The two fundamental conjectures on parabolic Hecke subalgebras}	\label{sec-2.3}
	A subgroup $P$ of a complex reflection group $W \subset {\rm GL}(V)$ is called a \emph{parabolic subgroup} if it is the pointwise stabiliser of a subset of $V$. As Steinberg has shown in \cite[Theorem 1.5]{St64}, $P$ is also a complex reflection group.
	
	Let $P$ be a parabolic subgroup of $W$.
	By suitably specialising the parameters of the generic Hecke algebra of $P$, the algebra $\cH(P) \otimes R(W)$ embeds into $\cH(W)$ (for a detailed construction, see \cite[\S 2.5]{ChCh}). We call the image of this embedding a \emph{parabolic Hecke subalgebra} of  $\cH(W)$, and to simplify notation we also denote it by $\cH(P)$ . The following conjecture appeared together with the BMM symmetrising trace conjecture   \cite[\S 2.1, Assumption 2(2)(c)]{BMM}:
	
	\begin{conjecture}[Parabolic trace conjecture]\label{parconj}
		The restriction of the canonical symmetrising trace $\tau$ of $\cH(W)$ to the subalgebra $\cH(P)$ is the canonical symmetrising trace on $\cH(P)$.
	\end{conjecture}

	The above conjecture is obviously true for real reflection groups, because of the way $\tau$ is defined with the use of the standard basis. It is also true for $G(de,e,n)$ given that Conjecture \ref{bmmconj} holds, as  remarked in \cite[\S 5.4]{Ch26}. Finally, it has been proved for the following non-real exceptional complex reflection groups:
	\begin{itemize}
		\item the groups $G_4, G_5, G_6, G_7,  G_8, G_{12}, G_{13},G_{22}$ by  \cite{ChCh}; \smallbreak
		\item the groups $G_9, G_{14}$ by \cite{Ch26}.
	\end{itemize}
	
	There is also a counterpart of the BMR freeness conjectures for parabolic subgroups, which was stated by Malle and Rouquier  \cite[Conjecture 2.1]{malrou}, and  is known to hold for real reflection groups \cite[\S 4.4.7, Relation (b)]{gepf}:
	
	\begin{conjecture}[Parabolic freeness conjecture]\label{pfc} The algebra $\mathcal{H}(W)$ is a free left and right $\mathcal{H}(P)$-module of rank $|W|/|P|$.
	\end{conjecture}
	
	The above conjecture has been proved for the following non-real complex reflection groups:
	\begin{itemize}
		\item the groups  $G_4$,  $G_5$, $G_6$, $G_7$, $G_8$,  $G_{12}$, $G_{13}$, $G_{14}$, $G_{16}$, $G_{22}$ by \cite{ChCh}; \smallbreak
		\item the groups $G_9$, $G_{11}$,  $G_{15}$ by \cite{Ch26}.
	\end{itemize}	
	It is also known to hold for the following pairs $(W,P)$:
	\begin{itemize}
		\item  ($G(de,e,n), G(de,e,n-1))$ by   \cite{Sch}; \smallbreak
		\item $(G_9, \mathbb{Z}/4\mathbb{Z})$, $(G_{10},\mathbb{Z}/4\mathbb{Z})$, $(G_{11},\mathbb{Z}/3\mathbb{Z})$, 
		$(G_{15}, \mathbb{Z}/2\mathbb{Z})$, $(G_{15}, \mathbb{Z}/3\mathbb{Z})$ by  \cite{ChCh};\smallbreak
		\item $(G_{25},G_4)$, $(G_{25},\mathbb{Z}/3\mathbb{Z}\times \mathbb{Z}/3\mathbb{Z})$, 
		$(G_{32},G_{25})$ 
		by  \cite{Mar41};\smallbreak
		\item $(G_{26},{G_4})$ 	by  \cite{Mar43};\smallbreak	
		\item $(G_{24}, B_2)$, $(G_{27},B_2)$, $(G_{29},B_3)$, $(G_{31}, A_3)$, $(G_{33}, A_4)$, $(G_{33},D_4)$, $(G_{34},G_{33})$ by  \cite{MaPf}.	
	\end{itemize}

	\section{Reducing the number of parameters}\label{sec:reduce}

Computations in the generic Hecke algebra involve polynomials in all
the parameters $a_{s,j}$, which is their main cost. In this section, we
show that it suffices to compute after specialising each $a_{s,0}$ to
$1$, provided that all the elements, bases and linear forms involved
are homogeneous for a grading of $\mathcal H(W)$ due to Brou\'e, Malle
and Michel. All computations in Sections~\ref{standard}, \ref{bases} and~\ref{rem} are carried out in
this way. We also describe two modular tests, which we use to discard
candidates before any exact computation.
	
	\subsection{The specialisation $\mathrm{pr}$ and the BMM grading}\label{sec-pr}
	
Let $\bar R(W)=\mathbb Z[a_{s,j}\mid j\ge1]$, and let
$\mathrm{pr}\colon R(W)\to\bar R(W)$ be the specialisation that sends
each $a_{s,0}$ to $1$ and fixes the $a_{s,j}$ with $j\ge1$. Let
$\bar{\mathcal H}(W)=\mathcal H_{\mathrm{pr}}$ be the corresponding
specialised algebra. We also write $\mathrm{pr}$ for the induced map
$\mathcal H(W)\to\bar{\mathcal H}(W)$, and, for an $R(W)$-linear form
$E$, we write $\mathrm{pr}(E)$ for the $\bar R(W)$-linear form given by
$\mathrm{pr}(E)(\mathrm{pr}(h))=\mathrm{pr}(E(h))$. The kernel of
$\mathrm{pr}$ is generated by the elements $a_{s,0}-1$, and it is
preserved by the involution $\vee$ of Conjecture~\ref{bmmconj}, since
$a_{s,0}^{-1}-1=-a_{s,0}^{-1}(a_{s,0}-1)$. Hence $\vee$ induces an
involution of $\bar R(W)$, which we denote again by $\vee$; it is given
by $a_{s,j}\mapsto -a_{s,e_s-j}$ for $1\le j\le e_s-1$, and satisfies
$\mathrm{pr}(r^\vee)=\mathrm{pr}(r)^\vee$. Since $\mathcal H(W)$ is free
of rank $|W|$, so is $\bar{\mathcal H}(W)$ over $\bar R(W)$.

A statement proved over $\bar R(W)$ need not hold over $R(W)$. For
instance, the non-zero element $a_{s,0}-1$ maps to $0$ under
$\mathrm{pr}$, and the non-invertible element $2-a_{s,0}$ maps to $1$.
However, this cannot happen for elements that are homogeneous with
respect to the grading of the generic Hecke algebra introduced by
Brou\'e, Malle and Michel \cite[\S1C]{BMM}, which we now recall in
terms of the parameters $a_{s,j}$.

We consider a $\mathbb Z^q$-grading, where $q$ is the number of $W$-conjugacy classes of distinguished pseudo-reflections, that is, the number of $W$-orbits of reflecting hyperplanes. Equivalently, $q$ is the number of conjugacy classes of braided reflections in $B(W)$, since two braided reflections are conjugate in $B(W)$ if and only if their images are conjugate in $W$ \cite[Lemma~2.14 and Theorem~2.17]{BMR}. Choose representatives $\mathbf s_1,\dots,\mathbf s_q$ of these classes, let $s_i\in W$ be the image of $\mathbf s_i$, and write $a_{i,j}:=a_{s_i,j}$. Let $\varepsilon_1,\dots,\varepsilon_q$ denote the standard basis of $\mathbb Z^q$. We give weight $w(\mathbf s)=\varepsilon_i$ to every braided reflection
$\mathbf s$ conjugate to $\mathbf s_i$. Since the braided reflections
generate $B(W)$ \cite[Theorem~2.17(1)]{BMR}, this extends in at most one way
to a group homomorphism $w\colon B(W)\to\mathbb Z^q$, and such an
extension exists by \cite[Proposition~2.16]{BMR}: its components are the
homomorphisms $\pi_1(\delta_{\mathcal C})$ induced by the
discriminants, one for each $W$-orbit $\mathcal C$ of reflecting
hyperplanes. Let $e_i$ be the order of $s_i$, so that $\mathbf s_i$ satisfies the
Hecke relation
\[
\mathbf s_i^{e_i}-a_{i,e_i-1}\mathbf s_i^{e_i-1}-\dots
-a_{i,1}\mathbf s_i-a_{i,0}=0.
\]

In \cite[\S 1C]{BMM}, the ring $\mathbb Z[u_{s,j}^{\pm1}]$ is graded by giving
each $u_{s_i,j}$ degree $\varepsilon_i$. Since $a_{i,j}$ is, up to sign,
the elementary symmetric polynomial of degree $e_i-j$ in
$u_{s_i,1},\dots,u_{s_i,e_i}$, it is homogeneous of degree
$w(a_{i,j})=(e_i-j)\varepsilon_i$, for $0\le j<e_i$, so $R(W)$ is a
graded subring.
In particular, the Hecke relations are homogeneous. They extend additively to monomials, with inverse generators and parameters assigned the negatives of their respective weights. Finally,
for such a monomial $r$ and $b\in B(W)$, we set $w(rb)=w(r)+w(b)$.  By convention, $0$ is regarded as homogeneous of every degree.

\begin{example}
	Let $W=G_4=\langle s_1,s_2\mid s_1s_2s_1=s_2s_1s_2,\ s_1^3=s_2^3=1\rangle$. The generators $s_1,s_2$ are
	conjugate, so $q=1$ and the grading takes values in $\mathbb Z$, with $w(\mathbf s_1)=w(\mathbf s_2)=1$.
	The Hecke relation $\mathbf s^3=a_2\mathbf s^2+a_1\mathbf s+a_0$ gives $w(a_2)=1$, $w(a_1)=2$ and
	$w(a_0)=3$, so that every term of the relation has degree $3$. The braid relation
	$\mathbf s_1\mathbf s_2\mathbf s_1=\mathbf s_2\mathbf s_1\mathbf s_2$ has degree $3$ on both sides.
\end{example}

In \cite[\S1C]{BMM}, it is shown that the generic Hecke algebra over
$\mathbb Z[u_{s,j}^{\pm1}]$ is $\mathbb Z^q$-graded, with $B(W)$ graded
by $w$. The following lemma is the analogous statement over $R(W)$,
with the grading of $R(W)$ induced by that of
$\mathbb Z[u_{s,j}^{\pm1}]$; the proof is the same.

\begin{lemma}\label{4.13}
	The weight function $w$ makes $\cH(W)$ a $\mathbb{Z}^q$-graded $R(W)$-algebra.
	\end{lemma}
\begin{proof}
Since $w$ is a group homomorphism on $B(W)$, the group algebra
$R(W)[B(W)]$ is $\mathbb Z^q$-graded. The Hecke relations are
homogeneous, since each term $a_{i,j}\mathbf s_i^{\,j}$ has degree
$e_i\varepsilon_i$, so the grading descends to the quotient
$\mathcal H(W)$.
\end{proof}

\subsection{Lifting}\label{sec-3.2}
We first show that homogeneous elements of $R(W)$ are determined by their images under $\mathrm{pr}$, in the following sense:

\begin{lemma}\label{liftlem}
	Let $f\in R(W)$ be homogeneous. Then, the following are true:
	\begin{enumerate}
		\item If $\mathrm{pr}(f)=0$, then $f=0$. \smallbreak
		\item  If $\mathrm{pr}(f)=\pm1$, then
		$f=\pm\prod_ia_{i,0}^{n_i}$ for some integers $n_i$; in particular, $f$ is a unit of $R(W)$.
	\end{enumerate} 
\end{lemma}
\begin{proof}
	Let $f$ be homogeneous of degree $\delta\in\mathbb Z^q$. The ring $R(W)$ has a $\mathbb Z$-basis consisting of the
	monomials
	\[
	\mu=\prod_ia_{i,0}^{n_{i,0}}\prod_{j\geq1}a_{i,j}^{n_{i,j}},\qquad n_{i,0}\in\mathbb Z,\ n_{i,j}\geq0\ (j\geq1),
	\]
	and $\mu$ has degree $\sum_i\bigl(e_in_{i,0}+\sum_{j\geq1}(e_i-j)n_{i,j}\bigr)\varepsilon_i$. Hence $f$ is a
	$\mathbb Z$-linear combination of monomials of degree $\delta$. For monomials of a fixed degree $\delta$, the
	exponents $n_{i,j}$ with $j\geq1$ determine $n_{i,0}$, since the $\varepsilon_i$-coordinate of $\delta$ equals
	$e_in_{i,0}+\sum_{j\geq1}(e_i-j)n_{i,j}$. Since $\mathrm{pr}$ sends $\mu$ to $\prod_{i,\,j\geq1}a_{i,j}^{n_{i,j}}$,
	it follows that $\mathrm{pr}$ sends distinct monomials of degree $\delta$ to distinct monomials of $\bar R(W)$.
	
	(1) Write $f=\sum_\mu c_\mu\mu$ with $c_\mu\in\mathbb Z$, the sum running over monomials of degree $\delta$. By the
	above, $\mathrm{pr}(f)=\sum_\mu c_\mu\mathrm{pr}(\mu)$ is a combination of pairwise distinct monomials, which are linearly
	independent over $\mathbb Z$. So $\mathrm{pr}(f)=0$ forces $c_\mu=0$ for all $\mu$, that is, $f=0$.
	
	(2) Suppose $\mathrm{pr}(f)=\epsilon\in\{\pm1\}$. Since the $\mathrm{pr}(\mu)$ are pairwise distinct monomials, $\mathrm{pr}(f)=\epsilon$
	forces $c_\mu=0$ for all $\mu$ except one monomial $\mu_0$ with $\mathrm{pr}(\mu_0)=1$, and $c_{\mu_0}=\epsilon$. The
	monomials with $\mathrm{pr}(\mu)=1$ are those without factors $a_{i,j}$ for $j\geq1$, that is,
	$\mu_0=\prod_ia_{i,0}^{n_i}$. Hence $f=\epsilon\prod_ia_{i,0}^{n_i}$, which is a unit of $R(W)$ since each
	$a_{i,0}$ is.
	\end{proof}
	
Lemma~\ref{liftlem} concerns single elements of $R(W)$. To lift the
statement that a collection spans $\bar{\mathcal H}(W)$, we need a
version for graded modules. Let $I$ be the kernel of
$\mathrm{pr}\colon R(W)\to\bar R(W)$, so that
$\bar{\mathcal H}(W)=\mathcal H(W)/I\mathcal H(W)$.

\begin{lemma}\label{lem:gradedNak}
Let $Q$ be a $\mathbb Z^q$-graded $R(W)$-module. If $Q=IQ$, then $Q=0$.
\end{lemma}

\begin{proof}
Let $L\subseteq\mathbb Z^q$ be the lattice spanned by the degrees
$e_i\varepsilon_i$ of the $a_{i,0}$, and let $\Gamma$ be a set of
representatives of $\mathbb Z^q/L$. For
$n=(n_1,\dots,n_q)\in\mathbb Z^q$, write $a^n=\prod_ia_{i,0}^{n_i}$.
Since the $a_{i,0}$ are units, multiplication by $a^n$ is a bijection
from $Q_\gamma$ onto $Q_{\gamma+\sum_in_ie_i\varepsilon_i}$, so
$Q=\bigoplus_{\gamma\in\Gamma}\bigoplus_{n\in\mathbb Z^q}a^nQ_\gamma$.
Hence, there is a well-defined additive map
$\varphi\colon Q\to\bigoplus_{\gamma\in\Gamma}Q_\gamma$ sending $a^nx$
to $x$ for $x\in Q_\gamma$, $\gamma\in\Gamma$. It is surjective, and it
satisfies $\varphi(a_{i,0}y)=\varphi(y)$ for all $y\in Q$, so it
vanishes on $(a_{i,0}-1)Q$ for every $i$. Since $I$ is generated by the
elements $a_{i,0}-1$, $\varphi$ vanishes on $IQ=Q$. Hence
$\bigoplus_{\gamma\in\Gamma}Q_\gamma=0$, and therefore $Q=0$.
\end{proof}

Braid lifts are homogeneous, since $w$ is a group homomorphism on
$B(W)$. We write $R(W)_\delta$ and $\mathcal H(W)_\gamma$ for the
homogeneous components of degree $\delta$ and $\gamma$.

\begin{definition}{\rm
A linear form $t \in\mathrm{Hom}_{R(W)}(\mathcal H(W),R(W))$ is
\emph{homogeneous of degree $\epsilon\in\mathbb Z^q$} if
$t(\mathcal H(W)_\gamma)\subseteq R(W)_{\gamma+\epsilon}$ for
every $\gamma\in\mathbb Z^q$.}
\end{definition}

The next lemma extends \cite[Lemma~2.4]{BMM}, which treats the standard
basis of an Iwahori--Hecke algebra, to arbitrary bases consisting of
braid lifts.

\begin{lemma}\label{exp}
Let $\mathcal B$ be an $R(W)$-basis of $\mathcal H(W)$ consisting of
braid lifts. Then, for every $b\in\mathcal B$, the linear form
$\lambda^{\mathcal B}_b$ is homogeneous of degree $-w(b)$. In
particular, if $1\in\mathcal B$, then $\lambda^{\mathcal B}_1$ is
homogeneous of degree~$0$.
\end{lemma}

\begin{proof} Let $h\in\mathcal H(W)$ be homogeneous of degree $\gamma$. 
	We write 
	$
	h=\sum_{b\in\cB}c_b b$, where 
	$c_b=\lambda_b^{\cB}(h).$
	Since $R(W)$ is graded, each coefficient $c_b$ admits a unique decomposition into homogeneous components
	$
	c_b=\sum_{\delta}c_{b,\delta},$
	where $c_{b,\delta}$ is homogeneous of degree $\delta$. Substituting these decompositions into the expression for $h$, we obtain
	$
	h=\sum_{b\in\cB}\sum_{\delta}c_{b,\delta}b.
	$
	Since each $b\in\cB$ is homogeneous of degree $w(b)$, every nonzero term $c_{b,\delta}b$ is homogeneous of degree $\delta+w(b)$. Thus, grouping the terms according to their degrees gives
	\[
	h=\sum_{\gamma'}\left(
	\sum_{b\in\cB}c_{b,\gamma'-w(b)}b
	\right),
	\]
	where the expression in parentheses is homogeneous of degree $\gamma'$.
	Now, since $h$ is homogeneous of degree $\gamma$ and  $\mathcal H(W)$ is the direct sum of its homogeneous components (Lemma \ref{4.13}), the uniqueness of the homogeneous decomposition implies that
	\[
	\sum_{b\in\cB}c_{b,\gamma'-w(b)}b=0
	\qquad\text{for every }\gamma'\neq \gamma.
	\]
	On the other hand, since $\mathcal B$ is an $R(W)$-basis, the expansion of any element in this basis is unique. Consequently,
	$
	c_{b,\gamma'-w(b)}=0$
	for every $b\in\cB$
	and every $\gamma'\neq \gamma.$
	Equivalently,
	$
	c_{b,\delta}=0$
	whenever $\delta+w(b)\neq \gamma.$
	It follows that the only possibly nonzero homogeneous component of $c_b$ is the one of degree $\gamma-w(b)$. Therefore,
	$
	c_b=c_{b,\gamma-w(b)},
	$
	and hence $c_b$ is either zero or homogeneous of degree $\gamma-w(b)$, as claimed. This is the statement that $\lambda^{\mathcal B}_b$ is homogeneous of degree $-w(b)$.
\end{proof}

The properties we want to lift, being a basis and being a symmetrising
trace, are expressed by determinants of matrices whose entries are
coefficients as in the above lemma. The following lemma shows
that such determinants are homogeneous.

\begin{lemma}\label{homog}
	Let $M$ be an $n\times n$ matrix with entries in $R(W)$. Suppose that there exist weights
	$\alpha_1,\ldots,\alpha_n$ and $\beta_1,\ldots,\beta_n$ in $\mathbb{Z}^q$ such that every entry $M_{ij}$ is
	either zero or homogeneous of degree $\alpha_i+\beta_j$. Then $\det(M)$ is homogeneous of degree
	$\sum_{i=1}^n(\alpha_i+\beta_i)$.
\end{lemma}
\begin{proof}
	By the Leibniz formula,
	\[
	\det(M)=\sum_{\sigma\in\mathfrak S_n}\operatorname{sgn}(\sigma)\prod_{i=1}^nM_{i,\sigma(i)}.
	\]
	We fix $\sigma\in\mathfrak S_n$. Each factor $M_{i,\sigma(i)}$ is either zero or homogeneous of degree
	$\alpha_i+\beta_{\sigma(i)}$, and a product of homogeneous elements is homogeneous of degree the sum of the
	degrees. Hence, if the product $\prod_iM_{i,\sigma(i)}$ is nonzero, it is homogeneous of degree
	\[
	\sum_{i=1}^n(\alpha_i+\beta_{\sigma(i)})=\sum_{i=1}^n\alpha_i+\sum_{i=1}^n\beta_{\sigma(i)}.
	\]
	This degree does not depend on $\sigma$, since $\sigma$ permutes the indices and hence
	$\sum_i\beta_{\sigma(i)}=\sum_i\beta_i$. Therefore $\det(M)$ is a sum of homogeneous elements of the same
	degree, and so it is either zero or homogeneous of degree $\sum_{i=1}^n(\alpha_i+\beta_i)$.
\end{proof} 
We can now show that bases and symmetrising traces can be detected
after specialisation, provided that all the objects involved are
homogeneous.
	
	\begin{proposition}\label{prop:lift}
\begin{enumerate}
\item Let $X$ be a collection of $|W|$ homogeneous elements of
  $\mathcal H(W)$. If $\mathrm{pr}(X)$ spans $\bar{\mathcal H}(W)$, then
  $X$ is an $R(W)$-basis of $\mathcal H(W)$.
\item Let $\mathcal B$ be an $R(W)$-basis of $\mathcal H(W)$ consisting
  of homogeneous elements, and let $t$ be a homogeneous linear
  form on $\mathcal H(W)$. If $\mathrm{pr}(t)$ is a trace
  function, then so is $t$; if $\mathrm{pr}(t)$ is a
  symmetrising trace, then so is $t$.
\end{enumerate}
\end{proposition}

\begin{proof}
(1) Let $N$ be the $R(W)$-submodule of $\mathcal H(W)$ spanned by $X$.
Since $X$ consists of homogeneous elements, $N$ is a graded submodule,
and $Q=\mathcal H(W)/N$ is a graded $R(W)$-module. Since
$\mathrm{pr}(X)$ spans $\bar{\mathcal H}(W)=\mathcal H(W)/I\mathcal H(W)$,
we have $\mathcal H(W)=N+I\mathcal H(W)$, that is, $Q=IQ$. By
Lemma~\ref{lem:gradedNak}, $Q=0$, so $X$ spans $\mathcal H(W)$, and it
is an $R(W)$-basis by Proposition~\ref{span}.

(2) Let $t$ be homogeneous of degree $\epsilon$. For braid lifts
$x,y$, the elements $t(xy)$ and $t(yx)$ are zero or
homogeneous of degree $w(x)+w(y)+\epsilon$, and their images under
$\mathrm{pr}$ agree if $\mathrm{pr}(t)$ is a trace function. By
Lemma~\ref{liftlem}(1), they are equal, and since braid lifts span
$\mathcal H(W)$, $t$ is a trace function. If moreover
$\mathrm{pr}(t)$ is a symmetrising trace, the entries
$t(bb')$ of $\mathrm{Gram}_t(\mathcal B,\mathcal B)$ are zero
or homogeneous of degree $w(b)+w(b')+\epsilon$, so its determinant is
homogeneous by Lemma~\ref{homog}. Its image under $\mathrm{pr}$ is
$\det\mathrm{Gram}_{\mathrm{pr}(t)}(\mathrm{pr}(\mathcal B),
\mathrm{pr}(\mathcal B))$, a unit of $\bar R(W)$. By
Lemma~\ref{liftlem}(2), it is a unit, so $t$ is a symmetrising
trace.
\end{proof}

Combining Proposition~\ref{prop:lift} with the reformulation~\eqref{extra con2} of
Condition~(3), we obtain a criterion for the BMM symmetrising trace
conjecture that only involves computations over $\bar R(W)$. We
formalise the conditions to be verified as follows.

\begin{definition}\label{def:pradm}{\rm
A collection $\mathcal B\subseteq\mathcal H(W)$ is
\emph{$\mathrm{pr}$-admissible} if the following hold:
\begin{enumerate}
\item $\mathcal B$ is the set of lifts of a section of
  $B(W)\twoheadrightarrow W$; \smallbreak
\item $\mathrm{pr}(\mathcal B)$ is a basis of $\bar{\mathcal H}(W)$;\smallbreak
\item $\bar\tau=\lambda_1^{\mathrm{pr}(\mathcal B)}$ is a symmetrising
  trace; and\smallbreak
\item $\bar\tau(\mathrm{pr}(b^{-1}\boldsymbol\pi))=0$ for all
  $b\in\mathcal B\setminus\{1\}$.
\end{enumerate}}
\end{definition}

The following theorem shows that these conditions suffice.

\begin{theorem}\label{thm:liftall}
Let $\mathcal B\subseteq\mathcal H(W)$ be $\mathrm{pr}$-admissible. Then
$\mathcal B$ is a section basis of $\mathcal H(W)$, and
$\lambda_1^{\mathcal B}$ is the canonical symmetrising trace of
$\mathcal H(W)$.
\end{theorem}

\begin{proof}
By Proposition~\ref{prop:lift}(1), $\mathcal B$ is an $R(W)$-basis of
$\mathcal H(W)$, hence a section basis. Then
$\tau:=\lambda_1^{\mathcal B}$ is homogeneous of degree $0$ by
Lemma~\ref{exp}, and $\mathrm{pr}(\tau)=\bar\tau$; by
Proposition~\ref{prop:lift}(2), $\tau$ is a symmetrising trace. Since
$\mathcal B$ is a section basis, $\tau$ satisfies Condition~(2)
of Conjecture \ref{bmmconj}. For $b\in\mathcal B\setminus\{1\}$, the element
$\tau(b^{-1}\boldsymbol\pi)$ is zero or homogeneous, and its image
under $\mathrm{pr}$ is $0$; hence it is $0$ by Lemma~\ref{liftlem}(1).
Condition~(3) now follows from \eqref{extra con2}.
\end{proof}

In Sections~\ref{bases} and~\ref{rem}, we construct a $\mathrm{pr}$-admissible
collection for every non-real exceptional complex reflection group;
Theorem~\ref{thm:liftall} then yields Conjecture \ref{bmmconj}.

\subsection{Modular tests}\label{modtests}
	
	Our strategy for proving the BMM symmetrising trace conjecture is to find a suitable section basis $\cB$ of $\cH(W)$  such that $\lambda_1^{\cB}$ is the canonical symmetrising trace. Since $\cB$ is a section basis, Condition~(2) of Conjecture~\ref{bmmconj} is automatically satisfied. It therefore suffices to show that the Gram matrix $\gram_{\lambda^{\cB}_1}(\cB,\cB)=\bigl(\lambda_1^{\cB}(bb')\bigr)_{b,b'\in\cB}$	is symmetric and has determinant a unit in $R(W)$, and to verify Condition~\eqref{extra con2}. This is the approach adopted in all previously established cases of the conjecture for exceptional complex reflection groups.  In this subsection we provide two \emph{modular tests}
	that help us find such bases. The tests are performed in a specialised algebra,
	where the calculations are much faster because there are no parameters. They
	provide only necessary conditions: if a set fails one of the tests, it cannot be
	such a basis and we discard it, whereas if it passes both tests, it has to be examined further by exact computations over $\bar R(W)$, whose results lift to $R(W)$ by \S \ref{sec-3.2}.
	
	Let $\cA=\{a_1,\ldots,a_{|W|}\}$ be an $R(W)$-basis of $\cH(W)$, which we use as
	a reference basis. Let $\bk$ be a field (usually $\bk=\mathbb F_p$ for some big prime $p$) and let $\theta\colon R(W)\to\bk$ be a ring homomorphism such that
	$\theta(a_{s,0})=1$ for all $s$. Such a $\theta$ factors through the specialisation $\mathrm{pr}$ of \S\ref{sec-pr}: we have $\theta=\bar\theta\circ\mathrm{pr}$ for a unique ring homomorphism $\bar\theta\colon\bar R(W)\to \bk$. We first introduce some notation. 
	
	For an ordered collection
	$X=\{x_1,\ldots,x_{|W|}\}\subseteq\cH(W)$, we define $P_{\cA}(X)$ to be the
	matrix whose columns are the coordinates of the elements of $X$ with respect to
	$\cA$, namely
	\[
	P_{\cA}(X)=
	\begin{pmatrix}
		\lambda^{\cA}_{a_1}(x_1) & \cdots & \lambda^{\cA}_{a_1}(x_{|W|})\\
		\vdots & \ddots & \vdots\\
		\lambda^{\cA}_{a_{|W|}}(x_1) & \cdots & \lambda^{\cA}_{a_{|W|}}(x_{|W|})
	\end{pmatrix}.
	\]
	For a linear form $t\in\Hom_{R(W)}(\cH(W),R(W))$ and two ordered collections
	$X,Y$ of $|W|$ elements, we define $\gram_t(X,Y)=\bigl(t(xy)\bigr)_{x\in X,\,y\in Y}$.
	We denote by $P_{\cA}^{\theta}(X)$ and $\gram_t^{\theta}(X,Y)$ the corresponding
	matrices in the specialised algebra. If $X=(x_1,\dots,x_{|W|})$ is a collection of invertible elements of $\mathcal H(W)$, later on we will use  $Y=X^{-1}=(x_1^{-1},\dots,x_{|W|}^{-1})$. Since $R(W)$ is a Laurent polynomial ring in the $a_{s,0}$ over the polynomial ring $\mathbb Z[a_{s,j}\mid j\ge1]$, which is a domain, its units are the elements $\pm\prod_s a_{s,0}^{m_s}$, $m_s\in\mathbb Z$; hence $\theta$ maps every unit of $R(W)$ to $\pm1$. We now give a test, which checks whether a collection $\cB$ of $|W|$ elements of $\cH(W)$
	can be an $R(W)$-basis of $\cH(W)$.
	
	\begin{lemma}[Modular test 1]\label{mtest}
		Let $\cB$ be an ordered collection of $|W|$ elements of $\cH(W)$. If $\cB$ is an
		$R(W)$-basis of $\cH(W)$, then $\det P_{\cA}^{\theta}(\cB)=\pm1$.
	\end{lemma}
	\begin{proof}
		Since $\cA$ and $\cB$ are both $R(W)$-bases, $P_{\cA}(\cB)$ is invertible over
		$R(W)$, so its determinant is a unit of $R(W)$. Hence
		$\det P_{\cA}^{\theta}(\cB)=\theta(\det P_{\cA}(\cB))=\pm1$.
	\end{proof}
	
	We now provide a second test, which checks whether the linear map $\tau=\lambda^{\cB}_1$ can be
	the canonical symmetrising trace of $\cH(W)$.
	
	\begin{lemma}[Modular test 2]\label{mtest2}
		Let $\cB$ be an $R(W)$-basis of $\cH(W)$ containing $1$ such that $\lambda^{\cB}_1$
		is a symmetrising trace on $\cH(W)$. Then
		$\gram^{\theta}_{\lambda^{\cB}_1}(\cB,\cB)$ is symmetric and has determinant $\pm1$.
	\end{lemma}
	\begin{proof}
		Since $\lambda^{\cB}_1$ is a symmetrising trace, $\gram_{\lambda^{\cB}_1}(\cB,\cB)$
		is symmetric, and the associated bilinear form induces an isomorphism
		$\cH(W)\to\Hom_{R(W)}(\cH(W),R(W))$. Therefore this matrix is invertible over
		$R(W)$ and its determinant is a unit. Applying $\theta$ preserves symmetry and
		gives determinant $\pm1$.
	\end{proof}
	
	\begin{remark}
		In our computations we use two primes. The default is the Mersenne prime
		$p=2^{61}-1$: it is large, so that a candidate which is not a basis is very
		unlikely to pass the tests by accident, and arithmetic modulo $p$ is fast on
		64-bit machines. When faster computations are needed we use the smaller prime
		$p=15013$, for which products of residues fit in 32 bits; the price is that
		accidental passes are more likely. In both cases we take $\theta(a_{s,0})=1$.
		Lemmas \ref{mtest} and \ref{mtest2} hold for
		any choice of the values $\theta(a_{s,i})$, $i\ge1$; we choose them 
		pseudo-randomly, so that a candidate is unlikely
		to pass the tests by accident. The modular tests only serve to
		discard candidates.	\end{remark}
	
	\section{Standard bases for complex reflection groups}\label{standard}
	
	As we saw in  \S \ref{sec-2.2}, Iwahori--Hecke algebras of real reflection groups admit a standard basis, whose construction relies on the existence of a length function $\ell$ that associates to each element $w$ of the group the minimal number of generators needed to express $w$. For complex reflection groups, however, there is no analogous notion of length, and hence no immediate counterpart of this basis. In this section, we introduce a natural standard basis for the generic Hecke algebras of exceptional complex reflection groups, through 
	a combinatorial reinterpretation of the standard basis in the real case. 
		
	Let $W$ be a complex reflection group, $U$ a subgroup of $W$,
	and $A$ a finite generating set for $W$, not necessarily minimal and not
	necessarily equal to the set $S$ of standard generators. Each element of
	$A$ that is not a standard generator is given as a fixed word in the
	standard generators and their inverses. 
	
	For any word $w$ in the elements of $A$, we will denote by bold letters the exact same word in the generic Hecke algebra $\cH(W)$, i.e., the element of $\cH(W)$ obtained by
	replacing each letter by the corresponding word in the braided
	reflections and their inverses. This is also an element of $B(W)$, and we use the same symbol for both. 
	By abuse of notation, we 
	write $w$ for the image of $\bw$ in $W$. Thus, an equality $w=w'$ between two words always means equality of their images in $W$, whereas $\bw=\bw'$ means equality in $B(W)$, or in $\cH(W)$, as specified. In particular, $e$ denotes both the
	empty word and the identity element of $W$.
	For example, in the case of $G_4=\langle s_1,s_2 \,|\, s_1s_2s_1=s_2s_1s_2,\, s_1^3=s_2^3=1\rangle$, if $w=s_1s_2^2s_1^{-1}$, then $\bw=\bs_1\bs_2^2\bs_1^{-1}$, 
	and  if $w'=s_1s_2^{-1}s_1^{2}$, then $\bw'=\bs_1\bs_2^{-1}\bs_1^{2}$. Even though $w=w'$ in $W$, we have $\bw \neq \bw'$ in $B(W)$ and in $\cH(W)$.

	\begin{definition}\label{treec}{\rm
			A \emph{coset tree} with respect to $(A,U)$ is an ordered tuple
			\[
			\mathcal{T}=(w_1,\ldots,w_{[W:U]})
			\]
			of words in the elements of $A$ such that
			\begin{enumerate}
				\item $w_1=e$;\smallbreak
				\item $Uw_i\neq Uw_j$ whenever $i\neq j$;\smallbreak
				\item (\emph{Schreier property}) for every $i>1$, there exist $j<i$ and $a\in A$ such that $w_i$ is the word $w_ja$.
		\end{enumerate}}
	\end{definition}

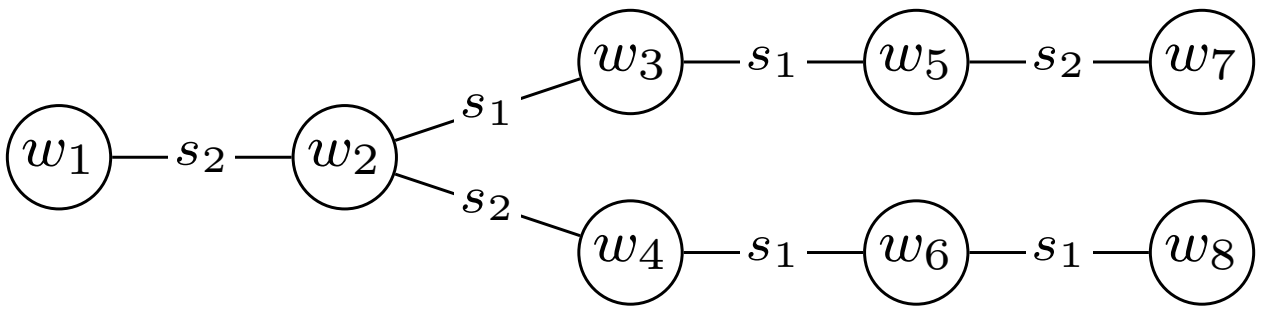
	\begin{figure}[H]
		\centering
		\scalebox{1.4}{\begin{tikzpicture}[x=1cm,y=1cm,
  v/.style={circle,draw,inner sep=1pt,font=\scriptsize},
  e/.style={font=\tiny,fill=white,inner sep=0.6pt}]
  \node[v] (n1) at (0.0,0.0) {$w_{1}$};
  \node[v] (n2) at (1.35,0.0) {$w_{2}$};
  \node[v] (n3) at (2.7,0.45) {$w_{3}$};
  \node[v] (n4) at (2.7,-0.45) {$w_{4}$};
  \node[v] (n5) at (4.05,0.45) {$w_{5}$};
  \node[v] (n6) at (4.05,-0.45) {$w_{6}$};
  \node[v] (n7) at (5.4,0.45) {$w_{7}$};
  \node[v] (n8) at (5.4,-0.45) {$w_{8}$};
  \draw (n1) -- node[e] {$s_2$} (n2);
  \draw (n2) -- node[e] {$s_1$} (n3);
  \draw (n2) -- node[e] {$s_2$} (n4);
  \draw (n3) -- node[e] {$s_1$} (n5);
  \draw (n4) -- node[e] {$s_1$} (n6);
  \draw (n5) -- node[e] {$s_2$} (n7);
  \draw (n6) -- node[e] {$s_1$} (n8);
\end{tikzpicture}}
		\[
		(w_1,\ldots,w_8)
		=
		(e,\ s_2,\ s_2s_1,\ s_2^2,\ s_2s_1^2,\ s_2^2s_1,\ 
		s_2s_1^2s_2,\ s_2^2s_1^2).
		\]
		\caption{A coset tree for $G_4$ with respect to $(\{s_1,s_2\}, \{e,s_1,s_1^2\})$.}
		\label{fig:one}
	\end{figure}

	As illustrated in Figure~\ref{fig:one}, where an example for $G_4$ is given,
	the elements of $\mathcal{T}$ can be viewed as the vertices of a tree rooted at $e$, where each edge corresponds to right multiplication by an element of $A$. This interpretation motivates the term ``tree''.
	
	We denote by $\cM_A^U$ the set of all coset trees with respect to $(A,U)$. Notice that $\cM_A^U$ is finite. Indeed, every element of a coset tree is obtained from the root by successively multiplying by elements of $A$. Since a coset tree contains exactly $[W:U]$ elements, we need at most $[W:U]-1$ such multiplications to reach any of its elements. As $A$ is finite, there are only finitely many possibilities at each step, and therefore only finitely many coset trees. 
		
	Our goal in this section is to use coset trees to construct basis candidates for the associated Hecke algebra, by taking $U$ to be from now on a parabolic subgroup $P$ of $W$. 
	This way, rather than working with elements of $W$, we choose one representative for each right coset of $P$ in $W$. This idea has already been used in \cite{Ch17, ChPf1, ChPf}.

	\subsection{$\BFS_o$ bases for real reflection groups}
	We will now see  how the standard basis arises naturally from a coset tree with respect to the trivial parabolic subgroup $P=\{e\}$.
	Let $(W,A)$ be a finite Coxeter system with $|A|=m$ and
	$
	o=(s_1,\ldots,s_m)$
	an ordering of $A$. We construct a coset tree $\mathcal{T}_o$ with respect to $(A,\{e\})$ by applying the breadth-first search (BFS) algorithm as follows:
	
	\begin{enumerate}
		\item Start with the root $e$, at level $0$. \smallbreak
		
		\item Suppose that all nodes up to level $L-1$ have been constructed. 
		Following the ordering established by the tree, for each node $w$ at level $L-1$, consider
		\[
		ws_1,\ldots,ws_m
		\]
		in the order prescribed by $o$. If the corresponding element of $W$ has already occurred in the tree, discard it. Otherwise, add it to the tree as a child of $w$, at level $L$. \smallbreak
		
		\item Terminate when all elements of $W$ have been reached.
	\end{enumerate}

	\begin{definition}{\rm
			For the coset tree $\mathcal{T}_o$ constructed above, we define the \emph{$\BFS_o$ basis candidate} by
			\[
			\BFS_o:=\{{\mathbf{w}}\mid w\in\mathcal{T}_o\}\subseteq\cH(W).
			\]}
	\end{definition}
	\begin{example}
		Let $W=S_3$ with generators
		$
		s_1=(1\,2),\, s_2=(2\,3),$
		and consider the ordering $o=(s_1,s_2)$. Starting from $w_1=e$, the  level 1 consists of $w_2=s_1$ and $w_3=s_2$.
		
		At level 2, the algorithm considers
		$
		s_1^2,\,s_1s_2,\, s_2s_1,\, s_2^2,$ in that order.
		Since $s_1^2=s_2^2=e$, these elements are discarded. The elements $w_4=s_1s_2$ and $w_5=s_2s_1$ are new and are added to the tree.  
		
		At level $3$, starting from $s_1s_2$, we obtain $s_1s_2s_1$, which has not yet appeared in the tree, while $s_1s_2^2=s_1$ is discarded. Next, from $s_2s_1$, we obtain $s_2s_1^2=s_2$ and $s_2s_1s_2=s_1s_2s_1$, which are both discarded. So $w_6=s_1s_2s_1$.
		
		We have reached all elements of $S_3$ and, therefore, the construction terminates. Thus
		\[
		\BFS_o=
		\{1,{\mathbf{s}_1},{\mathbf{s}_2},
		{\mathbf{s}_1\mathbf{s}_2},
		{\mathbf{s}_2\mathbf{s}_1},
		{\mathbf{s}_1\mathbf{s}_2\mathbf{s}_1}\}.
		\]
	\end{example}
	
	As the above example suggests, in the real case, every $\BFS_o$ basis candidate coincides with the standard basis of the corresponding Iwahori--Hecke algebra:
	
	\begin{proposition}\label{real}
		Let $(W,S)$ be a finite Coxeter system and let $o$ be an ordering of $S$. The elements occurring at level $L$ of $\mathcal T_o$ are exactly the elements of Coxeter length $L$. Consequently, $\BFS_o$ coincides with the standard basis of $\cH(W)$.
	\end{proposition}
	
	\begin{proof}
		We proceed by induction on $L$. At level $0$, the only element is $e$, and  $\ell(e)=0$.
		At level $1$, we obtain all elements of $S$, and these are exactly the elements of Coxeter length $1$.
		
		Suppose now that $L>1$ and that the statement is true for all levels up to level $L-1$. Let $w$ be an element added at level $L$. By construction, $w=vs$ for some element $v$ at level $L-1$ and some $s\in S$. By induction hypothesis, $\ell(v)=L-1$. Since
		$
		\ell(vs)=\ell(v)\pm1,
		$
		we have
		$
		\ell(w)\in\{L-2,L\}.
		$	
		If $\ell(w)=L-2$, then, by induction hypothesis, $w$ should have  been reached at  level $L-2$ and therefore would not have been added at level $L$, which is a contradiction. Hence $\ell(w)=L$.
		
		Thus, the expression for $w$ determined by its path from the root is reduced. Since the algorithm reaches every element of $W$, it determines a reduced expression for each $w\in W$. By Matsumoto's theorem, the corresponding Hecke algebra element is independent of the chosen reduced expression. Consequently, $\BFS_o$ coincides with the standard basis of $\cH(W)$.
	\end{proof}
	
	\subsection{$\BFS_o$ bases for complex reflection groups}
	
	Since the breadth-first construction does not rely on a length function, it can also be applied to complex reflection groups.
	
	Let $W$, $P$ and $A$ be as at the beginning of this section. Let
	$
	o=(a_1,\ldots,a_r)
	$
	be an ordering of  $A$. Since $A$ generates $W$, every right coset of $P$ in $W$ is reached from $P$ by successive right multiplications by elements of $A$.
	We construct a coset tree $\mathcal T_o$ with respect to $(A,P)$ by applying the breadth-first search (BFS) algorithm as follows:
	
	\begin{enumerate}	
		\item Start with the root $e$, at level $0$.\smallbreak
		\item Suppose that all nodes up to level $L-1$ have been constructed. 
		Following the ordering established by the tree, for each node $w$ at level $L-1$, consider
		\[
		wa_1,\ldots,wa_r
		\]
		in the order prescribed by $o$. If a right coset representative of $Pwa_i$ has already occurred in the tree, discard it. Otherwise, add it to the tree as a child of $w$, at level $L$. \smallbreak
		
		\item Terminate when the tree contains representatives of all right cosets in $P\backslash W$.
	\end{enumerate}

		We now associate a basis candidate of $\cH(W)$ with the tree $\mathcal T_o$. By the validity of the BMR freeness conjecture for $P$, the generic Hecke algebra of $P$ is free of rank $|P|$ over $R(P)$. Extending scalars from $R(P)$ to $R(W)$, we obtain a free $R(W)$-module of rank $|P|$, which we continue to denote by $\cH(P)$, as we did in \S\ref{sec-2.3}. We fix an $R(W)$-basis $\mathcal{B}_P$ of $\cH(P)$.
		
		\begin{definition}{\rm
				The \emph{$\BFS_o$ basis candidate} associated with the ordering $o$ and the fixed basis $\cB_P$ is the following collection of elements of $\cH(W)$:
				$$\BFS_o:=\mathcal{B}_P \cdot \mathcal{T}_o= \{b\mathbf{w} \mid b \in \mathcal{B}_P, \ w\in\mathcal{T}_o\}.$$}
		\end{definition}

		\begin{example}\label{g4bfs} The group
			$W=G_4=\langle s_1,s_2 \,|\, s_1s_2s_1=s_2s_1s_2,\, s_1^3=s_2^3=1\rangle$ has order $24$. We take $P=\langle s_1\rangle$ and $A=\{s_1,s_2\}$. We have $[W:P]=8$.
			For the ordering $o=(s_1,s_2)$, the breadth-first construction proceeds as follows:
			
			\begin{itemize}
				\item At level $0$, we start with $w_1=e$.
				
				\item At level $1$, we consider the elements $es_1=s_1$ and $es_2=s_2$. Since $s_1\in P$, the first element is discarded. The second represents a new coset, so we add
				$
				w_2=s_2.
				$
				
				\item At level $2$, both elements
				$
				w_3=w_2s_1=s_2s_1,$ and
				$w_4=w_2s_2=s_2^2
				$
				represent new cosets and are added to the tree as children of $w_2$.
				
				\item At level $3$, the algorithm considers the elements $w_3s_1,w_3s_2, w_4s_1,w_4s_2$, in that order. We first add
				$
				w_5=w_3s_1=s_2s_1^2
				$ as a child of $w_3$.
				On the other hand,
				$
				w_3s_2=s_2s_1s_2=s_1s_2s_1=s_1w_3\in Pw_3,$
				so this element is discarded.
				Similarly, we add
				$
				w_6=w_4s_1=s_2^2s_1
				$ as a child of $w_4$,
				while
				$
				w_4s_2=s_2^3=e
				$
				is discarded.
				
				\item At level $4$, the algorithm considers the elements $w_5s_1,w_5s_2, w_6s_1,w_6s_2$, in that order. 
				$
				w_5s_1=s_2s_1^3=s_2=w_2,
				$
				so this element is discarded. We then add
				$
				w_7=w_5s_2=s_2s_1^2s_2
				$ as 	a child of $w_5$.	
				We also add
				$
				w_8=w_6s_1=s_2^2s_1^2
				$ as a child of $w_6$.
			\end{itemize}		
		We have now reached all eight cosets of $P$ in $W$, so the algorithm terminates. The resulting coset tree is
			\[
			\mathcal T_o=
			\{e,\ s_2,\ s_2s_1,\ s_2^2,\ s_2s_1^2,\ s_2^2s_1,
			s_2s_1^2s_2,\ s_2^2s_1^2\},
			\]
			which is precisely the tree illustrated in Figure~\ref{fig:one} at the beginning of this section.
			Taking $\mathcal{B}_P=\{1,\bs_1,\bs_1^2\}$, the corresponding $\BFS_o$ basis candidate is
			\[
			\begin{aligned}
				\BFS_o=\big\{&
				1,\ {\mathbf{s}_2},
				{\mathbf{s}_2}{\mathbf{s}_1},
				{\mathbf{s}_2}^2,
				{\mathbf{s}_2}{\mathbf{s}_1}^2,
				{\mathbf{s}_2}^2{\mathbf{s}_1},
				{\mathbf{s}_2}{\mathbf{s}_1}^2{\mathbf{s}_2},
				{\mathbf{s}_2}^2{\mathbf{s}_1}^2,\\
				&
				{\mathbf{s}_1},
				{\mathbf{s}_1}{\mathbf{s}_2},
				{\mathbf{s}_1}{\mathbf{s}_2}{\mathbf{s}_1},
				{\mathbf{s}_1}{\mathbf{s}_2}^2,
				{\mathbf{s}_1}{\mathbf{s}_2}{\mathbf{s}_1}^2,
				{\mathbf{s}_1}{\mathbf{s}_2}^2{\mathbf{s}_1},
				{\mathbf{s}_1}{\mathbf{s}_2}{\mathbf{s}_1}^2{\mathbf{s}_2},
				{\mathbf{s}_1}{\mathbf{s}_2}^2{\mathbf{s}_1}^2,\\
				&
				{\mathbf{s}_1}^2,
				{\mathbf{s}_1}^2{\mathbf{s}_2},
				{\mathbf{s}_1}^2{\mathbf{s}_2}{\mathbf{s}_1},
				{\mathbf{s}_1}^2{\mathbf{s}_2}^2,
				{\mathbf{s}_1}^2{\mathbf{s}_2}{\mathbf{s}_1}^2,
				{\mathbf{s}_1}^2{\mathbf{s}_2}^2{\mathbf{s}_1},
				{\mathbf{s}_1}^2{\mathbf{s}_2}{\mathbf{s}_1}^2{\mathbf{s}_2},
				{\mathbf{s}_1}^2{\mathbf{s}_2}^2{\mathbf{s}_1}^2
				\big\}.
			\end{aligned}	
			\]
		\end{example}
		
		\begin{remark}In contrast to the real case, where all orderings yield the same standard basis by Proposition~\ref{real}, different orderings may give rise to different basis candidates in the complex case.
		\end{remark}
		
		Our goal in the remainder of this section is to show that, for every irreducible exceptional complex reflection group, a suitable choice of ordering $o$ yields an $R(W)$-basis $\BFS_o$ of $\cH(W)$. This provides a uniform proof of the BMR freeness conjecture for all these groups.

		\subsection{Completing coset tables}\label{coset}

		We will now apply Proposition \ref{span} to a basis candidate $\BFS_o$, associated with a fixed ordering $o$ of $A$. A standard approach is to consider the $R(W)$-submodule
		$\operatorname{Span}_{R(W)}(\BFS_o)
		$
		and prove that it is a right ideal of $\cH(W)$. Since $1\in \BFS_o$, it then follows that $\cH(W)=\operatorname{Span}_{R(W)}(\BFS_o)$. For this purpose, it suffices to verify that
		$
		x{\mathbf{s}}\in \operatorname{Span}_{R(W)}(\BFS_o)$
		for all $x\in \BFS_o$ and $s\in S,$
		where $S$ is the set of standard generators.
		
		Rather than working directly with all $|W|$ elements of $\BFS_o$, we take advantage of its parabolic structure. Recall that
		$
		\BFS_o=\{b{\mathbf{w}}\mid b \in \mathcal{B}_P, \ w\in\mathcal{T}_o\},
		$
		where $\mathcal{B}_P$ is a fixed basis of the parabolic Hecke subalgebra $\cH(P)$. We therefore work with the elements ${\mathbf w}$, for $w\in\mathcal T_o$, and seek to express their products with the elements of $S$ as linear combinations of these elements with coefficients in $\cH(P)$.
		
		Our aim is not only to show that these products can be written as linear combinations of elements of $\BFS_o$, but also to determine the corresponding coefficients explicitly. This information will also be useful later in our study of the BMM symmetrising trace conjecture. These coefficients are recorded in a coset table, in a method developed in \cite{ChPf1,ChPf}.

		We denote by $z\in W$  the image of the generator $\mathbf z$ of $Z(B(W))$. We set $\Sigma:=S\cup S^{-1}\cup\{z,z^{-1}\}$, and we assume that $A \subseteq\Sigma$.
		We aim to construct a \emph{coset table} $\mathcal{C}_{\BFS_o}$ with $[W:P]$ rows and $|\Sigma|$ columns, indexed by the elements of $\mathcal T_o$ and the elements of $\Sigma$ respectively.  
		Let  $w\in\mathcal T_o$ and $\sigma\in\Sigma$. If 
		\[
		{\mathbf w}{\boldsymbol{\sigma}}
		=
		\sum_{w'\in\mathcal T_o}
		h_{w,\sigma,w'}{\mathbf{w'}},
		\,\quad
		{\rm with }\,\, h_{w,\sigma,w'}\in\cH(P),
		\]
		then $$\mathcal C_{\BFS_o}(w,\sigma):=\sum_{w'\in\mathcal T_o}
		h_{w,\sigma,w'}{\mathbf{w'}}.$$
		Thus, completing the coset table amounts to finding such expressions for all the corresponding products. Before $\mathrm{\BFS}_o$ is known to be a basis, these expressions need not be unique; the table records the ones found by the procedure below.
		Although the generators in $S$ suffice to establish that $\operatorname{Span}_{R(W)}(\BFS_o)$ is a right ideal, the additional elements in $\Sigma$ can facilitate the completion of the coset table.

		A \emph{rewriting rule} is an identity
		$\boldsymbol\sigma=\sum_i c_i\mathbf v_i$ in $H(W)$, with $\sigma\in\Sigma$,
		$c_i\in R(W)$ and $v_i$ words in $\Sigma$. Thus a rewriting rule replaces a single letter, the one indexing a column of the table, by a combination of words whose letters also index columns. We use three kinds of rules:
		\begin{description}
			\item[\textsf{gen}] the positive Hecke relations \eqref{Hecker2}, solved for $\mathbf s$:
			$$\mathbf s=a_{s,e_s-1}+a_{s,e_s-2}\mathbf s^{-1}+\dots+a_{s,0}\mathbf s^{-(e_s-1)};$$
			\item[\textsf{inv}] the inverse Hecke relations \eqref{invhecke}, expressing $\mathbf s^{-1}$ in
			terms of $1,\mathbf s,\dots,\mathbf s^{e_s-1}$;\smallbreak
			\item[\textsf{rel}] the braid relations, and the relations expressing $\mathbf z$ as a
			word in the generators, solved for a single letter; for instance, in
			$G_4$ the braid relation gives
			$\mathbf s_2=\mathbf s_1^{-1}\mathbf s_2^{-1}\mathbf s_1\mathbf s_2\mathbf s_1$.
		\end{description}
		The coset table $\mathcal C:=\mathcal C_{\BFS_o}$ is then completed in three steps.
		\begin{description}
			\item[Initialisation] For each edge of $\cT_o$, that is, whenever $w'$ is the word $wa$ with $w, w' \in  \cT_o$ and $a \in A$, set $\mathcal C(w, a): = \bw'$ and $\mathcal C(w', a^{-1}) := \bw$. In addition, for each $\sigma \in  \Sigma \cap P$, set $\cC(w_1, \sigma) := \boldsymbol{\sigma}\, \bw_1$ (here the coefficient of $\bw_1$ is $\boldsymbol{\sigma} \in \cH(P)$).\smallbreak
			
			\item[Sweeps] Repeatedly run through all undefined entries $\mathcal C(w, \sigma)$. Whenever some rewriting rule $\boldsymbol\sigma=\sum_i c_i\mathbf v_i$ allows us to compute coefficients for $\bw \cdot \sum_i c_i\mathbf v_i$ letter by letter from entries already defined, we take these as $\mathcal C(w, \sigma)$. Each new entry can be used as soon as it is defined.\smallbreak
			
			\item[Termination] Stop when $\mathcal C$ is complete, or when a full sweep defines no new entry.
			
		\end{description}
		
		If the columns indexed by $S$ are completed, then $\operatorname{Span}_{R(W)}(\BFS_o)
		$
		is a right ideal of $\cH(W)$ containing $1$, and so $\BFS_o$ is an $R(W)$-basis of $\cH(W)$. The above procedure is related to the vector enumeration of Linton~\cite{Lin}, the linear analogue of Todd--Coxeter coset enumeration~\cite{toco, HEO}, and coset tables of this kind were used in \cite{ChPf1, ChPf} (the sweeps correspond to the colours in the hand computation of \cite[Example~3.1]{ChPf}). However, it is much simpler. The candidate basis is fixed in advance by the tree $\cT_o$, so no new vectors are introduced, and by Proposition~\ref{span} none can turn out to be redundant. In particular, we need neither the treatment of redundant cosets nor the handling of coincidences of Linton's enumerator.

		\begin{example}\label{ex:G4table}{\rm
				Let $W=G_4$, $P=\langle s_1\rangle$, and let $\cT_o=(w_1,\dots,w_8)$ be
				the tree of Example \ref{g4bfs} and Figure \ref{fig:one}, so that
				\[
				w_1=e,\ w_2=s_2,\ w_3=s_2s_1,\ w_4=s_2^2,\ w_5=s_2s_1^2,\
				w_6=s_2^2s_1,\ w_7=s_2s_1^2s_2,\ w_8=s_2^2s_1^2.
				\]
				We write the positive Hecke relation as
				$\mathbf s_i^3=a_2\mathbf s_i^2+a_1\mathbf s_i+a_0$, so that the rules
				\textsf{gen} and \textsf{inv} read
				\[
				\mathbf s_i=a_2+a_1\mathbf s_i^{-1}+a_0\mathbf s_i^{-2},\qquad
				\mathbf s_i^{-1}=a_0^{-1}\mathbf s_i^2-a_0^{-1}a_2\mathbf s_i-a_0^{-1}a_1.
				\]
				For $G_4$, the columns indexed by
				$s_1^{\pm1}$ and $s_2^{\pm1}$ can be completed using only \textsf{gen},
				\textsf{inv} and the braid relation, without the columns indexed by
				$z^{\pm1}$; we therefore omit the latter. We list the entries in the
				order in which our programs define them; each entry uses only entries
				defined before it.
				
				\smallskip
				\noindent\emph{Initialisation}. The edges of $\cT_o$ give
				\begin{align*}
					&\cC(w_1,s_2)=\mathbf w_2, && \cC(w_2,s_2^{-1})=\mathbf w_1, &
					&\cC(w_2,s_1)=\mathbf w_3, && \cC(w_3,s_1^{-1})=\mathbf w_2,\\
					&\cC(w_2,s_2)=\mathbf w_4, && \cC(w_4,s_2^{-1})=\mathbf w_2, &
					&\cC(w_3,s_1)=\mathbf w_5, && \cC(w_5,s_1^{-1})=\mathbf w_3,\\
					&\cC(w_4,s_1)=\mathbf w_6, && \cC(w_6,s_1^{-1})=\mathbf w_4, &
					&\cC(w_5,s_2)=\mathbf w_7, && \cC(w_7,s_2^{-1})=\mathbf w_5,\\
					&\cC(w_6,s_1)=\mathbf w_8, && \cC(w_8,s_1^{-1})=\mathbf w_6, &
					&\cC(w_1,s_1)=\mathbf s_1\mathbf w_1, && \cC(w_1,s_1^{-1})=\mathbf s_1^{-1}\mathbf w_1.
				\end{align*}

				\smallskip
				\noindent\emph{First sweep.}
				\begin{enumerate}
					\item $(w_5,s_1)$, \textsf{gen}: since $\mathbf w_5\mathbf s_1^{-1}=\mathbf w_3$
					and $\mathbf w_3\mathbf s_1^{-1}=\mathbf w_2$,
					$\cC(w_5,s_1)=a_2\mathbf w_5+a_1\mathbf w_3+a_0\mathbf w_2.$
					\item $(w_8,s_1)$, \textsf{gen}: since $\mathbf w_8\mathbf s_1^{-1}=\mathbf w_6$
					and $\mathbf w_6\mathbf s_1^{-1}=\mathbf w_4$,
					$\cC(w_8,s_1)=a_2\mathbf w_8+a_1\mathbf w_6+a_0\mathbf w_4.$
					\item $(w_3,s_2)$, \textsf{rel}, using
					$\mathbf s_2=\mathbf s_1^{-1}\mathbf s_2^{-1}\mathbf s_1\mathbf s_2\mathbf s_1$:
					the letters take\\ $\mathbf w_3\mapsto\mathbf w_2\mapsto\mathbf w_1\mapsto
					\mathbf s_1\mathbf w_1\mapsto\mathbf s_1\mathbf w_2\mapsto\mathbf s_1\mathbf w_3$, so
					$\cC(w_3,s_2)=\mathbf s_1\mathbf w_3.$
					\item $(w_4,s_2)$, \textsf{gen}: since $\mathbf w_4\mathbf s_2^{-1}=\mathbf w_2$
					and $\mathbf w_2\mathbf s_2^{-1}=\mathbf w_1$,
					$\cC(w_4,s_2)=a_2\mathbf w_4+a_1\mathbf w_2+a_0\mathbf w_1.$
					\item $(w_6,s_2)$, \textsf{rel}, with the same word: the letters take\\
					$\mathbf w_6\mapsto\mathbf w_4\mapsto\mathbf w_2\mapsto\mathbf w_3\mapsto
					\mathbf s_1\mathbf w_3\mapsto\mathbf s_1\mathbf w_5$, using~(3), so
					$\cC(w_6,s_2)=\mathbf s_1\mathbf w_5.$
					\item $(w_1,s_2^{-1})$, \textsf{inv}: since $\mathbf w_1\mathbf s_2=\mathbf w_2$
					and $\mathbf w_1\mathbf s_2^2=\mathbf w_4$,\\
					$\cC(w_1,s_2^{-1})=a_0^{-1}\mathbf w_4-a_0^{-1}a_2\mathbf w_2-a_0^{-1}a_1\mathbf w_1.$
					\item $(w_3,s_2^{-1})$, \textsf{inv}: since $\mathbf w_3\mathbf s_2=\mathbf s_1\mathbf w_3$
					and $\mathbf w_3\mathbf s_2^2=\mathbf s_1^2\mathbf w_3$,\\
					$\cC(w_3,s_2^{-1})=a_0^{-1}(\mathbf s_1^2-a_2\mathbf s_1-a_1)\mathbf w_3
					=\mathbf s_1^{-1}\mathbf w_3.$
					\item $(w_5,s_2^{-1})$, \textsf{rel}, using
					$\mathbf s_2^{-1}=\mathbf s_1^{-1}\mathbf s_2^{-1}\mathbf s_1^{-1}\mathbf s_2\mathbf s_1$:
					the letters take\\ $\mathbf w_5\mapsto\mathbf w_3\mapsto\mathbf s_1^{-1}\mathbf w_3
					\mapsto\mathbf s_1^{-1}\mathbf w_2\mapsto\mathbf s_1^{-1}\mathbf w_4\mapsto
					\mathbf s_1^{-1}\mathbf w_6$, using~(7), so
					$\cC(w_5,s_2^{-1})=\mathbf s_1^{-1}\mathbf w_6.$
					\item $(w_6,s_2^{-1})$, \textsf{inv}: since $\mathbf w_6\mathbf s_2=\mathbf s_1\mathbf w_5$
					and $\mathbf w_6\mathbf s_2^2=\mathbf s_1\mathbf w_7$,\\
					$\cC(w_6,s_2^{-1})=a_0^{-1}(\mathbf s_1\mathbf w_7-a_2\mathbf s_1\mathbf w_5-a_1\mathbf w_6).$
					\item $(w_2,s_1^{-1})$, \textsf{inv}: since $\mathbf w_2\mathbf s_1=\mathbf w_3$
					and $\mathbf w_2\mathbf s_1^2=\mathbf w_5$,
					$\cC(w_2,s_1^{-1})=a_0^{-1}(\mathbf w_5-a_2\mathbf w_3-a_1\mathbf w_2).$
					\item $(w_4,s_1^{-1})$, \textsf{inv}: since $\mathbf w_4\mathbf s_1=\mathbf w_6$
					and $\mathbf w_4\mathbf s_1^2=\mathbf w_8$,
					$\cC(w_4,s_1^{-1})=a_0^{-1}(\mathbf w_8-a_2\mathbf w_6-a_1\mathbf w_4).$
					\item $(w_7,s_1^{-1})$, \textsf{rel}, using
					$\mathbf s_1^{-1}=\mathbf s_2^{-1}\mathbf s_1^{-1}\mathbf s_2^{-1}\mathbf s_1\mathbf s_2$:
					the letters take\\ $\mathbf w_7\mapsto\mathbf w_5\mapsto\mathbf w_3\mapsto
					\mathbf s_1^{-1}\mathbf w_3\mapsto\mathbf s_1^{-1}\mathbf w_5\mapsto
					\mathbf s_1^{-1}\mathbf w_7$, using~(7), so
					$\cC(w_7,s_1^{-1})=\mathbf s_1^{-1}\mathbf w_7.$
				\end{enumerate}
				
				\smallskip
				\noindent\emph{Second sweep.}
				\begin{enumerate}\setcounter{enumi}{12}
					\item $(w_7,s_1)$, \textsf{gen}: by~(12),
					$\mathbf w_7\mathbf s_1^{-1}=\mathbf s_1^{-1}\mathbf w_7$ and
					$\mathbf w_7\mathbf s_1^{-2}=\mathbf s_1^{-2}\mathbf w_7$, so\\
					$\cC(w_7,s_1)=(a_2+a_1\mathbf s_1^{-1}+a_0\mathbf s_1^{-2})\mathbf w_7
					=\mathbf s_1\mathbf w_7.$
					\item $(w_7,s_2)$, \textsf{gen}: since $\mathbf w_7\mathbf s_2^{-1}=\mathbf w_5$
					and, by~(8), $\mathbf w_5\mathbf s_2^{-1}=\mathbf s_1^{-1}\mathbf w_6$,\\
					$\cC(w_7,s_2)=a_2\mathbf w_7+a_1\mathbf w_5+a_0\mathbf s_1^{-1}\mathbf w_6.$
					\item $(w_8,s_2)$, \textsf{rel}, using
					$\mathbf s_2=\mathbf s_1^{-1}\mathbf s_2\mathbf s_1\mathbf s_2\mathbf s_1^{-1}$:
					we have $\mathbf w_8\mathbf s_1^{-1}=\mathbf w_6$, then, by~(5), (1) and~(3),
					\begin{align*}
						\mathbf w_6\mathbf s_2&=\mathbf s_1\mathbf w_5,\\
						\mathbf s_1\mathbf w_5\mathbf s_1&=a_2\mathbf s_1\mathbf w_5+a_1\mathbf s_1\mathbf w_3+a_0\mathbf s_1\mathbf w_2,\\
						(a_2\mathbf s_1\mathbf w_5+a_1\mathbf s_1\mathbf w_3+a_0\mathbf s_1\mathbf w_2)\mathbf s_2
						&=a_2\mathbf s_1\mathbf w_7+a_1\mathbf s_1^2\mathbf w_3+a_0\mathbf s_1\mathbf w_4,
					\end{align*}
					and finally, multiplying by $\mathbf s_1^{-1}$ using~(12) and~(11),\\
					$\cC(w_8,s_2)=\mathbf s_1\mathbf w_8-a_2\mathbf s_1\mathbf w_6-a_1\mathbf s_1\mathbf w_4
					+a_2\mathbf w_7+a_1\mathbf s_1^2\mathbf w_2.$
					\item $(w_8,s_2^{-1})$, \textsf{inv}: multiplying $\cC(w_8,s_2)$ by
					$\mathbf s_2$ once more gives
					\begin{align*}
						\mathbf w_8\mathbf s_2^2={}&\mathbf s_1^2\mathbf w_8+a_2(a_2+\mathbf s_1)\mathbf w_7
						-a_2(a_1+a_2\mathbf s_1)\mathbf w_6-a_2(\mathbf s_1^2-a_1)\mathbf w_5\\
						&-a_1a_2\mathbf s_1\mathbf w_4+a_1(a_0+a_2\mathbf s_1^2)\mathbf w_2-a_0a_1\mathbf s_1\mathbf w_1,
					\end{align*}
					and hence
					\begin{align*}
						\cC(w_8,s_2^{-1})={}&a_0^{-1}\big(\mathbf w_8\mathbf s_2^2-a_2\mathbf w_8\mathbf s_2-a_1\mathbf w_8\big)\\
						={}&\mathbf s_1^{-1}\mathbf w_8+a_0^{-1}a_2\mathbf s_1\mathbf w_7-a_0^{-1}a_1a_2\mathbf w_6
						-a_0^{-1}a_2(\mathbf s_1^2-a_1)\mathbf w_5+a_1\mathbf w_2-a_1\mathbf s_1\mathbf w_1.
					\end{align*}
				\end{enumerate}
				
				The following table summarises the computation: for each entry, it
				shows the sweep in which it is defined, where $0$ refers to the
				initialisation, and the rule used.
				\[
				\begin{array}{l|cccc}
					& s_1 & s_2 & s_1^{-1} & s_2^{-1}\\\hline
					w_1=e & 0 & 0 & 0 & 1\ \textsf{inv}\\
					w_2=s_2 & 0 & 0 & 1\ \textsf{inv} & 0\\
					w_3=s_2s_1 & 0 & 1\ \textsf{rel} & 0 & 1\ \textsf{inv}\\
					w_4=s_2^2 & 0 & 1\ \textsf{gen} & 1\ \textsf{inv} & 0\\
					w_5=s_2s_1^2 & 1\ \textsf{gen} & 0 & 0 & 1\ \textsf{rel}\\
					w_6=s_2^2s_1 & 0 & 1\ \textsf{rel} & 0 & 1\ \textsf{inv}\\
					w_7=s_2s_1^2s_2 & 2\ \textsf{gen} & 2\ \textsf{gen} & 1\ \textsf{rel} & 0\\
					w_8=s_2^2s_1^2 & 1\ \textsf{gen} & 2\ \textsf{rel} & 0 & 2\ \textsf{inv}
				\end{array}
				\]
				Thus $16$ entries are known after the initialisation, $12$ are defined
				in the first sweep and the remaining $4$ in the second.}
		\end{example}

		Our programs record the sequence of definitions performed, so that the table can be recomputed without any search; these sequences are available on the project's webpage \cite{web}. Whether the table can be completed does not depend on the order in which entries are defined: once the table is complete, $\BFS_o$ is a basis, so all coefficients are uniquely determined, and any entry that can be defined at some stage can still be defined at every later stage.

		\begin{remark}
			In the programs, each entry of the table is
			recorded as a command 
			\begin{center}
				\texttt{setRow(row, column, value, rule)},
			\end{center}
			where
			the column of $s_i^{\pm1}$ is $\pm i$. The initialisation is recorded
			by commands \texttt{makeEdge(x, y, i)}, one for each edge
			$w_x\to w_y$ of $\cT_o$ labelled $s_i$; such a command sets
			$\cC(w_x,s_i)=\mathbf w_y$ and $\cC(w_y,s_i^{-1})=\mathbf w_x$, using a
			rule named \textsf{img}. For instance, the edge $w_2\to w_3$ in
			Example~\ref{ex:G4table} is recorded as \texttt{makeEdge(2, 3, 1)}, which executes \texttt{setRow(2, 1, basicVec(3), img)}
			and \texttt{setRow(3, -1, basicVec(2), img)}. For \textsf{rel}, the third argument is the word used, while for the rules
			\textsf{gen} and \textsf{inv}  the
			relation used is determined by the column, and the third
			argument is $0$.
			Thus, in Example~\ref{ex:G4table}, entry~(1) is recorded as
			\texttt{setRow(5, 1, 0, gen)}, entry~(6) as \texttt{setRow(1, -2, 0, inv)} and entry~(3) as
			$
			\texttt{setRow(3, 2, [-1, -2, 1, 2, 1], rel)},
			$
			where \texttt{[-1, -2, 1, 2, 1]} encodes the word
			$\mathbf s_1^{-1}\mathbf s_2^{-1}\mathbf s_1\mathbf s_2\mathbf s_1$.	
		\end{remark}

		In practice, we proceed in two stages. We first complete the coset table in a
		specialised algebra $\cH_\theta$ where $\theta\colon R(W)\to\bk$ is a ring homomorphism and $\bk$ is a large finite field, as in \S\ref{modtests}. Since we do not have to manipulate polynomials in the parameters, completing the coset table (if possible) is fast. This stage produces a list of $\texttt{setRow}$
		commands that completes the coset table of the candidate. We then replay the same list of commands over $\bar R(W)$. If the
replay succeeds, it completes the coset table over $\bar R(W)$, so
$\mathrm{pr}(\mathrm{BFS}_o)$ spans $\bar{\mathcal H}(W)$. Since the
elements of $\mathrm{BFS}_o$ are braid lifts, hence homogeneous,
$\mathrm{BFS}_o$ is then an $R(W)$-basis of $\mathcal H(W)$ by
Proposition~\ref{prop:lift}(1).		
		
					\begin{figure}[tbp]
		\ifdefined\bfscol\else\newlength{\bfscol}\fi
		\centering
		\footnotesize
		\renewcommand{\arraystretch}{1.15}
		\setlength{\tabcolsep}{4pt}
		\settowidth{\bfscol}{$A=S\setminus\{s_i\}\cup\{z\}$}\setlength{\bfscol}{\dimexpr\bfscol/4-2\tabcolsep+0.5pt\relax}\ifdim\bfscol<1.6em\setlength{\bfscol}{1.6em}\fi
		\begin{adjustbox}{max width=\textwidth}
			\begin{tabular}{c|r|c|r@{${}={}$}l|c||w{c}{\bfscol}w{c}{\bfscol}w{c}{\bfscol}w{c}{\bfscol}||w{c}{\bfscol}w{c}{\bfscol}w{c}{\bfscol}w{c}{\bfscol}||w{c}{\bfscol}w{c}{\bfscol}w{c}{\bfscol}w{c}{\bfscol}}
				\multicolumn{6}{c||}{} & \multicolumn{4}{c||}{$A=S$} & \multicolumn{4}{c||}{$A=S\cup\{z\}$} & \multicolumn{4}{c}{$A=S\setminus\{s_i\}\cup\{z\}$} \\
				$W$ & $|W|$ & $S$ & \multicolumn{2}{c|}{$P$} & $z$ & \rotatebox{90}{orderings} & \rotatebox{90}{candidates} & \rotatebox{90}{complete} & \rotatebox{90}{Modular test 2} & \rotatebox{90}{orderings} & \rotatebox{90}{candidates} & \rotatebox{90}{complete} & \rotatebox{90}{Modular test 2} & \rotatebox{90}{orderings} & \rotatebox{90}{candidates} & \rotatebox{90}{complete} & \rotatebox{90}{Modular test 2} \\
				\hline
				$G_{4}$ & $24$ & $\{s_1,s_2\}$ & \multicolumn{2}{c|}{$\langle s_1\rangle$} & $(s_1s_2)^3$ & $2$ & $1$ & $1$ & $1$ & $6$ & $4$ & $1$ & $1$ & \multicolumn{4}{c}{---} \\
				$G_{5}$ & $72$ & $\{s_1,s_2\}$ & \multicolumn{2}{c|}{$\langle s_1\rangle$} & $(s_1s_2)^2$ & $2$ & $2$ & $2$ & $2$ & $6$ & $4$ & $4$ & $4$ & \multicolumn{4}{c}{---} \\
				$G_{6}$ & $48$ & $\{s_1,s_2\}$ & \multicolumn{2}{c|}{$\langle s_2\rangle$} & $(s_1s_2)^3$ & $2$ & $1$ & $1$ & $1$ & $6$ & $4$ & $4$ & $4$ & \multicolumn{4}{c}{---} \\
				$G_{7}$ & $144$ & $\{s_1,s_2,s_3\}$ & \multicolumn{2}{c|}{$\langle s_2\rangle$} & $s_1s_2s_3$ & $6$ & $5$ & $5$ & $0$ & $24$ & $15$ & $15$ & $0$ & $18$ & $13$ & $13$ & $8$ \\
				$G_{8}$ & $96$ & $\{s_1,s_2\}$ & \multicolumn{2}{c|}{$\langle s_1\rangle$} & $(s_1s_2)^3$ & $2$ & $2$ & $2$ & $2$ & $6$ & $3$ & $3$ & $3$ & \multicolumn{4}{c}{---} \\
				$G_{9}$ & $192$ & $\{s_1,s_2\}$ & \multicolumn{2}{c|}{$\langle s_2\rangle$} & $(s_1s_2)^3$ & $2$ & $2$ & $2$ & $2$ & $6$ & $3$ & $3$ & $3$ & \multicolumn{4}{c}{---} \\
				$G_{10}$ & $288$ & $\{s_1,s_2\}$ & \multicolumn{2}{c|}{$\langle s_1\rangle$} & $(s_1s_2)^2$ & $2$ & $2$ & $2$ & $2$ & $6$ & $4$ & $4$ & $4$ & \multicolumn{4}{c}{---} \\
				$G_{11}$ & $576$ & $\{s_1,s_2,s_3\}$ & \multicolumn{2}{c|}{$\langle s_3\rangle$} & $s_1s_2s_3$ & $6$ & $5$ & $5$ & $0$ & $24$ & $18$ & $18$ & $0$ & $18$ & $15$ & $13$ & $2$ \\
				$G_{12}$ & $48$ & $\{s_1,s_2,s_3\}$ & \multicolumn{2}{c|}{$\langle s_1\rangle$} & $(s_1s_2s_3)^4$ & $6$ & $4$ & $4$ & $4$ & $24$ & $9$ & $0$ & $0$ & \multicolumn{4}{c}{---} \\
				$G_{13}$ & $96$ & $\{s_1,s_2,s_3\}$ & \multicolumn{2}{c|}{$\langle s_1\rangle$} & $(s_1s_2s_3)^3$ & $6$ & $6$ & $4$ & $1$ & $24$ & $12$ & $0$ & $0$ & \multicolumn{4}{c}{---} \\
				$G_{14}$ & $144$ & $\{s_1,s_2\}$ & \multicolumn{2}{c|}{$\langle s_2\rangle$} & $(s_1s_2)^4$ & $2$ & $2$ & $2$ & $2$ & $6$ & $4$ & $1$ & $1$ & \multicolumn{4}{c}{---} \\
				$G_{15}$ & $288$ & $\{s_1,s_2,s_3\}$ & \multicolumn{2}{c|}{$\langle s_2\rangle$} & $s_2s_3s_1s_2s_1$ & $6$ & $6$ & $0$ & $0$ & $24$ & $21$ & $18$ & $3$ & $18$ & $9$ & $4$ & $4$ \\
				$G_{16}$ & $600$ & $\{s_1,s_2\}$ & \multicolumn{2}{c|}{$\langle s_1\rangle$} & $(s_1s_2)^3$ & $2$ & $2$ & $1$ & $0$ & $6$ & $5$ & $2$ & $1$ & \multicolumn{4}{c}{---} \\
				$G_{17}$ & $1200$ & $\{s_1,s_2\}$ & \multicolumn{2}{c|}{$\langle s_2\rangle$} & $(s_1s_2)^3$ & $2$ & $2$ & $1$ & $0$ & $6$ & $5$ & $5$ & $3$ & \multicolumn{4}{c}{---} \\
				$G_{18}$ & $1800$ & $\{s_1,s_2\}$ & \multicolumn{2}{c|}{$\langle s_2\rangle$} & $(s_1s_2)^2$ & $2$ & $2$ & $2$ & $0$ & $6$ & $5$ & $5$ & $0$ & \multicolumn{4}{c}{---} \\
				$G_{19}$ & $3600$ & $\{s_1,s_2,s_3\}$ & \multicolumn{2}{c|}{$\langle s_3\rangle$} & $s_1s_2s_3$ & $6$ & $6$ & $6$ & $0$ & $24$ & $21$ & $21$ & $0$ & $18$ & $16$ & $9$ & $0$ \\
				$G_{20}$ & $360$ & $\{s_1,s_2\}$ & \multicolumn{2}{c|}{$\langle s_1\rangle$} & $(s_1s_2)^5$ & $2$ & $2$ & $2$ & $2$ & $6$ & $4$ & $0$ & $0$ & \multicolumn{4}{c}{---} \\
				$G_{21}$ & $720$ & $\{s_1,s_2\}$ & \multicolumn{2}{c|}{$\langle s_1\rangle$} & $(s_1s_2)^5$ & $2$ & $2$ & $2$ & $1$ & $6$ & $4$ & $0$ & $0$ & \multicolumn{4}{c}{---} \\
				$G_{22}$ & $240$ & $\{s_1,s_2,s_3\}$ & \multicolumn{2}{c|}{$\langle s_1\rangle$} & $(s_1s_2s_3)^5$ & $6$ & $6$ & $1$ & $1$ & $24$ & $14$ & $0$ & $0$ & \multicolumn{4}{c}{---} \\
				$G_{24}$ & $336$ & $\{s_1,s_2,s_3\}$ & \multicolumn{2}{c|}{$\langle s_1\rangle$} & $(s_1s_2s_3)^7$ & $6$ & $6$ & $6$ & $6$ & $24$ & $22$ & $0$ & $0$ & \multicolumn{4}{c}{---} \\
				$G_{25}$ & $648$ & $\{s_1,s_2,s_3\}$ & \multicolumn{2}{c|}{$\langle s_1\rangle$} & $(s_1s_2s_3)^4$ & $6$ & $6$ & $6$ & $6$ & $24$ & $17$ & $0$ & $0$ & \multicolumn{4}{c}{---} \\
				$G_{26}$ & $1296$ & $\{s_1,s_2,s_3\}$ & \multicolumn{2}{c|}{$\langle s_2\rangle$} & $(s_1s_2s_3)^3$ & $6$ & $6$ & $6$ & $6$ & $24$ & $20$ & $2$ & $2$ & \multicolumn{4}{c}{---} \\
				$G_{27}$ & $2160$ & $\{s_1,s_2,s_3\}$ & \multicolumn{2}{c|}{$\langle s_1\rangle$} & $(s_1s_2s_3)^5$ & $6$ & $6$ & $6$ & $0$ & $24$ & $22$ & $0$ & $0$ & \multicolumn{4}{c}{---} \\
				$G_{29}$ & $7680$ & $\{s_1,\dots,s_4\}$ & \multicolumn{2}{c|}{$\langle s_1\rangle$} & $(s_1s_2s_4s_3)^5$ & $24$ & $21$ & $21$ & $13$ & $120$ & $96$ & $0$ & $0$ & \multicolumn{4}{c}{---} \\
				\hline\hline
				$G_{31}$ & $46080$ & $\{s_1,\dots,s_5\}$ & $G(4,2,3)$ & $\langle s_1,s_2,s_3,s_5\rangle$ & $(s_1s_4s_2s_3s_5)^6$ & $120$ & $98$ & $98$ & $38$ & \multicolumn{4}{c||}{---} & \multicolumn{4}{c}{---} \\
				$G_{32}$ & $155520$ & $\{s_1,\dots,s_4\}$ & $G_{25}$ & $\langle s_1,s_2,s_3\rangle$ & $(s_1s_2s_3s_4)^5$ & $24$ & $22$ & $22$ & $0$ & \multicolumn{4}{c||}{---} & \multicolumn{4}{c}{---} \\
				$G_{33}$ & $51840$ & $\{s_1,\dots,s_5\}$ & $G(3,3,4)$ & $\langle s_1,\dots,s_4\rangle$ & $(s_1s_2s_3s_4s_5)^9$ & $120$ & $62$ & $62$ & $62$ & \multicolumn{4}{c||}{---} & \multicolumn{4}{c}{---} \\
				&  &  & $G(3,3,4)$ & $\langle s_2,\dots,s_5\rangle$ &  & $120$ & $62$ & $62$ & $62$ & \multicolumn{4}{c||}{---} & \multicolumn{4}{c}{---} \\
				&  &  & $A_4$ & $\langle s_1,s_2,s_4,s_5\rangle$ &  & $120$ & $76$ & $76$ & $76$ & \multicolumn{4}{c||}{---} & \multicolumn{4}{c}{---} \\
				$G_{34}$\footnotemark & $39191040$ & $\{s_1,\dots,s_6\}$ & $G_{33}$ & $\langle s_1,\dots,s_5\rangle$ & $(s_1s_2s_3s_4s_5s_6)^7$ & $1$ & $1$ & $1$ & $1$ & \multicolumn{4}{c||}{---} & \multicolumn{4}{c}{---}
			\end{tabular}
		\end{adjustbox}
		\caption{$\mathrm{BFS}_o$ candidates for the non-real exceptional complex
			reflection groups. }
		\label{bfs}\label{tbp}
	\end{figure}
	
	\footnotetext{For $G_{34}$  we did not search among candidates:
		we take the basis of $\cH(G_{33})$ and the coset tree of~\cite{MaPf}.}
	
		Figure~\ref{bfs} summarises our computations with breadth-first candidates.
		For the real groups $G_{23}$, $G_{28}$, $G_{30}$, $G_{35}$, $G_{36}$,
		$G_{37}$, every candidate $\mathrm{BFS}_o$ is the standard basis by
		Proposition~\ref{real}, so we only list the non-real groups. For each of them
		we fix a parabolic subgroup $P$ and consider up to three generating
		sets: $A=S$, $A=S\cup\{z\}$ and $A=(S\setminus\{s_i\})\cup\{z\}$ for
		$i=1,2,3$ (for the last choice, the numbers in the figure are summed over
		the three values of $i$). For each generating set, the figure gives the
		number of orderings $o$ of $A$, and the number of distinct candidates
		$\mathrm{BFS}_o$ obtained from those orderings whose tree reaches every
		coset of $P$; different orderings may give the same tree. It then gives
		the number of these candidates whose coset table is completed, and the
		number of the latter that pass Modular test~2 (Lemma~\ref{mtest2}); these are
		the candidates that may be suitable for proving the conjecture. For $G_{31}$,
		$G_{32}$, $G_{33}$ and $G_{34}$ we take for $P$ a maximal parabolic
		subgroup, $A=S$, and $\mathrm{BFS}_o=\cB_P\cdot \cT_o$, where $\cB_P$ is a
		$BFS_o$ basis of $\cH(P)$; for $G_{33}$ we list three choices of $P$.

		The coset tables of the $\mathrm{BFS}_o$ bases, which are available on
the project's webpage \cite{web}, have been completed over
$\bar R(W)$; as explained above, this establishes that these sets are
$R(W)$-bases. We deduce the following:
		
		\begin{theorem}\label{thm-BFSbasis}
			For every exceptional complex reflection group $W$, there exists a generating set $A \subseteq \Sigma$ and an ordering $o$ on $A$ such that $\BFS_o$ is an $R(W)$-basis of $\cH(W)$.
		\end{theorem}

		\begin{remark}{\rm
				Theorem~\ref{thm-BFSbasis} gives a new proof of the BMR freeness
				conjecture for the exceptional complex reflection groups. Previous
				proofs treated the groups case by case, or uniformly only within a
				family, such as the groups of rank~$2$ in \cite{Tsu} or the $2$-reflection
				groups in \cite{MaPf}. Our proof applies the same construction, and the
				same verification procedure, to all exceptional groups; it provides explicit bases for which every element of $\cH(W)$
				can be expressed in the basis by means of the coset table. }
		\end{remark}

	\section{$z^k$-bases and the BMM conjecture   for $G_{4},\ldots,G_{29}$}\label{bases}
		Now that we have obtained $\BFS_o$ bases for every exceptional
		complex reflection group $W$, recall that our strategy for proving the BMM symmetrising
		trace conjecture is to find a section
		basis $\cB$ of $\cH(W)$  such that $\lambda^{\cB}_1$ is the
		canonical symmetrising trace of $\cH(W)$. Modular test~2 (Lemma~\ref{mtest2}) allows us to identify bases that violate this condition. For example, Figure~\ref{bfs} shows that none of the $BFS_o$ bases of $G_{18}$ or $G_{19}$ satisfies the test.
		Moreover, even for bases that pass this test, the Gram
		matrix is in general too complicated to be computed directly, let alone calculate its determinant. We therefore look for basis candidates with
		additional structure that simplifies these computations, and we use two
		approaches, depending on the group.

For $G_4,\ldots,G_{29}$, except for the real reflection groups $G_{23}$ and $G_{28}$ (for which we already have the validity of the BMM symmetrising trace conjecture with the standard basis), we use a generalisation of \emph{$z$-bases}, which we call \emph{$z^k$-bases}. The name $z$-bases was given by the first and third authors in~\cite[Definition~2.4]{ChPf}. Using
the completed coset tables of the $\mathrm{BFS}_o$ bases, we find for
each of these groups a $z^k$-basis candidate, and we show that it is $\mathrm{pr}$-admissible (Definition~\ref{def:pradm}),
so that the conjecture over $R(W)$ follows from
Theorem~\ref{thm:liftall}. Although the conjecture is
already known for $G_4,\dots,G_{15}$, we treat these groups as well, so
that all the groups in this section are handled by the same method; the
$z$-bases we obtain for them differ from those constructed previously.
For $G_{31},\dots,G_{34}$, this approach requires too complicated
computations, because of the large orders of the groups, and in
Section~\ref{rem} we use instead a tower of parabolic subgroups.

To construct $z^k$-basis candidates, we will need the whole coset graph,
and not simply a coset tree, as we did for the $\mathrm{BFS}_o$ bases.  As in the previous section, let $W$ be a complex reflection group, $U$ a subgroup of $W$,
	and $A$ a finite generating set for $W$, not necessarily minimal and chosen
	independently of the generating sets used in
	Section~\ref{standard}.

\begin{definition}\label{def:cosetgraph}{\rm
The \emph{coset graph} $\Gamma_W(A,U)$ is the
directed multigraph whose vertices are the right cosets of $U$ in $W$,
with an edge labelled $a\in A$ from $Ux$ to $Uxa$ for every right coset
$Ux$ and every $a\in A$.}
\end{definition}

Several edges may join the same two vertices, when $Uxa=Uxa'$ for
distinct letters $a,a'\in A$, and an edge may be a loop, when $Uxa=Ux$.
We keep all of them, since the label of an edge becomes a letter of a
word, and distinct letters give distinct elements of $\cH(W)$.

\begin{example}\label{graph:G4}
	Consider the group $G_4=\langle s_1,s_2 \,|\, s_1s_2s_1=s_2s_1s_2,\, s_1^3=s_2^3=1\rangle$ and the parabolic subgroup
	$P=\langle s_1\rangle=\{1,s_1,s_1^2\}$,
	as in Example~\ref{g4bfs}.
	We take $U=P$ and $A=\{s_1,s_2\}$. The coset graph
	$\Gamma_{G_4}(A,P)$ has eight vertices, corresponding
	to the right cosets of $P$ in $G_4$.
	These cosets are represented by the words
	$w_1,\ldots,w_8$ of the coset tree
	$\mathcal T_o$ constructed in
	Example~\ref{g4bfs}.
	The edges describe right multiplication by
	$s_1$ and $s_2$ on these cosets.
	Figure~\ref{fig:cosetgraph_g4} shows the
	resulting coset graph.
	Edges labelled by $s_1$ are shown in blue,
	while those labelled by $s_2$ are shown in red.
	The edges corresponding to $\mathcal T_o$
	are highlighted by solid lines, while the
	remaining edges are shown as dashed lines. \qed
\end{example}

\begin{figure}[h]
	\centering
	\definecolor{cone}{RGB}{0,90,181}
	\definecolor{ctwo}{RGB}{204,51,17}
	
	\begin{tikzpicture}[scale=1.25,
		v/.style={circle,draw,fill=white,
			inner sep=1.2pt,font=\small},
		tree/.style={line width=1.5pt,-stealth},
		rest/.style={line width=0.5pt,
			densely dashed,-stealth}]
		
		\node[v] (n2) at (0,1.5)       {$w_{2}$};
		\node[v] (n4) at (1.5,0)       {$w_{4}$};
		\node[v] (n6) at (0,-1.5)      {$w_{6}$};
		\node[v] (n5) at (-1.5,0)      {$w_{5}$};
		\node[v] (n1) at (1.75,1.75)   {$w_{1}$};
		\node[v] (n3) at (-1.75,1.75)  {$w_{3}$};
		\node[v] (n7) at (-1.75,-1.75) {$w_{7}$};
		\node[v] (n8) at (1.75,-1.75)  {$w_{8}$};
		
		\draw[tree,ctwo] (n1) -- (n2);
		\draw[tree,cone] (n2) -- (n3);
		\draw[tree,ctwo] (n2) -- (n4);
		\draw[tree,cone] (n3) -- (n5);
		\draw[tree,cone] (n4) -- (n6);
		\draw[tree,ctwo] (n5) -- (n7);
		\draw[tree,cone] (n6) -- (n8);
		
		\draw[rest,ctwo] (n4) -- (n1);
		\draw[rest,cone] (n5) -- (n2);
		\draw[rest,cone] (n8) -- (n4);
		\draw[rest,ctwo] (n7) -- (n6);
		\draw[rest,ctwo] (n6) -- (n5);
		
		\draw[rest,cone] (n1.north)
		.. controls +(60:0.55) and +(30:0.55)
		.. (n1.east);
		
		\draw[rest,ctwo] (n3.west)
		.. controls +(150:0.55) and +(120:0.55)
		.. (n3.north);
		
		\draw[rest,cone] (n7.south)
		.. controls +(240:0.55) and +(210:0.55)
		.. (n7.west);
		
		\draw[rest,ctwo] (n8.east)
		.. controls +(330:0.55) and +(300:0.55)
		.. (n8.south);
		
		\begin{scope}[shift={(3.0,0.2)},font=\scriptsize]
			\draw[tree,black] (0,0) -- (0.5,0);
			\node[right] at (0.55,0)
			{Edges of $\mathcal T_o$};
			
			\draw[rest,black] (0,-0.4) -- (0.5,-0.4);
			\node[right] at (0.55,-0.4)
			{Remaining edges};
			
			\draw[cone,line width=1pt]
			(0,-0.85) -- (0.5,-0.85);
			\node[right] at (0.55,-0.85) {$s_1$};
			
			\draw[ctwo,line width=1pt]
			(0,-1.2) -- (0.5,-1.2);
			\node[right] at (0.55,-1.2) {$s_2$};
		\end{scope}
	\end{tikzpicture}
	
	\caption{The coset graph $\Gamma_{G_4}(A,P)$ for $P=\langle s_1\rangle$ and $A=\{s_1,s_2\}$.}
	\label{fig:cosetgraph_g4}
\end{figure}

As in Example \ref{graph:G4}, every coset tree with respect to $(A,U)$ determines a spanning tree of the coset graph, and conversely. More precisely, let $(w_1,\dots,w_{[W:U]})$ be a coset tree
with respect to $(A,U)$. For each $i>1$, writing $w_i$ as the word
$w_ja$ with $j<i$ and $a\in A$, the coset graph $\Gamma_W(A,U)$ has an edge
labelled $a$ from $Uw_j$ to $Uw_i$ (this edge is unique, since each vertex has exactly one outgoing edge with each label). These edges form a spanning tree of
$\Gamma_W(A,U)$, rooted at $U$ and directed away from it. Conversely, given
such a spanning tree, the labels along the path from $U$ to a vertex
$v$ form a word whose image in $W$ lies in $v$, since the path follows
right multiplication by its labels; and these words, listed so that
each comes after its parent, form a coset tree. Regarding coset trees
that differ only in the order of their words as the same, we obtain a
bijection between $\cM^U_A$ and the set of spanning trees of $\Gamma_W(A,U)$
rooted at $U$ and directed away from it. Since $\Gamma_W(A,U)$ is finite,
there are finitely many such trees, and we obtain every element of
$\cM^U_A$ by enumerating them.

	\subsection{Finding $z^k$-bases}
	Let $W$ be an irreducible complex reflection group. 
	Recall that we denote by $z\in Z(W)$  the image of the generator $\mathbf z$ of $Z(B(W))$.
	Let $l$ be the order of  $z$ in $W$, and let $k$ be a positive divisor of $l$ so that $z^k$ has
order $d:=l/k$. We set $m_k:=|W|/d$.
	
	In the previous section, in order to construct $\BFS_o$ bases, we took $U$ to be a parabolic subgroup $P$ of $W$.
In order to construct $z^k$-bases, we take $U$ to be $W_k:=PZ^k$, where $P=\langle s \rangle$ for some $s \in S$ and $Z^k:=\langle z^k\rangle$. 	If we denote by $e_s$ the order of $s$, then $W_k$ is a subgroup of $W$ of order $e_sd$ since $P \cap Z^k = \{1\}$ (every non-trivial power of $s$ is a pseudo-reflection and fixes a hyperplane pointwise, whereas every non-trivial element of $Z(W)$ acts on $V$ as a non-trivial scalar); in particular $m_k/e_s=[W:W_k]$ is an integer.

\begin{definition}\label{zk-basis}{\rm 
	An $R(W)$-basis $\mathcal{B}$ of $\cH(W)$ is called
	a \emph{$z^k$-basis} if there exist $s\in S$ and 	$\mathbf b_1,\ldots,\mathbf b_{m_k/e_s}\in B(W)$, with $\mathbf b_1=1$,
	such that
	\[
	\mathcal B
	=
	\{1,\mathbf z^{k},\ldots,\mathbf z^{k(d-1)}\}
	\cdot
	\{1,\mathbf s,\ldots,\mathbf s^{\,e_s-1}\}
	\cdot
	\{\mathbf b_1,\ldots,\mathbf b_{m_k/e_s}\}.
	\]
	When $k=1$ we simply call it a
	\emph{$z$-basis}.}
\end{definition}

Bases of this form first appeared in \cite{Ch17}, where the first author constructed
them for the groups $G_4,\dots,G_{15}$. The notion of a $z$-basis was introduced by the first and third
authors in \cite[Definition~2.4]{ChPf}, who constructed different $z$-bases for these
groups, and also for $G_{22}$ in
\cite{ChPf1}.

 A main advantage of these bases is their structure, which
	makes Gram matrices easier to compute. In particular,
	\cite[Lemma~3.2]{ChPf} gives a recursive formula that computes each
	block of the Gram matrix from the preceding ones. Moreover, for such bases Condition~(3) is easy to check: since they are section bases this condition is
	equivalent to \eqref{extra con2}, and the elements $\mathbf w^{-1}\boldsymbol\pi$
	occurring there are, thanks to the central factor $\mathbf z$, easily expressed in terms of
	the elements of the basis, so that the vanishing of $\lambda_1^{\cB}$ on them is immediate.
	
	Our goal was to find such bases for the remaining exceptional groups.
	However, for $G_{16}$ none of the $z$-basis candidates in our search
	is a basis (see Example~\ref{G16}). This led us to generalise the
	notion of a $z$-basis: we allow a larger class of candidates, but keep
	the structure that makes Gram matrices easy to compute.

We now turn to the construction of candidates for $z^k$-bases. If $\mathbf z^k$ and $\mathbf s$ are fixed, it remains to determine suitable elements $\mathbf b_1,\ldots,\mathbf b_{m_k/e_s}$. Under the specialisation $a_{s,0}\mapsto1$, $a_{s,i}\mapsto0$, a basis of $\cH(W)$ becomes a basis of the group algebra, so the $|W|$ elements $z^{k\alpha}s^{\beta}b_i$ of $W$, with $0\le\alpha<d$ and $0\le\beta<e_s$, where $b_i$ denotes the image of $\mathbf b_i$ in $W$, must be pairwise distinct. Since the $z^{k\alpha}s^{\beta}$ form the subgroup $W_k$, the $b_i$ must represent the right cosets of $W_k$ in $W$.

\begin{definition}\label{def:zkcandidate}{\rm
	Let $s\in S$, $P=\langle s\rangle$ and $k$ a positive divisor of $l$. For a coset tree
	$\mathcal T=\{w_1,\ldots,w_{m_k/e_s}\}\in\mathcal M_A^{W_k}$, let
	$\mathbf w_i\in B(W)$ be obtained from the word $w_i$ by replacing
	each letter of $A$ by its fixed lift to $B(W)$. We call the
	collection
	\[
	\{1,\mathbf z^{k},\ldots,\mathbf z^{k(d-1)}\}
	\cdot
	\{1,\mathbf s,\ldots,\mathbf s^{\,e_s-1}\}
	\cdot
	\{\mathbf w_1,\ldots,\mathbf w_{m_k/e_s}\}
	\]
	a \emph{$z^k$-basis candidate}. Its elements are braid lifts of the $|W|$ pairwise distinct elements of
$W$, so a $z^k$-basis candidate that is a basis is a section basis. When $k=1$ we simply call it a
	\emph{$z$-basis candidate}.}
\end{definition}

\begin{remark}
When $P=\langle s\rangle$ and $\cB_P=\{1,\mathbf s,\dots,\mathbf s^{e_s-1}\}$, a $\mathrm{BFS}_o$ basis candidate has the same form as a $z^k$-basis candidate for $k=l$: it is $\{1,\mathbf s,\dots,\mathbf s^{e_s-1}\}\cdot \cT_o$, where $\cT_o$ is a coset tree with respect to $(A,P)$. The $z^k$-basis candidates generalise this by adding the central factor and allowing arbitrary coset trees.
\end{remark}

By the correspondence between coset trees and spanning trees, we obtain
all coset trees in $\mathcal M^{W_k}_A$ by enumerating the spanning
trees of $\Gamma_W(A,W_k)$ rooted at $W_k$ and directed away from it, and
hence all $z^k$-basis candidates for the given $s$ and $k$. Among these
we discard those failing Modular test~1 (Lemma~\ref{mtest}) or Modular test~2
(Lemma~\ref{mtest2}). Both tests concern a $z^k$-basis candidate $\cB$ and are
carried out in the specialised algebra, using as reference basis a
$\BFS_o$ basis $\mathcal A$ of $\cH(W)$ with a completed coset
table. The coset table of $\mathcal A$ gives the products of elements
of $B$ in terms of $\mathcal A$; for test~1 we compute the
change-of-basis matrix from $\mathcal A$ to $\cB$, and for test~2 we use
its inverse to express these products in terms of $\cB$ and compute the
Gram matrix with respect to $\lambda^{\cB}_1$.

The following example shows what happens for $G_{16}$ when $k=1$, and
why values of $k>1$ are needed.

\begin{example}\label{G16}{\rm
Let $W=G_{16}=\langle s_1,s_2\mid s_1^5=s_2^5=1,\
s_1s_2s_1=s_2s_1s_2\rangle$, and take $s=s_2$ and $k=1$, so that
$P=\langle s_2\rangle=\{1,s_2,s_2^2,s_2^3,s_2^4\}$ and $W_1=\langle s_2\rangle Z(W)$, a subgroup of
order $50$ and index $12$ in $W$. Set $w=s_1s_2s_1$ and
$A=\{s_1^{\pm1},s_2^{\pm1},w\}$, where each element of $A$ is lifted to
$B(W)$ letter by letter (so $w$ lifts to
$\mathbf s_1\mathbf s_2\mathbf s_1$). For this choice,
$|\mathcal M^{W_1}_A|=405504$. We applied Modular test~1 to the
corresponding $z$-basis candidates, using a $\mathrm{BFS}_o$ basis of
$\cH(G_{16})$ with a completed coset table as reference basis, and all of
them fail the test. Hence none of them is an $R(W)$-basis of $\cH(W)$.}
\end{example}

The values of $k$ for which we found candidates that pass the modular
tests  (and turn out to be bases) are the following: For $G_{16}$ we found a $z^2$-basis candidate.
For $G_{17}$, $G_{18}$ and $G_{19}$ we found candidates with $k=4,6,12$
respectively, using the normal subgroup relations to guide the search,
as explained below. For $G_4,\dots,G_{15}$ and $G_{20},\ldots,G_{29}$ we found $z$-basis
candidates.

The group $G_{16}$ is a normal reflection subgroup of each of
$G_{17}$, $G_{18}$ and $G_{19}$, of index $i=2,3,6$ respectively.
Since the BMR freeness conjecture holds, a specialisation of the
parameters identifies $\mathcal H(G_{16})$ with a subalgebra of the
specialised Hecke algebra of the larger group (see \cite[Proposition 4.2]{Ma2} and
\cite[Section~3.2]{ChPf}). Moreover, if $z$ generates the centre of the
larger group, then $z^{i}$ generates $Z(G_{16})$ (see
\cite[Section~1]{ChPf}). The $z^2$-basis candidate for $G_{16}$
therefore suggests looking for $z^k$-basis candidates with $k=2i$, so
that $z^k=(z^{i})^2$ and the central factors of the two candidates
coincide. This gives $k=4$, $6$ and $12$ respectively.

The $z^k$-basis candidates we use, together with the programs used to find them are available on the
project's webpage~\cite{web}. In the rest of this section, we show that they are
$\mathrm{pr}$-admissible, and deduce the conjecture for these groups.

\subsection{Criteria for bases and symmetrising traces}\label{main1}
The results of this subsection rely only on the fact that $\cH(W)$ is free of rank $|W|$ over $R(W)$ and generated as an algebra by $
\mathbf S$. They therefore remain valid, with the same proofs, for any specialised algebra $\mathcal H_\theta$ (see \S\ref{sec-2.1}).

The following proposition shows that, to establish  that a linear form is a trace function, it suffices to check the trace condition on a spanning set for each standard generator. In this section, we use the same basis for every generator, but we state the result in this more general form since it will be needed in Section~\ref{rem}.

\begin{proposition}\label{symcond}
	Let $t\in\Hom_{R(W)}(\cH(W),R(W))$ be a linear form. Suppose that for every
	$\mathbf s\in\mathbf S$ there is a set $Y_{\mathbf s}$ that spans $\cH(W)$ over
	$R(W)$ and satisfies $t(y\mathbf s)=t(\mathbf sy)$ for all $y\in Y_{\mathbf s}$.
	Then $t$ is a trace function, that is, $t(hh')=t(h'h)$ for all $h,h'\in\cH(W)$.
	In particular, this holds if $\cB$ is an $R(W)$-basis of $\cH(W)$ and
	$t(b\mathbf s)=t(\mathbf sb)$ for all $b\in\cB$ and $\mathbf s\in\mathbf S$.
\end{proposition}

\begin{proof}
	Fix $\mathbf s\in\mathbf S$. Since $Y_{\mathbf s}$ spans $\cH(W)$ and the map
	$h\mapsto t(h\mathbf s)-t(\mathbf sh)$ is $R(W)$-linear, the hypothesis implies
	that $t(h\mathbf s)=t(\mathbf sh)$ for every $h\in\cH(W)$. This holds for every
	$\mathbf s\in\mathbf S$. Since $\mathbf S$ generates $\cH(W)$ as an
	$R(W)$-algebra, it suffices to extend this equality to products of generators.
	Let $\mathbf s_1,\ldots,\mathbf s_r\in\mathbf S$. Applying the equality
	successively, we obtain
	$$t(h\mathbf s_1\cdots \mathbf s_r)
	=t(\mathbf s_rh\mathbf s_1\cdots \mathbf s_{r-1})
	=t(\mathbf s_{r-1}\mathbf s_rh\mathbf s_1\cdots \mathbf s_{r-2})
	=\cdots
	=t(\mathbf s_1\cdots \mathbf s_rh).$$
	Hence $t(hx)=t(xh)$ for every $h\in\cH(W)$ and every product $x=\mathbf s_1\cdots\mathbf s_r$ of generators. Since $\mathbf S$ generates $\mathcal H(W)$ as an $R(W)$-algebra, such products span $\mathcal H(W)$ over $R(W)$, and by $R(W)$-linearity $t(hx)=t(xh)$ for all $h,x\in\mathcal H(W)$.
\end{proof}
We next examine how Gram matrices transform under a change of basis. This allows us to study the Gram matrix of a basis candidate by expressing its elements in terms of a reference basis. All notation used below has been introduced in \S\ref{modtests}.

\begin{proposition}\label{cbmg}
	Let $t\in\Hom_{R(W)}(\cH(W),R(W))$, let $\cB$ be an $R(W)$-basis of $\cH(W)$, and
	let $X,Y$ be two ordered collections of $|W|$ elements of $\cH(W)$. Then
	\[
	\gram_t(X,Y)=P_{\cB}(X)^T\gram_t(\cB,Y)=\gram_t(X,\cB)P_{\cB}(Y).
	\]
	In particular,
	$\gram_t(X,Y)=P_{\cB}(X)^T\gram_t(\cB,\cB)P_{\cB}(Y).$
\end{proposition}

\begin{proof}
Write $\cB=(b_1,\ldots,b_n)$, $X=(x_1,\ldots,x_n)$, and $Y=(y_1,\ldots,y_n)$, where $n=|W|$. By definition,
\[
x_i=\sum_{r=1}^n \bigl(P_{\cB}(X)\bigr)_{r,i}b_r.
\]
Thus, for all $i,j$,
$$
	t(x_i y_j)
	=\sum_{r=1}^n \bigl(P_{\cB}(X)\bigr)_{r,i}t(b_r y_j)\\
	=\bigl(P_{\cB}(X)^T\gram_t(\cB,Y)\bigr)_{i,j}.
$$
This establishes the first equality. The second follows by the same argument, expanding $y_j$ in the basis $\cB$. Applying the second equality with $X=\mathcal B$ gives $\mathrm{Gram}_t(\mathcal B,Y)=\mathrm{Gram}_t(\mathcal B,\mathcal B)P_{\mathcal B}(Y)$, and substituting this into the first yields the final formula.
\end{proof}

The following criterion allows us to establish simultaneously that two collections are bases and that a trace function is a symmetrising trace.

\begin{lemma}[Pairing lemma]\label{pair}
	Let $t\in\Hom_{R(W)}(\cH(W),R(W))$ be a trace function, and let $X,Y$ be two ordered collections of $|W|$ elements in $\cH(W)$.
	If $\det(\gram_t(X,Y))$ is invertible in $R(W)$, then
	\begin{enumerate}
		\item $X$ and $Y$ are $R(W)$-bases of $\cH(W)$.\smallbreak
		\item $t$ is a symmetrising trace.
	\end{enumerate}
\end{lemma}

\begin{proof}
	Let $\cB$ be a fixed $R(W)$-basis of $\cH(W)$. By Proposition~\ref{cbmg},
	$$
	\gram_t(X,Y)
	=
	P_{\cB}(X)^T
	\gram_t(\cB,\cB)
	P_{\cB}(Y).
	$$
	Taking determinants, we obtain
	$
	\det(\gram_t(X,Y))
	=
	\det(P_{\cB}(X))
	\det(\gram_t(\cB,\cB))
	\det(P_{\cB}(Y)).
	$
	Since the left-hand side is invertible in $R(W)$, each factor on the right-hand side is invertible. In particular, both $P_{\cB}(X)$ and $P_{\cB}(Y)$ are invertible over $R(W)$. Thus, $X$ and $Y$ are $R(W)$-bases of $\cH(W)$.
	Moreover, $\det(\gram_t(\cB,\cB))$ is invertible, so the bilinear form
	$(h,h')\longmapsto t(hh')$ is non-degenerate. Since $t$ is a trace function by assumption, it is a symmetrising trace.
\end{proof}

To prove that a $z^k$-basis candidate $\cB$ that passes the modular
tests is a basis, we apply the pairing lemma (Lemma \ref{pair}) to $\cB$ and
$\cB^{-1}$. This requires the determinant of the $|W|\times|W|$ matrix
$\mathrm{Gram}_t(\cB,\cB^{-1})$, which in general is too large to compute.
Since $\mathbf z$ is central, however, this matrix has a block
structure, and when some of the blocks vanish, its determinant reduces
to that of a matrix $d$ times smaller (where $d=l/k$).

We have $\cB=\mathbf Z^k\cdot Y$,
where $\mathbf Z^k=\{1,\mathbf z^k,\dots,\mathbf z^{k(d-1)}\}$ and
$Y=\{1,\mathbf s,\dots,\mathbf s^{e_s-1}\}\cdot\{\mathbf w_1,\dots,
\mathbf w_{m_k/e_s}\}$ consists of $N=m_k=|W|/d$ elements, which we
denote by $y_1,\dots,y_N$. We order the elements of $\cB$ by the powers
of $\mathbf z^k$, so that its elements are $\mathbf z^{k\alpha}y_i$,
with $0\le\alpha<d$ and $1\le i\le N$. The following lemma is
similar to \cite[\S 3.3.1]{ChPf}.

\begin{lemma}\label{gramdec}
	Let $t\in\Hom_{R(W)}(\cH(W),R(W))$ and, for each integer $r$, define the
	$N\times N$ matrix
	\[
	F_r=\big(t(\mathbf z^{kr}y_iy_j^{-1})\big)_{i,j=1}^N.
	\]
	Then
	\[
	\gram_{t}(\cB,\cB^{-1})=
	\begin{pmatrix}
		F_0 & F_{-1} & \cdots & F_{-(d-1)}\\
		F_1 & F_0 & \cdots & F_{-(d-2)}\\
		\vdots & \vdots & \ddots & \vdots\\
		F_{d-1} & F_{d-2} & \cdots & F_0
	\end{pmatrix}.
	\]
	In particular, if $F_r=0$ for all $0<r<d$, the Gram matrix is upper block
	triangular, and
	\[
	\det\bigl(\gram_{t}(\cB,\cB^{-1})\bigr)=(\det F_0)^d.
	\]
\end{lemma}

\begin{proof}
Since $\mathbf z$ is central, for $0\le\alpha,\beta<d$ we have
\[
t\big((\mathbf z^{k\alpha}y_i)(\mathbf z^{k\beta}y_j)^{-1}\big)
=t\big(\mathbf z^{k\alpha}y_iy_j^{-1}\mathbf z^{-k\beta}\big)
=t\big(\mathbf z^{k(\alpha-\beta)}y_iy_j^{-1}\big).
\]
	Thus, the block in position $(\alpha,\beta)$ is $F_{\alpha-\beta}$, giving the stated decomposition.
	
	If $F_r=0$ for all $0<r<d$, every block below the diagonal vanishes. The matrix is
	therefore upper block triangular with $d$ diagonal blocks equal to $F_0$, and its
	determinant is $(\det F_0)^d$.
\end{proof}

\subsection{Proof of the conjecture for $G_4,\dots,G_{29}$}\label{verif}
Let $W=G_j$, with $4\leq j\leq29$ and $j\neq23,28$. For each of these groups, we
have constructed a $\BFS_o$ basis $\cA_j$ of $\cH(W)$ with a completed coset table (see Section \ref{standard}).
 For each such group, we also have a $z^k$-basis candidate that passes
the modular tests, which we denote by $\cB_j$. These collections are available on the project's webpage \cite{web}.

Our goals are to prove that the collections $\cB_j$ are bases of $\cH(W)$ and that $\lambda_1^{\cB_j}$ is the canonical symmetrising trace. A direct approach to the first goal would require showing that the determinant of the $|W|\times|W|$ change-of-basis matrix from $\cA_j$ to $\cB_j$ is a unit in $R(W)$, which is computationally infeasible. Instead, we use the pairing lemma (Lemma~\ref{pair}), which allows us to establish both that $\cB_j$ is a basis and that a  trace function is a symmetrising trace. Moreover, the $z^k$-structure of $\cB_j$ (Lemma~\ref{gramdec}) reduces the computation to the determinant of a much smaller matrix.

To apply the pairing lemma, we first need a trace function. Since $\lambda_1^{\cB_j}$ is only defined once $\cB_j$ is known to be a basis, we instead work with a linear form $E_j$ defined using the $\mathrm{BFS}_o$ basis $\cA_j$, which can be evaluated via its completed coset table.  We prove that, after the specialisation $\mathrm{pr}$ of
\S\ref{sec-pr}, $E_j$ becomes a trace function equal to
$\lambda_1^{\mathrm{pr}(\mathcal B_j)}$ (Proposition~\ref{barconj});
by Theorem~\ref{thm:liftall}, this suffices for the conjecture over
$R(W)$.

Our candidates for $E_j$ are of two kinds. If $\mathcal A_j$ passes
Modular test~2 (Lemma~\ref{mtest2}), then $\lambda_1^{\mathcal A_j}$ is
a natural candidate, since its Gram matrix with respect to
$\mathcal A_j$ is then symmetric modulo~$p$. We take $E_j=\lambda_1^{\mathcal A_j}$ for all $W=G_j$ with $4\le j\le29$ and $j\neq23,28$, except for $W=G_{17},G_{18},G_{19},G_{27}$. Here $\mathcal A_j$ is the $\mathrm{BFS}_o$ basis of the corresponding row of Figure~\ref{bfs}: for $G_7$, $G_{11}$ and $G_{15}$ it is obtained from a generating set of the form $(S\setminus\{s_i\})\cup\{z\}$, for $G_{16}$ and $G_{19}$ from $A=S\cup\{z\}$, and for the rest of them from $A=S$.

For the remaining groups, we take $E_j$ of the form
$h\mapsto\lambda_a^{\mathcal A_j}(\mathbf z^{-kJ}h)$, where $a\in\mathcal A_j$ and
$J\in\mathbb Z$ are chosen, by modular computations, so that this
linear form is a trace function modulo~$p$; the choices are listed in
Table~\ref{tab:traces}. Although for $G_{17}$ some of the bases $\mathrm{BFS}_o$ with $A=S\cup\{z\}$ also pass
Modular test~2 (see Figure~\ref{bfs}), the choice of $E_{17}$ above, with $A=S$, leads to
considerably faster computations.
 In both cases, that $E_j$ is a trace function is not assumed: we verify it after specialisation in Proposition~\ref{barconj}, which suffices by Theorem~\ref{thm:liftall}.

We now explain where the second shape comes from; this is only a
heuristic.
Suppose that $\tau$ is a symmetrising trace. For an element $a\in\mathcal A_j$, let $a^\vee$ be the
element of the dual basis of $\mathcal A_j$ with respect to $\tau$,
that is, $\tau(a^\vee a')=\delta_{a,a'}$ for $a'\in\mathcal A_j$. Expanding an element $h\in\mathcal \cH(W)$ in the basis $\mathcal A_j$, we obtain $
\tau(a^\vee h)=\lambda_a^{\mathcal A_j}(h).
$
If $\lambda_a^{\mathcal A_j}$ is a trace function, then for all $h,h'\in\mathcal H(W)$ we have:
\[
\tau(a^\vee hh')=\lambda_a^{\mathcal A_j}(hh')=\lambda_a^{\mathcal A_j}(h'h)=\tau(a^\vee h'h)=\tau(ha^\vee h'),
\]
hence $\tau\bigl((a^\vee h-ha^\vee)h'\bigr)=0$ for all $h'\in\mathcal H(W)$.
Since $\tau$ is non-degenerate, we have that $a^\vee$ is central. It is natural to look for
elements $a\in\mathcal A_j$ for which $a^\vee$ is a power of $\mathbf z^k$, where $k$ is the
integer for which $\mathcal B_j$ is a $z^k$-basis candidate,
so that the resulting
linear form is compatible with the $z^k$-structure of $\mathcal B_j$ (Lemma~\ref{gramdec}).
If $a^\vee=\mathbf z^{kJ}$ for some $J\in\mathbb Z$, then $\lambda_a^{\mathcal A_j}(h)=\tau(\mathbf z^{kJ}h)$
for all $h\in\mathcal H(W)$. Replacing $h$ by $\mathbf z^{-kJ}h$, we obtain
\[
\tau(h)=\lambda_a^{\mathcal A_j}(\mathbf z^{-kJ}h).
\]
In particular, the bilinear form $(h,h')\mapsto\tau(hh')=\lambda_a^{\mathcal A_j}(\mathbf z^{-kJ}hh')$ is symmetric.
This motivates the choice of $E_j$ described above.

 Recall that $\mathcal A_j=\mathcal B_{\langle u \rangle}\cdot \cT_o$ is the
BFS$_o$ basis of Section~\ref{standard}, with respect to
a cyclic parabolic subgroup $\langle u \rangle$ of $W$ and $\cT_o=(w_1,w_2,\dots,w_{[W:\langle u \rangle]})$ is the associated
coset tree.
In Table \ref{tab:traces}, the pair $(J,a)$ is encoded by a triple $(J,i,\eta)$, where $a={\bf u}^{\eta}{\bw}_i$.
We define the $R(W)$-linear form
\[
E_j(h)=\lambda^{\mathcal A_j}_{{\bf u}^{\eta}{\bw}_i}\bigl(\mathbf z^{-kJ}h\bigr),\qquad \text{for all }h\in\mathcal H(W).
\]
For the groups not appearing in the table, if we take $(J,i,\eta)=(0,1,0)$, we recover $E_j=\lambda^{\mathcal A_j}_1$.
\begin{table}[H]
	\centering
	\renewcommand{\arraystretch}{1.3}
	\begin{adjustbox}{max width=\textwidth}
		\begin{tabular}{c|c|c|c|l}
			$G_j$ & $u$ & $k$ & $(J,i,\eta)$ & $a={\bf u}^{\eta}{\bw}_i$ \\
			\hline
			$G_{17}$ & $s_2$ & $4$ & $(1,237,4)$ &
			$\mathbf s_2^4\cdot \mathbf s_1\mathbf s_2\mathbf s_1\mathbf s_2^4(\mathbf s_1\mathbf s_2^3)^3\mathbf s_1$ \\
			$G_{18}$ & $s_2$ & $6$ & $(1,337,4)$ &
			$\mathbf s_2^4\cdot \mathbf s_1^2\mathbf s_2\mathbf s_1^2\mathbf s_2^3\mathbf s_1^2\mathbf s_2^2\mathbf s_1^2\mathbf s_2^3\mathbf s_1$ \\
			$G_{19}$ & $s_3$ & $12$ & $(3,704,4)$ &
			$\mathbf s_3^4\cdot \mathbf z(\mathbf s_1\mathbf s_3\mathbf s_2)^2\mathbf s_1\mathbf s_2(\mathbf s_1\mathbf s_3\mathbf s_2)^2$ \\
			$G_{27}$ & $s_1$ & $1$ & $(1,611,1)$ &
			$\mathbf s_1\cdot (\mathbf s_3\mathbf s_2\mathbf s_1)^4\mathbf s_3\mathbf s_2$
		\end{tabular}
	\end{adjustbox}
	\caption{The triples $(J,i,\eta) \neq (0,1,0)$ used in the construction of  $E_j$.}
	\label{tab:traces}
\end{table}

For each $G_j$, with $4\le j\le 29$ and $j\neq 23,28$, the programs on the project's webpage \cite{web} verify the following statements in $\bar{\cH}(W)$:
	\begin{enumerate}
		\item[(V1)] $\mathrm{pr}(E_j)(\mathrm{pr}(a)\mathbf s)=\mathrm{pr}(E_j)(\mathbf s\,\mathrm{pr}(a))$ for all $a\in\cA_j$ and $\mathbf s\in\mathbf S$; \smallbreak
		\item[(V2)] the blocks $F_r$ of $\gram_{\mathrm{pr}(E_j)}(\mathrm{pr}(\cB_j),\mathrm{pr}(\cB_j)^{-1})$ (see Lemma~\ref{gramdec})
		 satisfy $F_r=0$ for $0<r<d$ and $\det(F_0)=\pm1$;\smallbreak
		\item[(V3)] $\mathrm{pr}(E_j)(1)=1$ and $\mathrm{pr}(E_j)(\mathrm{pr}(b))=0$ for all $b\in\cB_j\setminus\{1\}$;\smallbreak
		\item[(V4)] $\mathrm{pr}(E_j)(\mathrm{pr}(b^{-1}\mathbf z^{|Z(W)|}))=0$ for all $b\in\cB_j\setminus\{1\}$.
	\end{enumerate}

\begin{proposition}\label{barconj}
Let $W=G_j$ with $4\le j\le 29$ and $j\neq 23,28$, and let $\cA_j$,
$\cB_j$ and $E_j$ be as above. Then $\mathrm{pr}(\cB_j)$ is an
$\bar R(W)$-basis of $\bar{\mathcal H}(W)$, and
$\mathrm{pr}(E_j)=\lambda_1^{\mathrm{pr}(B_j)}$ is a symmetrising trace
on $\bar{\mathcal H}(W)$.
\end{proposition}

\begin{proof}
		By (V1) and Proposition~\ref{symcond}, applied to the basis $\mathrm{pr}(\cA_j)$, $\mathrm{pr}(E_j)$ is a trace function. By (V2) and Lemma~\ref{gramdec},
	$\det\gram_{\mathrm{pr}(E_j)}(\mathrm{pr}(\cB_j),\mathrm{pr}(\cB_j)^{-1})=(\det F_0)^d=\pm1$, which is a unit of $\bar R(W)$. Lemma~\ref{pair}
	now shows that $\mathrm{pr}(\cB_j)$ is an $\bar R(W)$-basis of $\bar{\cH}(W)$ and that $\mathrm{pr}(E_j)$ is a
	symmetrising trace. Since $\mathrm{pr}(\cB_j)$ is a basis, (V3) says that $\mathrm{pr}(E_j)=\lambda_1^{\mathrm{pr}(\cB_j)}$.
	\end{proof}
	
	Together with (V4), Proposition~\ref{barconj} shows that the sets
$\mathcal B_j$ satisfy all the conditions of
Definition~\ref{def:pradm}. By Theorem~\ref{thm:liftall}, this yields
the conjecture for these groups.
	
\begin{theorem}\label{thm:Gj}
For every non-real group $W=G_j$ with $4\le j\le 29$, the set
$\mathcal B_j$ is a section basis of $\mathcal H(W)$, and
$\lambda_1^{\mathcal B_j}$ is the canonical symmetrising trace. In
particular, the BMM symmetrising trace conjecture holds for these
groups.
\end{theorem}

\begin{proof}
We show that $\mathcal B_j$ is $\mathrm{pr}$-admissible. By
Definition~\ref{def:zkcandidate}, $\mathcal B_j$ is the set of lifts of
a section of $B(W)\twoheadrightarrow W$. By
Proposition~\ref{barconj}, $\mathrm{pr}(\mathcal B_j)$ is a basis of
$\bar{\mathcal H}(W)$ and $\lambda_1^{\mathrm{pr}(\mathcal B_j)}$ is a
symmetrising trace. Finally, (V4) is condition~(4) of
Definition~\ref{def:pradm}, since $\boldsymbol\pi=\mathbf z^{|Z(W)|}$.
Theorem~\ref{thm:liftall} now gives the conclusion.
\end{proof}

\begin{remark}\label{rem:Ej}
The canonical trace $\tau_j=\lambda_1^{\mathcal B_j}$ of
Theorem~\ref{thm:Gj} can be computed directly from the completed coset
table of $\mathcal A_j$ over $\bar R(W)$. Indeed, $E_j$ is homogeneous
of degree $\epsilon_j=-w(a)-kJ\,w(\mathbf z)$ by Lemma~\ref{exp}, and
$\mathrm{pr}(E_j(1))=1$ by (V3), so
$E_j(1)=\prod_ia_{i,0}^{n_i}$ for some integers $n_i$ by
Lemma~\ref{liftlem}(2). The linear forms $E_j$ and $E_j(1)\,\tau_j$ are
then homogeneous of the same degree and agree after specialisation
(Proposition~\ref{barconj}), so $\tau_j=E_j(1)^{-1}E_j$ by
Lemma~\ref{liftlem}(1).

In practice, for $h$ given as a word in the generators, $\mathrm{pr}(E_j(h))$ is read off from the completed coset table, as the coefficient of $a$ in the product obtained by multiplying $1$ successively by the letters of $\mathbf z^{-kJ}$ and of $h$; lifting it to the homogeneous element of degree $w(h)+\epsilon_j$ gives $E_j(h)$. Hence, $\tau_j(h)$ can be computed without ever expressing $h$ in the basis $\mathcal B_j$.
\end{remark}

\section{Parabolic towers and the BMM conjecture for $G_{31},\ldots,G_{34}$}\label{rem}

In Section~\ref{bases}, we proved the BMM symmetrising trace conjecture for the
non-real groups $G_4,\dots,G_{29}$. Since $G_{30}$, $G_{35}$, $G_{36}$
and $G_{37}$ are real reflection groups, it remains to treat $G_{31}$,
$G_{32}$, $G_{33}$ and $G_{34}$, which we do in this section. As in
Section~\ref{bases}, all computations are carried out over $\bar R(W)$: we
construct a $\mathrm{pr}$-admissible collection for each of these
groups, and the conjecture follows from Theorem~\ref{thm:liftall}.

For $G_{31},G_{32},G_{33}$ and $G_{34}$, the approach of the previous section becomes computationally impractical, due to their large orders. We therefore adopt a different strategy, exploiting, for each such group $W$, a tower of parabolic subgroups $P_0<P_1<\dots<P_r=W$, where $P_0$ is cyclic and $P=P_{r-1}$ is a large parabolic subgroup of $W$. Going up the tower, we construct at each level a basis $\cB$ such that $\lambda_1^{\cB}$ is a symmetrising trace, so that at the next level $\bar{\cH}(P)$ comes with a symmetrising trace $\tau^P$. The remaining condition of the conjecture is needed only for $W$ itself.

The key idea is to decompose $\bar{\cH}(W)$ using a coset tree $X$ with respect to $(S,P)$, so that $\bar{\cH}(W)=\bigoplus_{x\in X}\bar{\cH}(P)x$, where $|X|=|W|/|P|$ is much smaller than $|W|$. This decomposition allows us to construct a linear form on $\bar{\cH}(W)$ from the known trace $\tau^P$ on $\bar{\cH}(P)$. If this form is a trace function, its Gram matrix factors into a block-diagonal matrix determined by the Gram matrix of $\tau^P$ and a smaller matrix over $\bar{\cH}(P)$ of size $|X|\times|X|$. Hence, once the form is known to be a trace function, showing that it is a symmetrising trace reduces to verifying the invertibility of this smaller matrix.

\subsection{A parabolic criterion for symmetrising traces}\label{parcr}
Let $W$ be an exceptional complex reflection group and $P$ a parabolic
subgroup of $W$, and suppose that $\bar{\cH}(P)$ has a section basis
$\cB_P$ such that $\tau^P=\lambda_1^{\cB_P}$ is a
symmetrising trace.
Throughout this subsection, we regard $\bar{\cH}(P)$ as a subalgebra of $\bar{\cH}(W)$ and we extend scalars from $\bar R(P)$ to $\bar R(W)$, so that $\tau^P$ becomes an $\bar R(W)$-linear form $\bar{\cH}(P)\to \bar R(W)$. We denote it again by $\tau^P$. As in the previous sections, $S$ denotes the set of standard generators of $W$, and we assume that $P$ is generated by a subset $S_P\subseteq S$.

Let $X$ be the set of elements of $\bar{\cH}(W)$ corresponding to a coset tree with respect to $(S,P)$, as in Section~\ref{standard}, and suppose that $\cB=\cB_P\cdot X$ is a basis of $\bar{\cH}(W)$ (in particular, $\bar{\cH}(P)$ embeds in
$\bar{\cH}(W)$, since the elements of $\cB_P$ are linearly
independent there). Since the elements of $\cB_P$ lift the elements of $P$ and $X$ is a set of right coset representatives of $P$ in $W$, $\cB$ is then a section basis. In particular, $\bar{\cH}(W)=\bigoplus_{x\in X}\bar{\cH}(P)x$ as a left $\bar{\cH}(P)$-module, so every $h\in\bar{\cH}(W)$ can be written uniquely as $h=\sum_{x\in X}h_x\,x$ with $h_x\in\bar{\cH}(P)$. We define
$$E_X\colon\bar{\cH}(W)\to\bar{\cH}(P),\quad h\mapsto h_1.$$
Since $1\in\cB_P$ and $1\in X$, we have $\lambda_1^{\cB}=\tau^P\circ E_X$. Our goal is to show that $\lambda_1^{\cB}$ is a symmetrising trace of $\bar{\cH}(W)$.

By construction, $E_X$ is left $\bar{\cH}(P)$-linear. The following lemma shows that, if $\lambda_1^{\cB}$ is a trace function, then $E_X$ is also right $\bar{\cH}(P)$-linear. This property allows us to separate the contributions of $\bar{\cH}(P)$ and the coset representatives to the Gram matrix of $\lambda_1^{\cB}$.

\begin{lemma}\label{bimod}
	If $\lambda_1^{\cB}$ is a trace function then $E_X$ is a homomorphism of $\bar{\cH}(P)$-bimodules.
\end{lemma}

\begin{proof}
	Since $E_X$ is left $\bar{\cH}(P)$-linear, it remains to prove right $\bar{\cH}(P)$-linearity.
	Let $h\in\bar{\cH}(W)$ and $p,p'\in\bar{\cH}(P)$. Since $\lambda_1^{\cB}$ is a trace function, we have $\lambda_1^{\cB}(p'hp)=\lambda_1^{\cB}(pp'h)$.
	Since $\lambda_1^{\cB}=\tau^P\circ E_X$ and $E_X$ is left $\bar{\cH}(P)$-linear, we have
	$$
	\lambda_1^{\cB}(p'hp)
	=\tau^P(E_X(p'hp))
	=\tau^P(p'E_X(hp)).
	$$
	Similarly,
	$$
	\lambda_1^{\cB}(pp'h)
	=\tau^P(E_X(pp'h))
	=\tau^P(pp'E_X(h))
	=\tau^P(p'E_X(h)p),
	$$
	where the last equality follows from the fact that $\tau^P$ is a trace function.
	Combining these identities, we obtain
	$
	\tau^P\bigl(p'(E_X(hp)-E_X(h)p)\bigr)=0
	$
	for every $p'\in\bar{\cH}(P)$.
	Since $\tau^P$ is non-degenerate, it follows that
	$E_X(hp)=E_X(h)p.$
	Thus $E_X$ is also right $\bar{\cH}(P)$-linear.
\end{proof}

We now give a criterion for $\lambda_1^{\cB}$ to be a symmetrising trace of $\bar{\cH}(W)$.
Let $Y_1=\{u_1,\ldots,u_N\}$ and $Y_2=\{v_1,\ldots,v_N\}$ be two sets of $N:=|W|/|P|$ elements of $\bar{\cH}(W)$, and consider the collections $\cB_1=\cB_P\cdot Y_1$ and $\cB_2=Y_2\cdot\cB_P$. We define
\[
M_X(Y_1,Y_2)=\bigl(E_X(u_iv_j)\bigr)_{i,j=1}^N\in\operatorname{Mat}_N(\bar{\cH}(P)).
\]
For $m\in\bar{\cH}(P)$, let $\rho(m)\in\operatorname{Mat}_{|P|}(\bar R(W))$ denote the matrix of right multiplication by $m$ on $\bar{\cH}(P)$ with respect to the basis $\cB_P$, that is,
\[
b\,m=\sum_{b'\in\cB_P}\rho(m)_{b,b'}\,b' \qquad (b\in\cB_P).
\]
We extend $\rho$ entrywise to block matrices, obtaining a map
\[
\rho\colon\operatorname{Mat}_N(\bar{\cH}(P))\longrightarrow\operatorname{Mat}_{|W|}(\bar R(W)),\qquad (a_{ij})\longmapsto\bigl(\rho(a_{ij})\bigr)_{i,j}.
\]
Note that $N\cdot|P|=|W|$, so that $\bigl(\rho(a_{ij})\bigr)_{i,j}$ is a square matrix of size $|W|$. Since $\rho$ preserves multiplication and the identity on $\bar{\cH}(P)$, its extension to matrices is an algebra homomorphism, sending the identity matrix to the identity matrix.

The purpose of introducing $\rho$ is to express the Gram matrix of $\lambda_1^{\cB}$ in terms of the matrix $M_X(Y_1,Y_2)$ and the Gram matrix of the parabolic trace $\tau^P$. More precisely, the following lemma shows that, if $E_X$ is an $\bar{\cH}(P)$-bimodule homomorphism, then the Gram matrix admits a factorisation involving $\rho(M_X(Y_1,Y_2))$ and a block-diagonal matrix whose blocks are copies of the Gram matrix of $\tau^P$.

\begin{lemma}\label{twostep}
	Let $G_P=\gram_{\tau^P}(\cB_P,\cB_P)$, and assume that $E_X$ is a homomorphism of $\bar{\cH}(P)$-bimodules. Then
	\[
	\gram_{\lambda_1^{\cB}}(\cB_1,\cB_2)=\rho\bigl(M_X(Y_1,Y_2)\bigr)\cdot\operatorname{diag}(G_P,\ldots,G_P),
	\]
	where the diagonal matrix has $N$ copies of $G_P$. Consequently,
	\[
	\det\gram_{\lambda_1^{\cB}}(\cB_1,\cB_2)=\det\rho\bigl(M_X(Y_1,Y_2)\bigr)\cdot(\det G_P)^N.
	\]
\end{lemma}

\begin{proof}
	Let $G=\gram_{\lambda_1^{\cB}}(\cB_1,\cB_2)$. We index the rows and columns by pairs $(i,b)$, where $1\leq i\leq N$ and $b\in\cB_P$, corresponding to the elements $bu_i\in\cB_1$ (rows) and $v_ib\in\cB_2$ (columns), respectively. Hence, for $b,b'\in\cB_P$, the entry in row $(i,b)$ and column $(j,b')$ is
$$		G_{(i,b),(j,b')}=\lambda_1^{\cB}(bu_i\cdot v_jb')
		=\tau^P\bigl(E_X(bu_iv_jb')\bigr)
		=\tau^P\bigl(bE_X(u_iv_j)b'\bigr),
	$$
	where the last equality follows from the $\bar{\cH}(P)$-bimodule linearity of $E_X$.
	
	By the definition of $\rho$, we have
	\[
	bE_X(u_iv_j)=\sum_{b''\in\cB_P}\rho\bigl(M_X(Y_1,Y_2)\bigr)_{(i,b),(j,b'')}\,b''.
	\]
	Substituting this into the preceding expression, we obtain
	\[
	\begin{aligned}
		G_{(i,b),(j,b')}&=\sum_{b''\in\cB_P}\rho\bigl(M_X(Y_1,Y_2)\bigr)_{(i,b),(j,b'')}\,\tau^P(b''b')\\
		&=\sum_{b''\in\cB_P}\rho\bigl(M_X(Y_1,Y_2)\bigr)_{(i,b),(j,b'')}\,(G_P)_{b'',b'}.
	\end{aligned}
	\]
	This is precisely the matrix identity $G=\rho(M_X(Y_1,Y_2))\cdot\operatorname{diag}(G_P,\ldots,G_P)$. Taking determinants gives $\det G=\det\rho(M_X(Y_1,Y_2))\cdot(\det G_P)^N$, as required.
\end{proof}

Note that $G_P$ is invertible over $\bar R(W)$, since $\tau^P$ is a symmetrising trace. Hence the factorisation above reduces the invertibility of $\gram_{\lambda_1^{\cB}}(\cB_1,\cB_2)$ to that of $\rho\bigl(M_X(Y_1,Y_2)\bigr)$.

The results of this section rely on suitable choices of the sets $Y_1$ and $Y_2$ satisfying the conditions of the following criterion.

\begin{proposition}\label{symcond2} 
	Assume that
	\begin{enumerate}
		\item $\lambda^{\cB}_1$ is a trace function, and
		\item there exist sets $Y_1$ and $Y_2$ of $N$ elements of
		$\bar{\cH}(W)$ such that $\det\rho\big(M_X(Y_1,Y_2)\big)=\pm1$.
	\end{enumerate}
	Then $\lambda^{\cB}_1$ is a symmetrising trace, and
	$\cB_P\cdot Y_1$ and $Y_2\cdot\cB_P$ are bases of
	$\bar{\cH}(W)$.
\end{proposition}

\begin{proof}
	By (1) and Lemma~\ref{bimod}, $E_X$ is a homomorphism of
	$\bar{\cH}(P)$-bimodules, so Lemma~\ref{twostep} gives
	$\det\gram_{\lambda_1^{\cB}}(\cB_P\cdot Y_1,
	Y_2\cdot\cB_P)=\det\rho(M_X(Y_1,Y_2))\cdot(\det G_P)^N$. The
	first factor is $\pm1$ by (2), and $\det G_P$ is a unit since $\tau^P$
	is a symmetrising trace. Lemma~\ref{pair} now gives the
	conclusion.
\end{proof}

We apply Proposition~\ref{symcond2} successively along a tower of
parabolic subgroups, starting from a cyclic one. At each level, the
conclusion provides the symmetrising trace $\tau^P$ needed at the next
level. We describe this procedure in the following subsections.

\subsection{The towers of parabolic subgroups}\label{towers}
For each $W\in\{G_{31},G_{32},G_{33},G_{34}\}$, we fix a tower of standard parabolic subgroups
\[
P_0<P_1<\dots<P_r=W,
\]
where $P_0=\langle s\rangle$ is cyclic with $s\in S$ and, for each $i\ge1$, $P_i$ is generated by a subset $S_i\subseteq S$, with $S_r=S$. We define
\[
\cB_{P_0}=\{1,s,\dots,s^{e-1}\},\qquad
\cB_{P_i}=\cB_{P_{i-1}}\cdot X_i\quad(i\ge1),
\] 
where $e$ is the order of $s$ and $X_i$ is the set of elements of $\bar{\cH}(P_i)$ corresponding to the coset tree with respect to $(S_i,P_{i-1})$ constructed from the ordering $o_i$ of $S_i$. Thus, $\cB_{P_i}$ is a $\BFS_{o_i}$ basis candidate for $\bar{\cH}(P_i)$. The only exception is the last level of $G_{32}$, where $X_r$ corresponds to a seeded tree (Section~\ref{g32}).
At each level, the programs verify that $\cB_{P_i}$ is indeed a basis by completing its coset table over  $\bar R(W)$, as described in Section~\ref{coset}. The orderings $o_i$ determine the corresponding trees, which, together with the coset tables, are available on the project's webpage~\cite{web}.

To apply Proposition~\ref{symcond2} at level $i$, we must verify its
two conditions: the existence of sets $Y_1^i$, $Y_2^i$ such that
$\det\rho\big(M_{X_i}(Y_1^i,Y_2^i)\big)=\pm1$, which we address in
Section~\ref{subsec:Y}, and that $\lambda_1^{\cB_{P_i}}$ is a trace function, which
we address in Section~\ref{subsec:symmetry}. Condition~(3) of Conjecture~\ref{bmmconj} is needed only at the top of the tower, that is, for $W$ itself, and is verified in Section~\ref{subsec:cond3}.

When $P_i$ is a finite Coxeter group, as for the level $A_2$ of $G_{31}$ and the levels $A_2,A_3,A_4$ of $G_{33}$ and $G_{34}$ (see Table~\ref{tab:towers}), breadth-first search reaches
each coset at the level given by the length of its distinguished
representative, along a reduced word for it. Hence, inductively,
$\cB_{P_i}$ is the standard basis of the Iwahori--Hecke algebra
$\bar{\cH}(P_i)$ (compare Proposition~\ref{real}), and
$\lambda_1^{\cB_{P_i}}$ is a symmetrising trace. For these
levels no separate verification is needed; they are also covered by Proposition~\ref{prop:symmetry}.

\subsection{The sets $Y^i_1,Y^i_2$}\label{subsec:Y} In this subsection, we describe how the sets $Y_1^i$ and $Y_2^i$ are chosen.
Since the matrix $M_{X_i}(Y_1^i,Y_2^i)$ is large, our main objective is
to choose these sets so that the resulting matrix is triangular, or at
least block triangular, with invertible diagonal blocks. This structure
allows us to establish the invertibility of the matrix by examining
only its diagonal blocks.

We start by giving some notation.
Let $u=a_1\cdots a_l$ be a word in $S$. We denote by $\wrev u$ the reverse word $a_l\cdots a_1$. If the generators have order $2$, then $\wrev u$ represents the inverse of the element of $W$ represented by $u$. If the generators have order $3$, we write $u=t_1^{c_1}\cdots t_m^{c_m}$ in maximal runs, with $t_j\in S$ and $c_j\in\{1,2\}$, and, since $s^{-c}=s^{3-c}$ in $W$, we denote by
\[
\winv u=t_m^{3-c_m}\cdots t_1^{3-c_1}
\]
the positive word representing the inverse of the element of $W$ represented by $u$.

Every element $x$ of $X_i$ is the lift of a unique word $u_x$ of the coset tree. We write $\wrev x$ and $\winv x$ for the lifts of the words $\wrev{u_x}$ and $\winv{u_x}$, respectively, and, for a subset $X$ of $X_i$, we set $\wrev X=\{\wrev x\mid x\in X\}$ and $\winv X=\{\winv x\mid x\in X\}$.

We now describe the sets $Y^i_1,Y^i_2$, by distinguishing the intermediate levels from the last level of each tower.

\medskip\noindent\textbf{The intermediate levels.} Let $1\le i<r$. By Proposition~\ref{symcond2}, if $\lambda_1^{\cB_{P_i}}$ is a
trace function, the set $\cB_{P_{i-1}}\cdot Y_1^i$ is a basis of
$\bar{\cH}(P_i)$ whenever
$\det\rho\big(M_{X_i}(Y_1^i,Y_2^i)\big)=\pm1$ for some $Y_2^i$. It is therefore natural to take $Y^i_1=X_i$.

For $Y^i_2$ we are inspired by the real case, which we first recall. Suppose that $P_i$ is a finite Coxeter group, as is the case for certain towers considered here. Then $\bar{\cH}(P_i)$ has the standard basis $\{T_w\}_{w\in P_i}$ (Section~\ref{sec-2.2}), with canonical trace $\tau$ satisfying
\[
\tau(T_wT_{w'^{-1}})=\delta_{w,w'}\,q_w, \qquad q_w\in \{\pm1\}.
\]

Assume moreover that $X_i$ comes from a coset tree whose words are reduced words of minimal length in their cosets $P_{i-1}x$, as in Section~\ref{standard}. Then $\cB_{P_i}=\cB_{P_{i-1}}\cdot X_i$ is the standard basis and $\lambda_1^{\cB_{P_i}}=\tau$. By the relation above, the dual basis of the standard basis with respect to the form $(h,h')\mapsto\tau(hh')$ consists of the elements $q_w^{-1}T_{w^{-1}}$, and these are obtained from the standard basis by reversing reduced words, since the reverse of a reduced word is a reduced word for the inverse element. Using Lemma~\ref{bimod}, one deduces that $M_{X_i}(X_i,\wrev{X_i})$ is diagonal with entries $\pm1$, so that $\det\rho\bigl(M_{X_i}(X_i,\wrev{X_i})\bigr)=\pm1$.

In our towers the groups $P_i$ are not all real and the matrix $M_{X_i}(Y^i_1,Y^i_2)$ is not diagonal in general. In accordance with the real case, we build $Y^i_2$ from the inverses of the elements of $X_i$, written as positive words: if the generators have order $2$ we take $Y^i_2=\wrev{X_i}$, and if they have order $3$ we take $Y^i_2=\winv{X_i}$. The levels $i<r$, with their orderings and the corresponding choices of $Y^i_2$, are listed in Table~\ref{tab:towers}.

\begin{table}[h]
	\centering
	\renewcommand{\arraystretch}{1.6}
	\begin{tabular}{llllll}
		$W$ & $i$ & $P_{i-1}<P_i$ & $|P_i:P_{i-1}|$ & $o_i$ & $Y^i_2$ \\
		\hline
		$G_{31}$ & $1$ & $\langle s_5\rangle<A_2=\langle s_3,s_5\rangle$ & $3$ & $(s_3,s_5)$ & $\wrev{X_i}$ \\
		& $2$ & $A_2<G(4,2,3)=\langle s_1,s_2,s_3,s_5\rangle$ & $32$ & $(s_3,s_2,s_1,s_5)$ & $\wrev{X_i}$ \\
		\hline
		$G_{32}$ & $1$ & $\langle s_1\rangle<G_4=\langle s_1,s_2\rangle$ & $8$ & $(s_1,s_2)$ & $\winv{X_i}$ \\
		& $2$ & $G_4<G_{25}=\langle s_1,s_2,s_3\rangle$ & $27$ & $(s_1,s_2,s_3)$ & $\winv{X_i}$ \\
		\hline
		$G_{33}$ & $1$ & $\langle s_1\rangle<A_2=\langle s_1,s_2\rangle$ & $3$ & $(s_1,s_2)$ & $\wrev{X_i}$ \\
		& $2$ & $A_2<A_3=\langle s_1,s_2,s_4\rangle$ & $4$ & $(s_1,s_2,s_4)$ & $\wrev{X_i}$ \\
		& $3$ & $A_3<A_4=\langle s_1,s_2,s_4,s_5\rangle$ & $5$ & $(s_1,s_2,s_4,s_5)$ & $\wrev{X_i}$ \\
		\hline
		$G_{34}$ & $1$ & $\langle s_1\rangle<A_2=\langle s_1,s_2\rangle$ & $3$ & $(s_1,s_2)$ & $\wrev{X_i}$ \\
		& $2$ & $A_2<A_3=\langle s_1,s_2,s_4\rangle$ & $4$ & $(s_1,s_2,s_4)$ & $\wrev{X_i}$ \\
		& $3$ & $A_3<A_4=\langle s_1,s_2,s_4,s_5\rangle$ & $5$ & $(s_1,s_2,s_4,s_5)$ & $\wrev{X_i}$ \\
		& $4$ & $A_4<G_{33}=\langle s_1,\dots,s_5\rangle$ & $432$ & $(s_1,\dots,s_5)$ & $\wrev{X_i}$ \\
		& & & & &\\
	\end{tabular}
	\caption{The intermediate levels $i<r$ of the towers}
	\label{tab:towers}
\end{table}

\noindent\textbf{The final level.} For the final level $P_{r-1}<P_r=W$, the rule above has to be modified for some of the groups; we describe $X_r$ and $Y^r_1,Y^r_2$ for each group in the following subsections. The programs verifying the properties of $M_{X_r}(Y^r_1,Y^r_2)$ are available on the project's webpage~\cite{web}.

\subsubsection{The groups $G_{33}$ and $G_{34}$}\label{g3334}
For $G_{33}$ and $G_{34}$ the last level is $A_4<G_{33}$ and $G_{33}<G_{34}$, respectively, with $|G_{33}:A_4|=432$ and $|G_{34}:G_{33}|=756$. As for the intermediate levels, $X_r$ corresponds to the coset tree with respect to $(S,P_{r-1})$ constructed from the ordering $o_r$ of $S$, where $o_r=(s_1,\dots,s_5)$ for $G_{33}$ and $o_r=(s_1,\dots,s_6)$ for $G_{34}$, respectively; these are the coset words of~\cite{MaPf}. Since the generators have order $2$, we take once again  $Y^r_1=X_r$ and $Y^r_2=\wrev{X_r}$.

\subsubsection{The group $G_{31}$}\label{g31}
In this case, $P_{r-1}=G(4,2,3)=\langle s_1,s_2,s_3,s_5\rangle$ and $|G_{31}:P_{r-1}|=240$. The set $X_r$ corresponds to the coset tree with respect to $(S,P_{r-1})$, constructed from the ordering $o_r=(s_2,s_3,s_4,s_5,s_1)$ of $S$. The resulting basis $\cB_{P_{r-1}}\cdot X_r$ is one of the $\BFS_o$ candidates of Figure~\ref{bfs} that pass Modular test~2. 

However, with $Y^r_1=X_r$ and $Y^r_2=\wrev{X_r}$ the matrix $M_{X_r}(Y^r_1,Y^r_2)$, of size $240\times240$, has $45$ diagonal blocks of size $1$ and one of size $195$, so that its block structure is of little use for the computation of its determinant.
 We therefore choose $Y^r_1$ to be the set of elements corresponding to a second coset tree with respect to $(S\cup\{z\},P_{r-1})$, constructed from the ordering $o'_r=(s_1,\dots,s_5,z)$, where the central element $z$ is considered as an additional letter. This tree consists of $240$ words, $140$ of which contain $z$. Each word is of the form $w_az^{k_a}$, with $w_a$ a word in $S$ and $k_a\in\{0,1,2,3\}$ (recall that $z$ has order $4$). We denote by $\mathbf w_a$ the lift of the word $w_a$, and by $\wrev{\mathbf w_a}$ the lift of the word $\wrev{w_a}$, as in Section~\ref{standard}. Writing
\[
Y^r_1=\{u_1,\dots,u_{240}\},\qquad u_a=\mathbf w_a\,\mathbf z^{k_a},
\]
we define
\[
Y^r_2=\{v_1,\dots,v_{240}\},\qquad v_a=\wrev{\mathbf w_a}\,\mathbf z^{-k_a}.
\]

The introduction of $z$ is crucial for obtaining a block triangular matrix. Indeed, since $\mathbf z$ is central, the $(a,b)$ entry of $M_{X_r}(Y^r_1,Y^r_2)$ is
\[
E_{X_r}(u_av_b)=E_{X_r}\bigl(\mathbf w_a\,\wrev{\mathbf w_b}\,\mathbf z^{\,k_a-k_b}\bigr)\in\bar{\cH}(P_{r-1}).
\]
Direct computations carried out by the programs show that this entry vanishes whenever $k_a>k_b$. Thus, by ordering the rows and columns according to the exponents $k_a$ and $k_b$, respectively, we obtain a block triangular matrix.

Moreover, the programs verify that the diagonal blocks, which correspond to equal exponents $k_a=k_b$, are themselves triangular, so that all $240$ diagonal blocks have size $1$, and that each diagonal entry is invertible in $\bar{\cH}(P_{r-1})$.  Since the units of $\bar R(W)$ are $\pm1$, we obtain $\det\rho\bigl(M_{X_r}(Y^r_1,Y^r_2)\bigr)=\pm1$.

\subsubsection{The group $G_{32}$}\label{g32}
In this case, $P_{r-1}=G_{25}=\langle s_1,s_2,s_3\rangle$ and $|G_{32}:P_{r-1}|=240$. Recall that the central element $z$ of $G_{32}$ is of order $6$, and that the normaliser of $P_{r-1}$ is $P_{r-1}\times\langle z\rangle$. Hence, the six cosets $P_{r-1}z^j$, $0\le j\le5$, are precisely those cosets $P_{r-1}g$ satisfying $P_{r-1}gP_{r-1}=P_{r-1}g$. We consider now the coset $P_{r-1}z^5$. Let $X'_r$ be the set of elements corresponding to a coset tree with respect to $(S,P_{r-1})$ constructed from an ordering of $S$, and let $\cB=\cB_{P_{r-1}}\cdot X'_r$. Modular computations show that, for every such ordering, the element $x\in X'_r$ representing this coset satisfies
\[
\theta\bigl(E_{X'_r}(x\,\mathbf s_2)\bigr)\neq0,
\]
where $\theta$ is the specialisation modulo $p$ introduced in Section~\ref{modtests}. Consequently, $E_{X'_r}(x\,\mathbf s_2)\neq0$. On the other hand, Lemma~\ref{bimod} implies that, if $\lambda_1^{\cB}$ is a trace function, then
\[
E_{X'_r}(x\,\mathbf s)=0
\qquad
\text{for all }x\in X'_r\setminus\{1\}
\text{ and }\mathbf s\in\{\mathbf s_1,\mathbf s_2,\mathbf s_3\}.
\]
Hence $\lambda_1^{\cB}$ is not a trace function, and $\cB$ does not pass Modular test~2 (Lemma~\ref{mtest2}). This rules out all coset trees obtained directly from an ordering of $S$, in accordance with Figure~\ref{bfs}, where none of the $22$ candidates for $G_{32}$ passes Modular test~2.

To avoid this, we follow the method of~\cite[\S4.2.4]{BCCK} for $G_7$: we replace the word of the coset $P_{r-1}z^5$ by another word representing the same coset. This removes the failure $E_{X'_r}(x\,\mathbf s_2)\neq0$. However, in our computations a corresponding problem then appears for the coset $P_{r-1}z^4$: the Gram matrix of $\lambda_1^{\cB}$ is still not symmetric, that is, the resulting set does not pass Modular test~2, until the word of $P_{r-1}z^4$ is replaced as well.  We therefore prescribe the words representing both cosets $P_{r-1}z^5$ and $P_{r-1}z^4$. Let
\begin{align*}
	\sigma_1&=s_4^2s_3s_2^2s_3s_4^2s_1s_2^2s_1^2s_2s_3^2s_4,\\
	\sigma_2&=s_4s_3^2s_2s_4s_3^2s_2s_1^2s_2s_3^2s_2s_4s_3^2s_4
\end{align*}
be words representing $P_{r-1}z^5$ and $P_{r-1}z^4$, respectively. These words were found by a computer search, as shortest positive words (of lengths $17$ and $18$, respectively) of the cosets $P_{r-1}z^5$ and $P_{r-1}z^4$ with the following two properties: their prefixes lie in pairwise distinct cosets, and the resulting coset tree passes Modular test~2. The choice among such words is not essential.

We now let $X_r$ be the set of elements corresponding to a coset tree with respect to $(S,P_{r-1})$ that contains all prefixes of the words $\sigma_1$ and $\sigma_2$, and that is completed by breadth-first search with respect to $o_r=(s_1,s_2,s_3,s_4)$. The programs verify that the coset table of $\cB_{P_{r-1}}\cdot X_r$ closes (see Section~\ref{coset}), so this is a basis of $\bar{\cH}(G_{32})$. 

We take $Y^r_1=X_r=\{x_1,\dots,x_{240}\}$ and
\[
Y^r_2=\{v_1,\dots,v_{240}\},\qquad v_i=\winv{x_{m_i}}\,\mathbf z^{k_i},
\]
with $m_i\in\{1,\dots,240\}$ and $k_i\in\{0,\dots,5\}$, where $v_i$ represents the coset $g_i^{-1}P_{r-1}$, with $g_i\in W$ the element represented by $x_i$. The pairs $(m_i,k_i)$ are listed on the project's webpage~\cite{web}.
\subsubsection{The invertibility of $M_{X_i}(Y_1^i,Y_2^i)$}
The following proposition collects the verification of Condition~(2) of Proposition~\ref{symcond2} at all the levels of the towers for which the choices of $Y_1^i$ and $Y_2^i$ were described above.
\begin{proposition}\label{prop:Minv}
	For every level $i\ge1$ of the towers of Section~\ref{towers} at which $P_i$ is
	not a Coxeter group, with $Y_1^i$ and $Y_2^i$ as above,
	$\det\rho\big(M_{X_i}(Y_1^i,Y_2^i)\big)=\pm1$.
\end{proposition}

\begin{proof}
	For each such level, the programs verify that there are permutations of the rows and of the columns of $M:=M_{X_i}(Y_1^i,Y_2^i)$ for which the matrix becomes block triangular, and that every diagonal block $D$ is invertible as a matrix over $\bar{\cH}(P_{i-1})$. More precisely, a diagonal block of size $1$ is an entry of $\bar{\cH}(P_{i-1})$ that the programs recognise as a unit (for instance, $\pm b$ with $b\in\cB_{P_{i-1}}$). The diagonal blocks of larger size are brought to triangular form by Gaussian elimination over $\bar{\cH}(P_{i-1})$, in which all pivots are units. Such blocks occur only at the level $A_4<G_{33}$ of the towers of $G_{33}$ and $G_{34}$, where there is a single block of size $6$, and at the last level of $G_{34}$, where there are two blocks, of sizes $18$ and $3$. At all other levels, including all the levels of the towers of $G_{31}$ and $G_{32}$, every diagonal block has size $1$. Since $\rho$ is a ring homomorphism, $\rho(D)$ is then invertible over $\bar R(W)$, and so $\det\rho(D)$ is a unit of $\bar R(W)$, that is, $\pm1$. Permuting the rows or the columns of $M$ permutes blocks of rows or columns of $\rho(M)$, so it changes the determinant only by a sign. Since the determinant of a block triangular matrix is the product of the determinants of its diagonal blocks, we obtain $\det\rho(M)=\pm1$.
\end{proof}

\subsection{Checking the trace property}\label{subsec:symmetry}
In this subsection we explain how Condition~(1) of Proposition~\ref{symcond2}, namely that $\lambda_1^{\cB}$ is a trace function, is verified at each level of the towers of Section~\ref{towers}.

Let $W$, $P$, $S$, $\cB_P$, $X$ and $\cB=\cB_P\cdot X$ be as in Section~\ref{parcr}, and let $\mathbf S$ and $\mathbf S_P$ denote the sets of braid generators of $W$ and $P$, lifting $S$ and $S_P=S\cap P$. We assume that $\tau^P=\lambda_1^{\cB_P}$ is a trace function on $\bar{\cH}(P)$; this is the only property of $\tau^P$ that we use in this subsection. It holds at the bottom of the tower, where $\bar{\cH}(P_0)$ is commutative, and at each subsequent level by the induction hypothesis. Recall that $\lambda_1^{\cB}=\tau^P\circ E_X$, that $E_X(x_m)=\delta_{m,1}$ for the elements $x_m$ of $X$, and that $E_X$ is left $\bar{\cH}(P)$-linear, since $\cB=\cB_P\cdot X$.

To prove that $\lambda_1^{\cB}$ is a trace function, we apply Proposition~\ref{symcond} with $Y_{\mathbf s}=\cB$ for every generator $\mathbf s\in\mathbf S$. Thus, it suffices to verify that
\begin{equation}\label{eq:symid}
	\lambda_1^{\cB}(b\,x\,\mathbf s)=\lambda_1^{\cB}(\mathbf s\,b\,x)\qquad\text{for all }\mathbf s\in\mathbf S,\ b\in\cB_P,\ x\in X.
\end{equation}
These are $|S|\cdot|W|$ identities. To verify them, we distinguish between the generators of $\bar{\cH}(W)$ that lie in $\bar{\cH}(P)$, that is, those in $\mathbf S_P$, and the remaining ones. For the first ones, the identities follow from a simple condition on the coset tree $X$, as the following lemma states.

\begin{lemma}\label{lem:genP}
	Let $\mathbf s\in\mathbf S_P$ and assume that $E_X(x\mathbf s)=\delta_{x,1}\,\mathbf s$ for all $x\in X$. Then \eqref{eq:symid} holds for $\mathbf s$.
\end{lemma}
\begin{proof}
Let $b\in\cB_P$ and $x\in X$. By the left $\bar{\cH}(P)$-linearity of $E_X$,
\[
\lambda_1^{\cB}(b\,x\,\mathbf s)=\tau^P\bigl(b\,E_X(x\mathbf s)\bigr)=\delta_{x,1}\,\tau^P(b\mathbf s).
\]
Since $\mathbf s b\in\bar{\cH}(P)$, we have $E_X(\mathbf s b\,x)=\mathbf s b\,E_X(x)=\delta_{x,1}\,\mathbf s b$, so $\lambda_1^{\cB}(\mathbf s\,b\,x)=\delta_{x,1}\,\tau^P(\mathbf s b)$. The two sides agree, because $\tau^P$ is a trace function.
\end{proof}
The following lemma holds for every generator $\mathbf s\in\mathbf S$, but we use it only for the generators not in $\mathbf S_P$, since the identities for the generators in $\mathbf S_P$ are already covered by Lemma~\ref{lem:genP}.

\begin{lemma}\label{lem:othergens}
	Let $\mathbf s\in\mathbf S$, $b\in\cB_P$ and $x\in X$, and write $\mathbf s\,b=\sum_{n}r_n\,x_n$ with $r_n\in\bar{\cH}(P)$, $x_n\in X$. Then
	\[
	\lambda_1^{\cB}(b\,x\,\mathbf s)=\tau^P\bigl(b\,E_X(x\mathbf s)\bigr),\qquad
	\lambda_1^{\cB}(\mathbf s\,b\,x)=\sum_n\tau^P\bigl(r_n\,E_X(x_nx)\bigr).
	\]
\end{lemma}
\begin{proof}
	Since $E_X$ is left $\bar{\cH}(P)$-linear, we have
	$
	\lambda_1^{\cB}(b\,x\,\mathbf s)=\tau^P\bigl(E_X(b\,x\mathbf s)\bigr)=\tau^P\bigl(b\,E_X(x\mathbf s)\bigr).
	$
	Similarly, using the expansion of $\mathbf s\,b$, we obtain
	\[
	\begin{aligned}
		\lambda_1^{\cB}(\mathbf s\,b\,x)
		&=\tau^P\bigl(E_X(\mathbf s\,b\,x)\bigr)
		=\tau^P\Bigl(E_X\Bigl(\sum_n r_n\,x_nx\Bigr)\Bigr)
		=\sum_n\tau^P\bigl(r_n\,E_X(x_nx)\bigr).
	\end{aligned}
	\]
\end{proof}

\begin{proposition}\label{prop:symmetry}
	Let $W\in\{G_{31},G_{32},G_{33},G_{34}\}$, let $P_0<\dots<P_r=W$ be the tower of Section~\ref{towers}, and let $1\le i\le r$. If $\lambda_1^{\cB_{P_{i-1}}}$ is a trace function, then so is $\lambda_1^{\cB_{P_i}}$. In particular, $\lambda_1^{\cB_{P_i}}$ is a trace function for every $0\le i\le r$.
\end{proposition}
\begin{proof}
	Let $P=P_{i-1}$, $X=X_i$ and $\cB=\cB_{P_i}$, so that $\lambda_1^{\cB}=\tau^P\circ E_X$, where $\tau^P=\lambda_1^{\cB_P}$. Let $\mathbf S_{i-1}\subseteq\mathbf S_i$ denote the sets of generators of $P_{i-1}$ and $P_i$, respectively. At every level of each tower, the programs verify that
	\[
	E_X(x\mathbf s)=\delta_{x,1}\mathbf s
	\qquad\text{for all }x\in X,\ \mathbf s\in\mathbf S_{i-1},
	\]
	by a finite check of the rows of the coset table. By Lemma~\ref{lem:othergens}, for each generator $\mathbf s\in\mathbf S_i$ that does not lie in $\mathbf S_{i-1}$, the identities~\eqref{eq:symid} become
	\[
	\tau^P\bigl(b\,E_X(x\mathbf s)\bigr)
	=\sum_n\tau^P\bigl(r_n\,E_X(x_nx)\bigr),
	\qquad b\in\cB_P,\ x\in X,
	\]
	where $\mathbf s\,b=\sum_n r_nx_n$. Both sides can be computed using only right multiplication, from the coset table of $\bar{\cH}(P_i)$ over $\bar{\cH}(P)$ and the Gram matrix of $\tau^P$ with respect to $\cB_P$. For example, for the last generator $\mathbf s_6$ of $G_{34}$, the column $(E_X(x_m\mathbf s_6))_m$ has exactly one nonzero entry. The programs verify these identities for every such generator $\mathbf s$ and are available on the project's webpage~\cite{web}.
	
	By Lemma~\ref{lem:genP}, the identities~\eqref{eq:symid} also hold for the generators in $\mathbf S_{i-1}$. Thus, they hold for every generator in $\mathbf S_i$. Since $\cB$ is a basis of $\bar{\cH}(P_i)$, Proposition~\ref{symcond}, applied with $Y_{\mathbf s}=\cB$, establishes the first claim.
	
	For the last claim, we proceed by induction on $i$. At the bottom of the tower, $\bar{\cH}(P_0)$ is commutative, so $\lambda_1^{\cB_{P_0}}$ is a trace function. The result then follows at each subsequent level from the first claim.
\end{proof}

\subsection{Checking the specialised form of Condition~(3)}\label{subsec:cond3}
As mentioned at the beginning of this section, Condition~(3) of Conjecture~\ref{bmmconj} is needed only for $W$ itself, and only in its specialised form; the condition itself then follows from Theorem~\ref{thm:liftall}. We may therefore simplify the notation in this subsection as follows. Let $\cB=\cB_P\cdot X$ be the basis of $\bar{\cH}(W)$ at the top level of the tower, so that $P=P_{r-1}$ and $X=X_r$, and let $\lambda_1^{\cB}=\tau^P\circ E_X$, which is a trace function by Proposition~\ref{prop:symmetry}. Once $\lambda^{\mathcal B}_1(\boldsymbol\pi)=1$, the specialised form of Condition~(3) amounts to the identity $\lambda_1^{\cB}(h\boldsymbol\pi)=\lambda_1^{\cB}(\iota(h))^\vee$ for $h$ in a basis of $\bar{\cH}(W)$ (see Lemma~\ref{lem:zeta} below), that is, to $|W|$ identities. Instead of verifying them one by one, we reduce their number by using a central element $\boldsymbol\zeta\in B(W)$ satisfying $\boldsymbol\zeta^2=\boldsymbol\pi$, namely a power of $\mathbf z$ (Table~\ref{tab:zeta}). Since $\boldsymbol\pi\boldsymbol\zeta^{-1}=\boldsymbol\zeta$, both sides can be expressed, with respect to a suitable basis, in terms of the linear form
\[
h\longmapsto\lambda_1^{\cB}(h\boldsymbol\zeta),
\]
which, in the cases under consideration, is supported on a single coset.

We first need an anti-automorphism of $\bar{\cH}(W)$ that implements the inversion in $B(W)$: it is the anti-automorphism of \cite[Proposition~1.26]{BMM},
restricted to $\mathcal H(W)$ and specialised by $\mathrm{pr}$.

\begin{lemma}\label{lem:iota}
	There is a bijection $\iota\colon\bar{\cH}(W)\to\bar{\cH}(W)$ such that $\iota(rh)=r^\vee\iota(h)$ for $r\in \bar R(W)$, $\iota(hh')=\iota(h')\iota(h)$, and $\iota(\mathbf b)=\mathbf b^{-1}$ for all $\mathbf b\in B(W)$.
\end{lemma}
\begin{proof}
	Let $\mathbf s\in\mathbf S$ have order $e$ and Hecke relation $\mathbf s^e=\sum_{j<e}a_j\mathbf s^j$. Multiplying it by $a_0^{-1}\mathbf s^{-e}$ gives $\mathbf s^{-e}=\sum_{j<e}a_j^\vee\mathbf s^{-j}$, with $a_0^\vee=a_0^{-1}$ and $a_j^\vee=-a_0^{-1}a_{e-j}$. Hence, the $\vee$-semilinear anti-homomorphism from the free algebra on $\mathbf S$ to $\bar{\cH}(W)$ sending $\mathbf s\mapsto\mathbf s^{-1}$ kills the Hecke relations. It also kills the braid relations $u=v$, since they are relations between positive words, which are sent to $u^{-1}=v^{-1}$. It therefore induces $\iota\colon\bar{\cH}(W)\to\bar{\cH}(W)$ with the stated properties on generators, and $\iota(\mathbf b)=\mathbf b^{-1}$ for all $\mathbf b\in B(W)$ follows. Finally, $\iota(\mathbf s^{-1})=\iota(\mathbf s)^{-1}=\mathbf s$, so $\iota^2$ is $\vee$-semilinear, $\vee^2=\mathrm{id}$ and $\iota^2$ fixes the generators; hence $\iota^2=\mathrm{id}$ and $\iota$ is bijective.
\end{proof}

Write $h=\sum_{x\in X}h_x\,x$ with $h_x\in\bar{\cH}(P)$. For $x_m\in X$ put $\psi_m=E_X(x_m\boldsymbol\zeta)$, so that $\lambda_1^{\cB}(h\boldsymbol\zeta)=\tau^P\bigl(\sum_mh_m\psi_m\bigr)$.

\begin{lemma}\label{lem:zeta}
	Let $\boldsymbol\zeta\in B(W)$ be central with $\boldsymbol\zeta^2=\boldsymbol\pi$, and assume that $\lambda_1^{\cB}(\boldsymbol\pi)=1$ and that $\psi_m=0$ for all $m\neq c$, for some $c$. Assume moreover that
	\begin{enumerate}
		\item $E_X(x_m^{-1}\boldsymbol\zeta)=0$ for all $m\ne c$, and \smallbreak
		\item $\tau^P\bigl(E_X(x_c^{-1}\boldsymbol\zeta)\,b^{-1}\bigr)=\tau^P(b\,\psi_c)^\vee$ for all $b\in\cB_P$.
	\end{enumerate}
	Then $\lambda^{\mathcal B}_1(\beta^{-1})^\vee=\lambda^{\mathcal B}_1(\beta\boldsymbol\pi)$
	for all $\beta\in B(W)$.
\end{lemma}
\begin{proof}
	Let $f(h)=\lambda_1^{\cB}(h\boldsymbol\pi)$ and $g(h)=\lambda_1^{\cB}(\iota(h))^\vee$. Both are $\bar R(W)$-linear, and, by Lemma~\ref{lem:iota} and $\lambda_1^{\cB}(\boldsymbol\pi)=1$, the conclusion says that $f(\beta)=g(\beta)$ for all $\beta\in B(W)$. Hence it is enough to prove $f=g$ on an $\bar R(W)$-basis of $\bar{\cH}(W)$. We take the basis $\{b\,x_m\boldsymbol\zeta^{-1}\mid b\in\cB_P,\ m\}$, that is, $\cB$ multiplied on the right by the unit $\boldsymbol\zeta^{-1}$. The elements of $\cB_P$ and of $X$ are images of elements of $B(W)$, so Lemma~\ref{lem:iota} applies to them.
	
	As $\boldsymbol\pi\boldsymbol\zeta^{-1}=\boldsymbol\zeta$ and $E_X$ is left $\bar{\cH}(P)$-linear, $f(b\,x_m\boldsymbol\zeta^{-1})=\tau^P(b\,\psi_m)$. On the other hand, $\iota(b\,x_m\boldsymbol\zeta^{-1})=\boldsymbol\zeta x_m^{-1}b^{-1}$, and, since $\boldsymbol\zeta$ is central and $\lambda_1^{\cB}$, $\tau^P$ are trace functions,
	\[
	\lambda_1^{\cB}\bigl(\boldsymbol\zeta x_m^{-1}b^{-1}\bigr)=\lambda_1^{\cB}\bigl(b^{-1}x_m^{-1}\boldsymbol\zeta\bigr)=\tau^P\bigl(b^{-1}E_X(x_m^{-1}\boldsymbol\zeta)\bigr)=\tau^P\bigl(E_X(x_m^{-1}\boldsymbol\zeta)\,b^{-1}\bigr).
	\]
	For $m\neq c$ we have $f=0$ since $\psi_m=0$, and $g=0$ by~(1). For $m=c$, the equality $f=g$ is~(2) after applying $\vee$. Hence $f=g$.
\end{proof}
For the groups of Table~\ref{tab:zeta}, the programs verify the hypotheses of Lemma~\ref{lem:zeta} as follows. They check that $\boldsymbol\zeta$ is central and that $\psi_m=0$ outside a single coset $c$ ($c=240,239,432,756$ for $G_{31},G_{32},G_{33},G_{34}$, respectively), and that $\lambda_1^{\cB}(\boldsymbol\pi)=\tau^P\bigl((\boldsymbol\zeta)_c\,\psi_c\bigr)=1$, where $(\boldsymbol\zeta)_c$ is the $x_c$-coordinate of $\boldsymbol\zeta$. They then compute the coordinates of $x_m^{-1}$ along the coset tree, using $E_X(x_m^{-1}\boldsymbol\zeta)=(x_m^{-1})_c\,\psi_c$, and verify~(1) for the other $|X|-1$ cosets and~(2) for the $|\cB_P|$ elements $b$. For $G_{31}$, for example, this is $239+192=431$ identities instead of $|W|=46\,080$.

\begin{table}[h]
	\centering
	\begin{tabular}{lccc}
		\hline
		$W$ & $z$ & order of $z$ & $\boldsymbol\zeta$\\
		\hline
		$G_{31}$ & $(s_1s_4s_2s_3s_5)^6$ & $4$ & $\mathbf z^2$\\
		$G_{32}$ & $(s_1s_2s_3s_4)^5$ & $6$ & $\mathbf z^3$\\
		$G_{33}$ & $(s_1s_2s_3s_4s_5)^9$ & $2$ & $\mathbf z$\\
		$G_{34}$ & $(s_1\cdots s_6)^7$ & $6$ & $\mathbf z^3$\\
		\hline
		& &  &\\
	\end{tabular}
	\caption{The  element $\boldsymbol\zeta$  used for Condition~(3).}
	\label{tab:zeta}
\end{table}
\subsection{Conclusion}\label{subsec:concl5}
We now summarise the results of this section, and deduce the conjecture for $G_{31},\dots,G_{34}$.

\begin{proposition}\label{prop:top}
	Let $W\in\{G_{31},G_{32},G_{33},G_{34}\}$, and let
	$\mathcal B=\mathcal B_{P_{r-1}}\cdot X_r\subseteq\mathcal H(W)$ be the
	set of braid lifts of Section~\ref{towers}. Then $\mathcal B$ is
	$\mathrm{pr}$-admissible.
\end{proposition}

\begin{proof}
	We verify the four conditions of Definition~\ref{def:pradm}, writing
	$\tau=\lambda_1^{\mathrm{pr}(\mathcal B)}$.
	
	\emph{Step 1: Conditions (1) and (2).} By construction, $\mathcal B$ is
	the set of lifts of a section of $B(W)\twoheadrightarrow W$. By
	Section~\ref{towers}, at each level $i$ the programs verify that
	$\mathrm{pr}(\mathcal B_{P_i})$ is a basis of $\bar{\mathcal H}(P_i)$ by
	completing its coset table over $\bar R(W)$. In particular,
	$\mathrm{pr}(\mathcal B)$ is a basis of $\bar{\mathcal H}(W)$.
	
	\emph{Step 2: Condition (3).} We prove by induction on $i$ that
	$\lambda_1^{\mathcal B_{P_i}}$ is a symmetrising trace of
	$\bar{\mathcal H}(P_i)$. For $i=0$ the result is obvious, since the
	algebra $\bar{\mathcal H}(P_0)$ is commutative and $\mathcal B_{P_0}$ is
	its standard basis. Let $i\ge1$. By Proposition~\ref{prop:symmetry},
	$\lambda_1^{\mathcal B_{P_i}}$ is a trace function. If $P_i$ is a finite
	Coxeter group, then $\lambda_1^{\mathcal B_{P_i}}$ is a symmetrising
	trace by Section~\ref{towers}. Otherwise,
	$\tau^{P_{i-1}}=\lambda_1^{\mathcal B_{P_{i-1}}}$ is a symmetrising
	trace by the induction hypothesis, and Proposition~\ref{prop:Minv} gives
	$\det\rho\big(M_{X_i}(Y_1^i,Y_2^i)\big)=\pm1$. Hence
	Proposition~\ref{symcond2} shows that $\lambda_1^{\mathcal B_{P_i}}$ is a
	symmetrising trace. For $i=r$, this is condition~(3).
	
\emph{Step 3: Condition (4).} By Lemma~\ref{lem:zeta} and the
verifications of \S\ref{subsec:cond3}, $\tau(\beta^{-1})^\vee=\tau(\beta\boldsymbol\pi)$ for all
$\beta\in B(W)$. Applying this to $\beta^{-1}$ gives
$\tau(\beta^{-1}\boldsymbol\pi)=\tau(\beta)^\vee$ for all $\beta\in B(W)$. For $b\in\mathcal B\setminus\{1\}$, let $\beta\in B(W)$ be an element mapping to $b$. Then $\tau(\beta^{-1}\boldsymbol\pi)=\tau(\beta)^\vee=0$, since $\tau=\lambda_1^{\mathrm{pr}(\mathcal B)}$.
\end{proof}

\begin{theorem}\label{thm:top}
For $W\in\{G_{31},G_{32},G_{33},G_{34}\}$, the BMM symmetrising trace
conjecture holds.
\end{theorem}

\begin{remark}\label{rem:parfree}
	The towers of \S \ref{towers} also yield the parabolic freeness conjecture
	(Conjecture~\ref{pfc}) for every pair of consecutive levels
	$P_{j-1}<P_j$, $1\le j\le r$, of each tower. Indeed,
	$\mathcal B_{P_j}=\mathcal B_{P_{j-1}}\cdot X_j$ consists of braid
	lifts and its image is a basis of $\bar{\mathcal H}(P_j)$, so by
	Proposition~\ref{prop:lift}(1), applied to $P_j$, it is a basis of
	$\mathcal H(P_j)$. Hence
	\[
	\mathcal H(P_j)=\bigoplus_{x\in X_j}\mathcal H(P_{j-1})\,x
	\]
	is a free left $\mathcal H(P_{j-1})$-module of rank
	$|X_j|=|P_j:P_{j-1}|$. Applying the anti-automorphism of \cite[Proposition~1.26]{BMM}, which
maps $\mathcal H(P_{j-1})$ onto itself and $x$ to $x^{-1}$, we obtain
$\mathcal H(P_j)=\bigoplus_{x\in X_j}x^{-1}\,\mathcal H(P_{j-1})$, so
$\mathcal H(P_j)$ is also a free right $\mathcal H(P_{j-1})$-module of
the same rank. By
	transitivity, the same holds for every pair $P_i<P_j$ of levels of a
	tower. This includes the pairs $(G_{31},G(4,2,3))$, $(G_{32},G_4)$
	and $(G_{33},A_3)$, which do not appear in the list of \S \ref{sec-2.3}.
\end{remark}

\section{The groups $G(de,e,2)$ with $e$ even}\label{sec-7}

For the groups $G(de,e,n)$ of the infinite series, a symmetrising trace
$\tau$ on $\mathcal H(W)$ satisfying Condition~(2) of
Conjecture~\ref{bmmconj} was constructed in \cite{BreMa, MM98}, and
Condition~(3) for this trace follows from \cite[Lemma~2.7]{BMM} and
\cite[Theorem~5.2]{GIM}, except when $n=2$ and $e$ is even. We will treat here this remaining case.
Note that the symmetrising trace for this case is only constructed in \cite{MM98}.

\subsection{Reduction to $e=2$}
Let $e=2f$ be even. For $d=1$ the group $G(e,e,2)$ is a
dihedral group, for which the conjecture is classical. Let $d \geq 2$ and let
$W'=G(de,e,2)$, a normal subgroup of index $f$ of $W=G(de,2,2)$.
The groups $W'$ and $W$ have the same reflecting
hyperplanes, hence the same hyperplane complement. Thus $B(W')$ is a
subgroup of $B(W)$, and the central element $\boldsymbol\pi$ is the same
for both.

The generic Hecke algebra $\mathcal H(W)$ is generated 
over a Laurent polynomial ring of the form
$$\mathbb{Z}[x_0,x_0^{-1},x_1,x_1^{-1},y_0,y_0^{-1},y_1,y_1^{-1},z_0,z_0^{-1},z_1,z_1^{-1},\ldots,z_{fd-1},z_{fd-1}^{-1}]$$
by three elements $\bf s$, $\bf t$, $\bf u$, 
subject
to ${\bf stu}={\bf tus}={\bf ust}$ and
\[
(\bs-x_0)(\bs-x_1)=(\bt-y_0)(\bt-y_1)=(\bu-z_0)(\bu-z_1)\cdots(\bu-z_{fd-1})=0.
\]
By \cite[Proposition~3.9]{Ma2}, which is deduced from
\cite[Proposition~1.6]{Ar} (see also \cite[\S 5]{MM98} and
\cite[Lemma~A.7.3]{Ch09}), if the parameters of $\bu$ are
specialised to all the $f$-th roots of $u_0,\dots,u_{d-1}$, keeping the
other parameters, then $\mathcal H(W')$ is the subalgebra of  this specialisation of $\mathcal{H}(W)$ generated by $\bs$, $\bt$ and $\bu^f$ over 
$$\mathbb{Z}[x_0,x_0^{-1},x_1,x_1^{-1},y_0,y_0^{-1},y_1,y_1^{-1},u_0,u_0^{-1},u_1,u_1^{-1},\ldots,u_{d-1},u_{d-1}^{-1}].$$
 By \cite[\S 5]{MM98}, the
trace $\tau'$ of  $\mathcal H(W')$ is the restriction of
the specialised trace $\tau$ of $\mathcal H(W)$.

\begin{proposition}\label{prop:reduction}
If $\tau$ satisfies Condition~(3) of Conjecture~\ref{bmmconj} on
$\mathcal H(W)$, then $\tau'$ satisfies Condition~(3) on
$\mathcal H(W')$.
\end{proposition}

\begin{proof}
Inverting the parameters of $\bu$ inverts their $f$-th powers
$u_0,\dots,u_{d-1}$, so the involution $\vee$ of $\mathcal H(W)$ is
compatible with the specialisation and induces the involution $\vee$ of
$\mathcal H(W')$. Specialising Condition~(3) for $\tau$ and restricting
it to elements $\mathbf b\in B(W')\subseteq B(W)$ gives
$\tau'(\mathbf b^{-1})^\vee\,\tau'(\boldsymbol\pi)
=\tau'(\mathbf b\boldsymbol\pi)$, which is Condition~(3) for $\tau'$.
\end{proof}

Thus it suffices to verify Condition~(3) for the groups $G(2n,2,2)$,
$n\ge2$.

\subsection{Palindromicity of the Schur elements of $G(2n,2,2)$}
Let $W=G(2n,2,2)$ with $n \geq 2$. 
We have that $\tau$ is a linear combination of the irreducible characters of $\mathcal H(W)$;  the coefficient of each character $\chi \in {\rm Irr}(\cH(W))$ is equal to the inverse of a Laurent polynomial $s_\chi$, which is called the \emph{Schur element} of $\chi$.
The algebra $\mathcal H(W)$
has $4n$ irreducible characters $\chi_{ijk}$ of degree $1$, with
$\bs \mapsto x_i$, $\bt \mapsto y_j$, $\bu \mapsto z_k$, for $i,j\in\{0,1\}$ and
$0\le k<n$, and $n^2-n$ irreducible characters $\chi^1_{kl}$,
$\chi^2_{kl}$ of degree $2$, for $0\le k<l<n$, on which $\bs$, $\bt$ and
$\bu$ have eigenvalues $\{x_0,x_1\}$, $\{y_0,y_1\}$ and $\{z_k,z_l\}$.
By \cite[Th\'eor\`eme~3.11]{Ma2} and \cite[Corollary 4.9]{MM98}, the Schur
elements of $\mathcal H(W)$ with respect to $\tau$ are
\begin{align*}
s_{\chi_{ijk}}&=\Phi_1(x_ix_{1-i}^{-1})\,\Phi_1(y_jy_{1-j}^{-1})
\prod_{l\neq k}\Phi_1(z_kz_l^{-1})\,
\Phi_1(x_ix_{1-i}^{-1}y_jy_{1-j}^{-1}z_kz_l^{-1}),\\
s_{\chi^{1,2}_{kl}}&=-2\prod_{m\neq k,l}\Phi_1(z_kz_m^{-1})\,
\Phi_1(z_lz_m^{-1})
\prod_{i=0}^1\Phi_1(X_iX_{1-i}^{-1}Y_iY_{1-i}^{-1}Z_kZ_l^{-1})\,
\Phi_1(X_iX_{1-i}^{-1}Y_{1-i}Y_i^{-1}Z_lZ_k^{-1}),
\end{align*}
where $\Phi_1(q)=q-1$, $X_i^2=x_i$, $Y_j^2=y_j$ and $Z_k^2=z_k$. 

For $\chi\in\mathrm{Irr}(\mathcal H(W))$, let $\omega_\chi$ denote its
central character. By \cite[p.~169]{BMM}, the central characters of $\boldsymbol\pi$ are
\[
\omega_{\chi_{ijk}}(\boldsymbol\pi)=(x_iy_jz_k)^{2n}
\qquad\text{and}\qquad
\omega_{\chi^{1,2}_{kl}}(\boldsymbol\pi)=(x_0x_1y_0y_1z_kz_l)^n.
\]
\begin{proposition}\label{prop:palindromic}
For every $\chi\in\mathrm{Irr}(\mathcal H(W))$,
\[
\omega_\chi(\boldsymbol\pi)\,s_\chi^\vee
=(x_0x_1y_0y_1)^n(z_0z_1\cdots z_{n-1})^2\,s_\chi,
\] where $\vee$ inverts the parameters $x_i$, $y_j$, $z_k$ (and also
$X_i$, $Y_j$, $Z_k$).
\end{proposition}

\begin{proof}
Since $\Phi_1(q^{-1})=-q^{-1}\Phi_1(q)$, inverting the parameters
multiplies $s_\chi$ by $\prod(-q^{-1})$, where $q$ runs over the
arguments of the factors $\Phi_1(q)$ of $s_\chi$. There are $2n$ such
factors in each case, so the signs cancel. For $\chi=\chi_{ijk}$, the
product of the inverses of the arguments is
\[
\frac{x_{1-i}}{x_i}\cdot\frac{y_{1-j}}{y_j}\cdot
\prod_{l\neq k}\Bigl(\frac{z_l}{z_k}\cdot
\frac{x_{1-i}y_{1-j}z_l}{x_iy_jz_k}\Bigr)
=\frac{(x_0x_1y_0y_1)^n(z_0\cdots z_{n-1})^2}{(x_iy_jz_k)^{2n}}.
\]
For $\chi=\chi^{1,2}_{kl}$, the factors involving only the $z$'s
contribute $(z_0\cdots z_{n-1})^2/(z_kz_l)^n$. For each $i$, the
arguments of the two remaining factors multiply to
$X_i^2X_{1-i}^{-2}=x_i/x_{1-i}$, and the product over $i=0,1$ is $1$. In
both cases the claim follows from the values of
$\omega_\chi(\boldsymbol\pi)$ above.
\end{proof}

\subsection{Conclusion} 
Now, the map $\vee$ appearing in Proposition \ref{prop:palindromic} is the same as the one appearing in Conjecture \ref{bmmconj}. Moreover, \cite[Lemma~2.7]{BMM} states that Condition (3) of Conjecture \ref{bmmconj} is equivalent to the quotient $\omega_\chi(\boldsymbol\pi)\,s_\chi^\vee/s_\chi$ being constant for all $\chi$ (and this constant is then equal to $\tau(\boldsymbol\pi)$, see the proof of \cite[Lemma 2.7]{BMM}).

\begin{theorem}\label{thm:Gde2}
Let $e$ be even and $d \geq 2$. The BMM symmetrising trace conjecture holds for
$G(de,e,2)$.
\end{theorem}

\begin{proof} 
By \cite[Lemma~2.7]{BMM}, Proposition~\ref{prop:palindromic} implies that
the symmetrising trace $\tau$ of \cite{MM98} on $\mathcal H(G(2n,2,2))$ satisfies
Condition~(3). Since $\tau$ also satisfies Condition~(2) by \cite[\S 4]{MM98}, the conjecture holds for $G(2n,2,2)$. For $G(de,e,2)$,
we apply Proposition~\ref{prop:reduction} with $2n=de$.
\end{proof}

Together with the cases treated in \cite{BreMa, MM98, GIM}, as well as the classical case of the real groups, this
proves the BMM symmetrising trace conjecture for all groups of the
infinite series $G(de,e,n)$:

\begin{theorem}
The BMM symmetrising trace conjecture holds for all groups of the infinite series $G(de,e,n)$.
\end{theorem}

\appendix

\section{Computing the Gram determinants}

The exact values of the Gram determinants are not necessary for the proof of the BMM symmetrising trace conjecture for the exceptional complex reflection groups, as it suffices that they are units. Nonetheless, we give them here, since they easily follow from our techniques.

For a basis $\cB$ of $\cH(W)$ consisting of braid lifts, we define its total weight as
$w(\cB)=\sum_{b\in\cB}w(b)=(w_1(\cB),\dots,w_q(\cB))\in\bZ^q$. Let $N$ be the number of elements $w\in W$ with
$w^2=1$.

\begin{lemma}\label{gramdet}
	Let $\cB$ be a section basis of $\cH(W)$ such that $\tau=\lambda_1^{\cB}$ is the canonical symmetrising trace,
	and let $l_i=2w_i(\cB)/e_i$, where $e_i$ is the order of $s_i$. Then
	\[
	\det\gram_{\tau}(\cB,\cB)=(-1)^{\frac{|W|-N}{2}}\prod_{i=1}^q a_{i,0}^{l_i}.
	\]
\end{lemma}

\begin{proof}
	Let $d_{\cB}=\det\gram_{\tau}(\cB,\cB)\in R(W)$. Since $\tau$ is symmetrising, $d_{\cB}$ is a unit of $R(W)$.
	Hence $\mathrm{pr}(d_{\cB})$ is a unit of the polynomial ring $\bar R(W)=\bZ[a_{i,j}\mid j\ge1]$, that is,
	$\mathrm{pr}(d_{\cB})=\pm1$.
	
	Since $1\in\cB$, $\tau$ is homogeneous of degree~$0$ by Lemma~\ref{exp}, so every entry $\tau(bb')$ of the Gram
	matrix is zero or homogeneous of degree $w(b)+w(b')$. By Lemma~\ref{homog}, applied with
	$\alpha_b=\beta_b=w(b)$, $d_{\cB}$ is homogeneous of degree $2w(\cB)\in \bZ^q$, thus $w_i(d_{\cB})=2w_i(\cB)=e_il_i$.
	
	By Lemma~\ref{liftlem}(2), $d_{\cB}=\pm\prod_ia_{i,0}^{n_i}$ for some integers $n_i$. Since
	$w(a_{i,0})=e_i\varepsilon_i$, we have $w_i(d_{\cB})=e_in_i$, so $n_i=l_i$.
	
	For the sign, specialise $a_{i,0}\mapsto1$ and $a_{i,j}\mapsto0$ for $j\ge1$. The monomial
	$\prod_ia_{i,0}^{l_i}$ maps to $1$. Since $\cB$ is a section basis, it maps to the basis $W$ of $\bZ[W]$, and
	$\tau$ maps to the canonical trace of $\bZ[W]$. Hence $\gram_\tau(\cB,\cB)$ maps to the matrix
	$(\delta_{ww',1})_{w,w'\in W}$, the permutation matrix of $w\mapsto w^{-1}$. This permutation fixes the $N$
	elements with $w^2=1$ and exchanges the other $|W|-N$ in pairs, so its determinant is $(-1)^{(|W|-N)/2}$.
\end{proof}

By Theorem~\ref{thm:Gj} for $j\le29$, and by Proposition~\ref{prop:top} and Theorem~\ref{thm:liftall} for
$j\ge31$, the bases $\cB_j$ (for $j\leq 29$ the $z^k$-bases and for $j\ge31$, the top level of the tower of
		\S\ref{towers}) satisfy the hypotheses of Lemma~\ref{gramdet}. Our programs compute the weights
$w(\cB_j)$ from the words of the bases, and $N$ from the group. Table~\ref{tab:gramdets} lists the
resulting determinants, where the column $S$ lists the generators of $W=G_j$ of Figure~\ref{bfs}, grouped into
		conjugacy classes.

\begin{table}[H]
	\centering
	\footnotesize
	\renewcommand{\arraystretch}{1.3}
	\begin{tabular}{lll|lll}
		$W$ & $S$ & $\det\mathrm\gram_{\tau}(\cB_j,\cB_j)$ & $W$ & $S$ & $\det\gram_{\tau}(\cB_j,\cB_j)$\\
		\hline
		$G_{4}$ & $\{s_1,s_2\}$ & $-a_{s_1,0}^{96}$ & $G_{18}$ & $\{s_1\}\cup\{s_2\}$ & $-a_{s_1,0}^{31700}a_{s_2,0}^{22230}$\\
		$G_{5}$ & $\{s_1\}\cup\{s_2\}$ & $-a_{s_1,0}^{324}a_{s_2,0}^{276}$ & $G_{19}$ & $\{s_1\}\cup\{s_2\}\cup\{s_3\}$ & $+a_{s_1,0}^{97600}a_{s_2,0}^{63600}a_{s_3,0}^{41900}$\\
		$G_{6}$ & $\{s_1\}\cup\{s_2\}$ & $+a_{s_1,0}^{264}a_{s_2,0}^{192}$ & $G_{20}$ & $\{s_1,s_2\}$ & $-a_{s_1,0}^{6912}$\\
		$G_{7}$ & $\{s_1\}\cup\{s_2\}\cup\{s_3\}$ & $+a_{s_1,0}^{936}a_{s_2,0}^{672}a_{s_3,0}^{528}$ & $G_{21}$ & $\{s_1\}\cup\{s_2\}$ & $+a_{s_1,0}^{14040}a_{s_2,0}^{21132}$\\
		$G_{8}$ & $\{s_1,s_2\}$ & $+a_{s_1,0}^{600}$ & $G_{22}$ & $\{s_1,s_2,s_3\}$ & $+a_{s_1,0}^{6288}$\\
		$G_{9}$ & $\{s_1\}\cup\{s_2\}$ & $+a_{s_1,0}^{2240}a_{s_2,0}^{1248}$ & $G_{24}$ & $\{s_1,s_2,s_3\}$ & $+a_{s_1,0}^{5588}$\\
		$G_{10}$ & $\{s_1\}\cup\{s_2\}$ & $+a_{s_1,0}^{2496}a_{s_2,0}^{1728}$ & $G_{25}$ & $\{s_1,s_2,s_3\}$ & $-a_{s_1,0}^{7698}$\\
		$G_{11}$ & $\{s_1\}\cup\{s_2\}\cup\{s_3\}$ & $+a_{s_1,0}^{7296}a_{s_2,0}^{4416}a_{s_3,0}^{3936}$ & $G_{26}$ & $\{s_1\}\cup\{s_2,s_3\}$ & $+a_{s_1,0}^{11016}a_{s_2,0}^{17412}$\\
		$G_{12}$ & $\{s_1,s_2,s_3\}$ & $-a_{s_1,0}^{412}$ & $G_{27}$ & $\{s_1,s_2,s_3\}$ & $+a_{s_1,0}^{105948}$\\
		$G_{13}$ & $\{s_1\}\cup\{s_2,s_3\}$ & $+a_{s_1,0}^{528}a_{s_2,0}^{1024}$ & $G_{29}$ & $\{s_1,\dots,s_4\}$ & $+a_{s_1,0}^{306048}$\\
		$G_{14}$ & $\{s_1\}\cup\{s_2\}$ & $-a_{s_1,0}^{1656}a_{s_2,0}^{1128}$ & $G_{31}$ & $\{s_1,\dots,s_5\}$ & $+a_{s_1,0}^{636096}$\\
		$G_{15}$ & $\{s_1\}\cup\{s_2\}\cup\{s_3\}$ & $+a_{s_1,0}^{3384}a_{s_2,0}^{2376}a_{s_3,0}^{1764}$ & $G_{32}$ & $\{s_1,\dots,s_4\}$ & $+a_{s_1,0}^{1922400}$\\
		$G_{16}$ & $\{s_1,s_2\}$ & $-a_{s_1,0}^{7610}$ & $G_{33}$ & $\{s_1,\dots,s_5\}$ & $+a_{s_1,0}^{859440}$\\
		$G_{17}$ & $\{s_1\}\cup\{s_2\}$ & $+a_{s_1,0}^{32050}a_{s_2,0}^{14740}$ & $G_{34}$ & $\{s_1,\dots,s_6\}$ & $+a_{s_1,0}^{1440348480}$\\
	\end{tabular}
	\caption{Gram determinants of the bases $\cB_j$.} 
	\label{tab:gramdets}
\end{table}

\end{document}